\documentclass[a4paper,11pt]{article}
\usepackage[latin1]{inputenc}
\usepackage[OT1,T1]{fontenc}
\usepackage{a4wide}
\usepackage{amsmath,amssymb}
\usepackage{mathrsfs}
\usepackage[ps,all]{xy}
\usepackage[pdfborder={0 0 0},pdfpagemode=None,pdfstartview=FitH]{hyperref}
\renewcommand\labelitemi{$\ \bullet$}
\renewcommand\theenumi{(\roman{enumi})}
\renewcommand\labelenumi{\theenumi}

\renewcommand\enumerate{\list\labelenumi
  {\setlength\leftmargin{0pt}\setlength\labelwidth{0pt}
  \usecounter{enumi}\def\makelabel##1{\kern\labelsep{##1}}}}
\renewcommand\itemize{\list\labelitemi
  {\setlength\leftmargin{0pt}\setlength\labelwidth{0pt}
  \def\makelabel##1{\kern\labelsep{##1}}}}

\newcommand\altfootnotetext[1]{\insert\footins{\normalfont\footnotesize
  \interlinepenalty\interfootnotelinepenalty\splittopskip\footnotesep
  \splitmaxdepth\dp\strutbox\hsize\columnwidth\rule{0pt}\footnotesep
  \ignorespaces#1\par}}
\newtheorem{theorem}{Theorem}
\newtheorem{lemma}[theorem]{Lemma}
\newtheorem{proposition}[theorem]{Proposition}
\newtheorem{corollary}[theorem]{Corollary}
\newenvironment{proof}{\trivlist
  \item[\hskip\labelsep{\itshape Proof.}]\upshape}{\nobreak\noindent$
  \square$\endtrivlist}
\newenvironment{other}[1]{\refstepcounter{theorem}\trivlist
  \item[\hskip\labelsep{\itshape #1~\arabic{theorem}.}]\upshape}
  {\endtrivlist\bigbreak}
\newenvironment{other*}[1]{\trivlist
  \item[\hskip\labelsep{\itshape #1.}]\upshape}
  {\endtrivlist\bigbreak}
\SelectTips{cm}{10}
\UseTips
\newcommand\SL{\mathord{\mathbf{SL}}}
\newcommand\GL{\mathord{\mathbf{GL}}}
\DeclareMathOperator\im{im}
\DeclareMathOperator\wt{wt}
\DeclareMathOperator\height{ht}

\DeclareMathOperator\End{End}
\DeclareMathOperator\val{val}
\DeclareMathOperator\Spec{Spec}
\DeclareMathOperator\Irr{Irr}
\DeclareMathOperator\Stab{Stab}
\DeclareMathOperator\proj{proj}
\newcommand\high{{\mathord{\mathrm{high}}}}
\newcommand\low{{\mathord{\mathrm{low}}}}
\newcommand\trop{{\mathord{\mathrm{trop}}}}
\newcommand\LS{{\mathord{\mathrm{LS}}}}
\newcommand\dom{{\mathord{\mathrm{dom}}}}
\newcommand\fund{{\mathord{\mathrm{fund}}}}
\newcommand\aff{{\mathord{\mathrm{aff}}}}

\begin{document}
\title{On Mirkovi\'c-Vilonen cycles and crystal combinatorics}
\author{Pierre Baumann and St\'ephane Gaussent\thanks{Both authors
are members of the European Research Training Network ``LieGrits'',
contract no.~MRTN-CT~2003-505078.}}
\date{}
\maketitle
\altfootnotetext{MSC: Primary 20G05, Secondary 05E15 14M15 17B10 22E67.}
\altfootnotetext{Keywords: affine Grassmannian, Mirkovi\'c-Vilonen cycle,
crystal.}
\begin{abstract}
Let $G$ be a complex connected reductive group and let $G^\vee$ be its
Langlands dual. Let us choose a triangular decomposition $\mathfrak
n^{-,\vee} \oplus\mathfrak h^\vee\oplus\mathfrak n^{+,\vee}$ of the Lie
algebra of $G^\vee$. Braverman, Finkelberg and Gaitsgory show that the
set of all Mirkovi\'c-Vilonen cycles in the affine Grassmannian
$G\bigl(\mathbb C((t))\bigr)/G\bigl(\mathbb C[[t]]\bigr)$ is a crystal
isomorphic to the crystal of the canonical basis of $U(\mathfrak
n^{+,\vee})$. Starting from the string parameter of an element of the
canonical basis, we give an explicit description of a dense subset of
the associated MV cycle. As a corollary, we show that the varieties
involved in Lusztig's algebraic-geometric parametrization of the
canonical basis are closely related to MV cycles. In addition, we prove
that the bijection between LS paths and MV cycles constructed by Gaussent
and Littelmann is an isomorphism of crystals.
\end{abstract}

\section{Introduction}
\label{se:Intro}
Let $G$ be a complex connected reductive group, $G^\vee$ be its
Langlands dual, and $\mathscr G$ be its affine Grassmannian. The
geometric Satake correspondence of Lusztig~\cite{Lusztig83},
Beilinson and Drinfeld~\cite{BeilinsonDrinfeld??} and Ginzburg
\cite{Ginzburg95} relates rational representations of $G^\vee$ to
the geometry of $\mathscr G$. More precisely, let us fix a pair of
opposite Borel subgroups in $G$, to enable us to speak of weights
and dominance. Each dominant weight $\lambda$ for $G^\vee$
determines a $G(\mathbb C[[t]])$-orbit $\mathscr G_\lambda$ in
$\mathscr G$. Then the geometric Satake correspondence identifies
the underlying space of the irreducible rational $G^\vee$-module
$L(\lambda)$ with highest weight~$\lambda$ with the intersection
cohomology of $\overline{\mathscr G_\lambda}$.

In~\cite{MirkovicVilonen04}, Mirkovi\'c and Vilonen present a proof of
the geometric Satake correspondence valid in any characteristic. Their
main tool is a class $\mathscr Z(\lambda)$ of subvarieties of
$\overline{\mathscr G_\lambda}$, the so-called MV cycles, which affords
a basis of the intersection cohomology of $\overline{\mathscr G_\lambda}$.
It is tempting to try to compare this construction with standard bases
in $L(\lambda)$, for instance with the canonical basis of Lusztig
\cite{Lusztig90} (also known as the global crystal basis of Kashiwara
\cite{Kashiwara91}).

Several works achieve such a comparison on a combinatorial level. More
precisely, let us recall that the combinatorial object that indexes
naturally the canonical basis of $L(\lambda)$ is the crystal
$\mathbf B(\lambda)$. In~\cite{BravermanGaitsgory01}, Braverman and
Gaitsgory endow the set $\mathscr Z(\lambda)$ with the structure of a
crystal and show the existence of an isomorphism of crystals $\Xi(\lambda):
\mathbf B(\lambda)\stackrel\simeq\longrightarrow\mathscr Z(\lambda)$.
In \cite{GaussentLittelmann05}, Gaussent and Littelmann introduce a set
$\Gamma^+_\LS(\gamma_\lambda)$ of ``LS~galleries''. They endow it with
the structure of a crystal and they associate an MV cycle $Z(\delta)
\in\mathscr Z(\lambda)$ to each LS~gallery $\delta\in\Gamma^+_\LS
(\gamma_\lambda)$. Finally they show the existence of an isomorphism of
crystals $\chi:\mathbf B(\lambda)\stackrel\simeq\longrightarrow
\Gamma^+_\LS(\gamma_\lambda)$ and they prove that the map
$Z:\Gamma^+_\LS(\gamma_\lambda)\to\mathscr Z(\lambda)$ is a bijection.
One of the results of the present paper (Theorem~\ref{th:CompLSBG})
says that Gaussent and Littelmann's map $Z$ is the composition
$\Xi(\lambda)\circ\chi^{-1}$; in particular $Z$ is an isomorphism of
crystals.

Let $\Lambda$ be the lattice of weights of $G^\vee$, let
$\mathfrak n^{-,\vee}\oplus\mathfrak h\oplus\mathfrak n^{+,\vee}$ be
the triangular decomposition of the Lie algebra of $G^\vee$ afforded by
the pinning of $G$, and let $\mathbf B(-\infty)$ be the crystal of
the canonical basis of $U(\mathfrak n^{+,\vee})$. Then for each dominant
weight $\lambda$, the crystal $\mathbf B(\lambda)$ can be embedded into
a shifted version $\mathbf T_{w_0\lambda}\otimes\mathbf B(-\infty)$
of $\mathbf B(-\infty)$, where $w_0\lambda$ is the smallest weight
of $\mathbf B(\lambda)$. It is thus natural to consider a big crystal
$\widetilde{\mathbf B(-\infty)}=\bigoplus_{\lambda\in\Lambda}
\mathbf T_\lambda\otimes\mathbf B(-\infty)$ in order to deal with
all the $\mathbf B(\lambda)$ simultaneously. The isomorphisms
$\Xi(\lambda):\mathbf B(\lambda)\stackrel\simeq\longrightarrow
\mathscr Z(\lambda)$ then assemble in a big bijection
$\Xi:\widetilde{\mathbf B(-\infty)}\stackrel\simeq\longrightarrow
\mathscr Z$. The set $\mathscr Z$ here collects subvarieties of
$\mathscr G$ that have been introduced by Anderson in \cite{Anderson03}.
These varieties are a slight generalization of the usual MV cycles;
indeed $\mathscr Z\supseteq\mathscr Z(\lambda)$ for each dominant weight
$\lambda$. Kamnitzer~\cite{Kamnitzer05a} calls the elements of
$\mathscr Z$ ``stable MV cycles'', but we will simply call them MV
cycles. The existence of $\Xi$ and of a crystal structure on $\mathscr Z$,
and the fact that $\Xi$ is an isomorphism of crystals are due to Braverman,
Finkelberg and Gaitsgory~\cite{BravFinkGait03}.

The crystal $\mathbf B(-\infty)$ can be parametrized in several
ways. Two families of parametrizations, usually called the Lusztig
parametrizations and the string parametrizations (see
\cite{BerensteinZelevinsky01}), depend on the choice of a reduced
decomposition of the longest element in the Weyl group of $G$; they
establish a bijection between $\mathbf B(-\infty)$ and tuples
of natural integers. On the contrary, Lusztig's algebraic-geometric
parametrization \cite{Lusztig96} is intrinsic and describes
$\mathbf B(-\infty)$ in terms of closed subvarieties in
$U^-\bigl(\mathbb C[[t]]\bigr)$, where $U^-$ is the unipotent
radical of the negative Borel subgroup of~$G$.

One of the main results of the present paper is Theorem~\ref{th:ParamOpen},
which describes very explicitly the MV cycle $\Xi(t_0\otimes b)$ starting
from the string parameter of $b\in\mathbf B(-\infty)$. In the course of
his work on MV polytopes \cite{Kamnitzer05a}, Kamnitzer obtains a similar
result, this time starting from the Lusztig parameter of $b$. Though both
results are related (see Section~\ref{ss:LinkKamnCons}), our approach is
foreign to Kamnitzer's methods. Our main ingredient indeed is a concrete
algebraic formula for Braverman, Finkelberg and Gaitsgory's crystal
operations on $\mathscr Z$ that translates the original geometric
definition (Proposition~\ref{pr:CrystOpAlg}). Moreover, our result
implies that Lusztig's algebraic-geometric parametrization is closely
related to MV cycles (Proposition~\ref{pr:CompLuszBFG}).

The paper consists of four sections (plus the introduction). Section
\ref{se:Prelims} fixes some notation and gathers facts and terminology
from the theory of crystals bases. Section~\ref{se:AffGrass} recalls
several standard constructions in the affine Grassmannian and presents
the known results concerning MV cycles. Section~\ref{se:CrystStruct}
defines Braverman, Finkelberg and Gaitsgory's crystal operations on
$\mathscr Z$ and presents our results concerning string parametrizations.
Section~\ref{se:BFG&GalMod} establishes that Gaussent and Littelmann's
bijection $Z:\Gamma^+_\LS(\gamma_\lambda)\to\mathscr Z(\lambda)$ is a
crystal isomorphism. Each section opens with a short summary which gives
a more detailed account of its contents.

We wish to thank M.~Ehrig, J.~Kamnitzer, P.~Littelmann, I.~Mirkovi\'c,
S.~Morier-Genoud and G.~Rousseau for fruitful conversations, vital
information and/or useful indications. We are also grateful to the
referee for his attentive reading and his skilful suggestions.

\section{Preliminaries}
\label{se:Prelims}
The task devoted to Section~\ref{ss:NotPinGps} is to fix the notation
concerning the pinned group $G$. In Section~\ref{ss:Crystals}, we fix
the notation concerning crystal bases for $G^\vee$-modules.

\subsection{Notations for pinned groups}
\label{ss:NotPinGps}
In the entire paper, $G$ will be a complex connected reductive algebraic
group. We assume that a Borel subgroup $B^+$ and a maximal torus
$T\subseteq B^+$ are fixed. We let $B^-$ be the opposite Borel subgroup
to $B^+$ relatively to $T$. We denote the unipotent radical of $B^\pm$
by $U^\pm$.

We denote the character group of $T$ by $X=X^*(T)$; we denote the lattice
of all one-parameter subgroups of $T$ by $\Lambda=X_*(T)$. A point
$\lambda\in\Lambda$ is a morphism of algebraic groups $\mathbb C^\times
\to T,\ a\mapsto a^\lambda$. We denote the root system and the coroot
system of $(G,T)$ by $\Phi$ and $\Phi^\vee=\{\alpha^\vee\mid
\alpha\in\Phi\}$, respectively. The datum of $B^+$ splits $\Phi$ into the
subset $\Phi_+$ of positive roots and the subset $\Phi_-$ of negative
roots. We set $\Phi_+^\vee=\{\alpha^\vee\mid\alpha\in\Phi_+\}$. We denote
by $X_{++}=\{\eta\in X\mid\forall\alpha^\vee\in\Phi_+^\vee,\;\langle\eta,
\alpha^\vee\rangle\geqslant0\}$ and $\Lambda_{++}=\{\lambda\in\Lambda\mid
\forall\alpha\in\Phi_+,\;\langle\alpha,\lambda\rangle\geqslant0\}$ the
cones of dominant weights and coweights. We index the simple roots as
$(\alpha_i)_{i\in I}$. The coroot lattice is the subgroup $\mathbb
Z\Phi^\vee$ generated by the coroots in $\Lambda$; the height of an
element $\lambda=\sum_{i\in I}n_i\alpha_i^\vee$ in $\mathbb Z\Phi^\vee$
is defined as $\height(\lambda)=\sum_{i\in I}n_i$. The dominance order on
$X$ is the partial order $\leqslant$ defined by
$$\eta\geqslant\theta\Longleftrightarrow\eta-\theta\in\mathbb N
\Phi_+.$$
The dominance order on $\Lambda$ is the partial order $\leqslant$
defined by
$$\lambda\geqslant\mu\Longleftrightarrow\lambda-\mu\in\mathbb N
\Phi_+^\vee.$$

For each simple root $\alpha_i$, we choose a non-trivial additive
subgroup $x_i$ of $U^+$ such that $a^\lambda x_i(b)a^{-\lambda}=
x_i\bigl(a^{\langle\alpha_i,\lambda\rangle}b\bigr)$ holds for all
$\lambda\in\Lambda$, $a\in\mathbb C^\times$, $b\in\mathbb C$.
Then there is a unique morphism $\varphi_i:\SL_2\to G$ such that
$$\varphi_i\begin{pmatrix}1&b\\0&1\end{pmatrix}=x_i(b)\quad\text{and}
\quad\varphi_i\begin{pmatrix}a&0\\0&a^{-1}\end{pmatrix}=a^{\alpha_i^\vee}$$
for all $a\in\mathbb C^\times$, $b\in\mathbb C$. We set
$$y_i(b)=\varphi_i\begin{pmatrix}1&0\\b&1\end{pmatrix}\quad\text{and}
\quad\overline{s_i}=\varphi_i\begin{pmatrix}0&1\\-1&0\end{pmatrix}.$$

Let $N_G(T)$ be the normalizer of $T$ in $G$ and let $W=N_G(T)/T$ be
the Weyl group of $(G,T)$. Each element $\overline{s_i}$ normalizes $T$;
its class $s_i$ modulo $T$ is called a simple reflection. Endowed with
the set of simple reflections, the Weyl group becomes a Coxeter system.
Since the elements $\overline{s_i}$ satisfy the braid relations, we may
lift each element $w\in W$ to an element $\overline w\in G$ so that
$\overline w=\overline{s_{i_1}}\cdots\overline{s_{i_l}}$ for any
reduced decomposition $s_{i_1}\cdots s_{i_l}$ of $w$. For any two
elements $w$ and $w'$ in $W$, there exists an element $\lambda\in
\mathbb Z\Phi^\vee$ such that $\overline{ww'}=(-1)^\lambda\,
\overline{w^{\vphantom\prime}}\,\overline{w'}$. We denote the longest
element of $W$ by $w_0$.

Let $\alpha$ be a positive root. We make the choice of a simple root
$\alpha_i$ and of an element $w\in W$ such that $\alpha=w\alpha_i$.
Then we define the one-parameter additive subgroups
\begin{equation}
\label{eq:EqCommut1}
x_\alpha:b\mapsto\overline w\,x_i(b)\,\overline w^{-1}\quad\text{and}
\quad x_{-\alpha}:b\mapsto\overline w\,y_i(b)\,\overline w^{-1}
\end{equation}
and the element $\overline{s_\alpha}=\overline w\,\overline{s_i}\,
\overline w^{-1}$.

Products in $G$ may then be computed using several commutation rules:
\begin{itemize}
\item For any $\lambda\in\Lambda$, any root $\alpha$, any $a\in\mathbb
C^\times$ and any $b\in\mathbb C$,
\begin{equation}
\label{eq:EqCommut2}
a^\lambda x_\alpha(b)=x_\alpha\bigl(a^{\langle\alpha,\lambda\rangle}b
\bigr)a^\lambda.
\end{equation}
\item For any root $\alpha$ and any $a,b\in\mathbb C$ such that
$1+ab\neq0$,
\begin{equation}
\label{eq:EqCommut3}
x_\alpha(a)x_{-\alpha}(b)=x_{-\alpha}\bigl(b/(1+ab)\bigr)
\bigl(1+ab\bigr)^{\alpha^\vee}
x_\alpha\bigl(a/(1+ab)\bigr).
\end{equation}
\item For any positive root $\alpha$ and any $a\in\mathbb C^\times$,
\begin{equation}
\label{eq:EqCommut4}
x_\alpha(a)\,x_{-\alpha}(-a^{-1})\,x_\alpha(a)=x_{-\alpha}(-a^{-1})
\,x_\alpha(a)\,x_{-\alpha}(-a^{-1})=a^{\alpha^\vee}\,\overline{s_\alpha}=
\overline{s_\alpha}\;a^{-\alpha^\vee}.
\end{equation}
\item (Chevalley's commutator formula) If $\alpha$ and $\beta$ are two
linearly independent roots, then there are numbers $C_{i,j,\alpha,\beta}
\in\{\pm1,\pm2,\pm3\}$ such that
\begin{equation}
\label{eq:EqCommut5}
x_\beta(b)^{-1}x_\alpha(a)^{-1}x_\beta(b)\,x_\alpha(a)=\prod_{i,j>0}
x_{i\alpha+j\beta}\bigl(C_{i,j,\alpha,\beta}(-a)^ib^j\bigr)
\end{equation}
for all $a$ and $b$ in $\mathbb C$. The product in the right-hand side
is taken over all pairs of positive integers $i,j$ for which
$i\alpha+j\beta$ is a root, in order of increasing $i+j$.
\end{itemize}

\subsection{Crystals}
\label{ss:Crystals}
Let $G^\vee$ be the Langlands dual of $G$. This connected reductive group
is equipped with a Borel subgroup $B^{+,\vee}$ and a maximal torus
$T^\vee\subseteq B^{+,\vee}$ so that $\Lambda$ is the weight lattice of
$T^\vee$ and $\Phi^\vee$ is the root system of $(G^\vee,T^\vee)$, the
set of positive roots being $\Phi_+^\vee$. The Lie algebra $\mathfrak
g^\vee$ of $G^\vee$ has a triangular decomposition $\mathfrak g^\vee=
\mathfrak n^{-,\vee}\oplus\mathfrak h^\vee\oplus\mathfrak n^{+,\vee}$.

A crystal for $G^\vee$ (in the sense of Kashiwara~\cite{Kashiwara95}) is
a set $\mathbf B$ endowed with maps
$$\tilde e_i,\tilde f_i:\mathbf B\to\mathbf B\sqcup\{0\},\quad
\varepsilon_i,\varphi_i:\mathbf B\to\mathbb Z\sqcup\{-\infty\},
\quad\text{and}\quad\wt:\mathbf B\to\Lambda,$$
where $0$ is a ghost element added to $\mathbf B$ in order that
$\tilde e_i$ and $\tilde f_i$ may be everywhere defined. These
maps are required to satisfy certain axioms, which the reader may
find in Section~7.2 of~\cite{Kashiwara95}. The map $\wt$ is called
the weight.

A morphism from a crystal $\mathbf B$ to a crystal $\mathbf B'$ is
a map $\psi:\mathbf B\sqcup\{0\}\to\mathbf B'\sqcup\{0\}$ satisfying
$\psi(0)=0$ and compatible with the structure maps $\tilde e_i$,
$\tilde f_i$, $\varepsilon_i$, $\varphi_i$ and $\wt$. The conditions
are written in full detail in~\cite{Kashiwara95}.

Given a crystal $\mathbf B$, one defines a crystal $\mathbf B^\vee$
whose elements are written $b^\vee$, where $b\in\mathbf B$, and whose
structure maps are given by
\begin{xalignat*}2
\varepsilon_i(b^\vee)&=\varphi_i(b),&
\tilde e_i(b^\vee)&=(\tilde f_ib)^\vee,\\
\varphi_i(b^\vee)&=\varepsilon_i(b),&
\tilde f_i(b^\vee)&=(\tilde e_ib)^\vee,\\[2pt]
\wt(b^\vee)&=-\wt(b),&&
\end{xalignat*}
where one sets $0^\vee=0$. The correspondence $\mathbf B\rightsquigarrow
\mathbf B^\vee$ is a covariant functor. (Caution: Usually in this paper,
the symbol $\vee$ is used to adorn coroots or objects related to the
Langlands dual. Here and in Section~\ref{ss:LuszAlgGeomParam} however,
it will also be used to denote contragredient duality for crystals.)

The most important crystals for our work are the crystal $\mathbf
B(\infty)$ of the canonical basis of $U(\mathfrak n^{-,\vee})$ and the
crystal $\mathbf B(-\infty)$ of the canonical basis of
$U(\mathfrak n^{+,\vee})$. The crystal $\mathbf B(\infty)$ is a highest
weight crystal; this means that it has an element annihilated by all
operators $\tilde e_i$ and from which any other element of $\mathbf
B(\infty)$ can be obtained by applying the operators $\tilde f_i$. This
element is unique and its weight is $0$; we denote it by $1$. Likewise
the crystal $\mathbf B(-\infty)$ is a lowest weight crystal; its lowest
weight element has weight $0$ and is also denoted by $1$.

The antiautomorphism of the algebra $U(\mathfrak n^{-,\vee})$ that fixes
the Chevalley generators leaves stable its canonical basis; it therefore
induces an involution $b\mapsto b^*$ of the set $\mathbf B(\infty)$.
This involution $*$ preserves the weight. The operators $\tilde f_i$ and
$b\mapsto(\tilde f_ib^*)^*$ correspond roughly to the left and right
multiplication in $U(\mathfrak n^{-,\vee})$ by the Chevalley generator
with index $i$ (see Proposition~5.3.1 in~\cite{Kashiwara93a} for a more
precise statement). One can therefore expect that $\tilde f_i$ and
$b\mapsto(\tilde f_jb^*)^*$ commute for all $i,j\in I$. This actually
holds only if $i\neq j$; and when $i=j$, one can analyze precisely the
mutual behavior of these operators. In return, one obtains a
characterization of $\mathbf B(\infty)$ as the unique highest weight
crystal generated by a highest weight element of weight $0$ and endowed
with an involution $*$ with specific properties (see Section~2 in
\cite{Kashiwara93b}, Proposition~3.2.3 in \cite{KashiwaraSaito97}, and
Section~12 in~\cite{BravFinkGait03} for more details).

For any weight $\lambda\in\Lambda$, we consider the crystal $\mathbf
T_\lambda$ with unique element $t_\lambda$, whose structure maps are
given by
$$\wt(t_\lambda)=\lambda,\quad\tilde e_it_\lambda=\tilde f_it_\lambda=
0\quad\text{and}\quad\varepsilon_i(t_\lambda)=\varphi_i(t_\lambda)=-\infty$$
(see Example~7.3 in~\cite{Kashiwara95}). There are two operations $\oplus$
and $\otimes$ on crystals (see Section~7.3 in~\cite{Kashiwara95}). We set
$\widetilde{\mathbf B(-\infty)}=\bigoplus_{\lambda\in\Lambda}\mathbf
T_\lambda\otimes\mathbf B(-\infty)$. Thus for any $b\in\mathbf B(-\infty)$,
any $\lambda\in\Lambda$ and any $i\in I$,
\begin{xalignat*}2
\varepsilon_i(t_\lambda\otimes b)&=\varepsilon_i(b)-
\langle\alpha_i,\lambda\rangle,&
\tilde e_i(t_\lambda\otimes b)&=t_\lambda\otimes\tilde e_i(b),\\
\varphi_i(t_\lambda\otimes b)&=\varphi_i(b),&
\tilde f_i(t_\lambda\otimes b)&=t_\lambda\otimes\tilde f_i(b),\\[2pt]
\wt(t_\lambda\otimes b)&=\wt(b)+\lambda.&&
\end{xalignat*}
We transport the involution $*$ from $\mathbf B(\infty)$ to
$\mathbf B(-\infty)$ by using the isomorphism $\mathbf B(-\infty)\cong
\mathbf B(\infty)^\vee$ and by setting $(b^\vee)^*=(b^*)^\vee$ for each
$b\in\mathbf B(\infty)$. Then we extend it to $\widetilde{\mathbf
B(-\infty)}$ by setting
$$(t_\lambda\otimes b)^*=t_{-\lambda-\wt(b)}\otimes b^*.$$

For $\lambda\in\Lambda$, we denote by $L(\lambda)$ the irreducible
rational representation of $G^\vee$ whose highest weight is the unique
dominant weight in the orbit $W\lambda$. We denote the crystal of the
canonical basis of $L(\lambda)$ by $\mathbf B(\lambda)$. It has a unique
highest weight element $b_\high$ and a unique lowest weight element
$b_\low$, which satisfy $\tilde e_ib_\high=\tilde f_ib_\low=0$ for any
$i\in I$. If $\lambda$ is dominant, there is a unique embedding of
crystals $\kappa_\lambda:\mathbf B(\lambda)\hookrightarrow\mathbf
B(\infty)\otimes\mathbf T_\lambda$; it maps the element $b_\high$ to
$1\otimes t_\lambda$ and its image is
$$\{b\otimes t_\lambda\mid b\in\mathbf B(\infty)\text{ such that }
\forall i\in I,\;\varepsilon_i(b^*)\leqslant\langle\alpha_i,\lambda
\rangle\}$$
(see Proposition~8.2 in~\cite{Kashiwara95}). If $\lambda$ is antidominant,
then the sequence
$$\mathbf B(\lambda)\cong\mathbf B(-\lambda)^\vee\xrightarrow{
(\kappa_{-\lambda})^\vee}\bigl(\mathbf B(\infty)\otimes\mathbf
T_{-\lambda}\bigr)^\vee\cong\mathbf T_\lambda\otimes\mathbf B(-\infty)$$
defines an embedding of crystals $\iota_\lambda:\mathbf B(\lambda)
\hookrightarrow\mathbf T_\lambda\otimes\mathbf B(-\infty)$;
it maps the element $b_\low$ to $t_\lambda\otimes1$ and its image is
$$\{t_{\lambda}\otimes b\mid b\in\mathbf B(-\infty)\text{ such
that }\forall i\in I,\;\varphi_i(b^*)\leqslant-\langle\alpha_i,
\lambda\rangle\}.$$

\section{The affine Grassmannian}
\label{se:AffGrass}
In Section~\ref{ss:AffGrDefs}, we recall the definition of an
affine Grassmannian. In Section~\ref{ss:Orbits}, we present several
properties of orbits in the affine Grassmannian of $G$ under the
action of the groups $G\bigl(\mathbb C[[t]]\bigr)$ and $U^\pm
\bigl(\mathbb C((t))\bigr)$. Section~\ref{ss:MVCycles} recalls the
notion of MV cycle, in the original version of Mirkovi\'c and Vilonen
and in the somewhat generalized version of Anderson. Finally Section
\ref{ss:ParabRet} introduces a map from the affine Grassmannian of $G$
to the affine Grassmannian of a Levi subgroup of $G$.

An easy but possibly new result in this section is Proposition
\ref{pr:IncidPropStrat}~\ref{it:PrIPSc}. Joint with Mirkovi\'c and
Vilonen's work, it implies the expected Proposition~\ref{pr:DimEstimate},
which provides the dimension estimates that Anderson needs for his
generalization of MV cycles.

\subsection{Definitions}
\label{ss:AffGrDefs}
We denote the ring of formal power series by $\mathscr O=\mathbb C[[t]]$
and we denote its field of fractions by $\mathscr K=\mathbb C((t))$. We
denote the valuation of a non-zero Laurent series $f\in\mathscr K^\times$
by $\val(f)$. Given a complex linear algebraic group $H$, we define the
affine Grassmannian of $H$ as the space $\mathscr H=H(\mathscr K)/
H(\mathscr O)$. The class in $\mathscr H$ of an element $h\in H(\mathscr
K)$ will be denoted by $[h]$.

\begin{other*}{Example}
If $H$ is the multiplicative group $\mathbf G_m$, then the valuation map
yields a bijection from $\mathscr H=\mathscr K^\times/\mathscr O^\times$
onto $\mathbb Z$. More generally, if $H$ is a torus, then the map
$\lambda\mapsto[t^\lambda]$ is a bijection from the lattice $X_*(H)$ of
one-parameter subgroups in $H$ onto the affine Grassmannian $\mathscr H$.
\end{other*}

The affine Grassmannian $\mathscr H$ is the set of $\mathbb C$-points of
an ind-scheme defined over $\mathbb C$ (see \cite{BeauvilleLaszlo94} for
$H=\GL_n$ or $\SL_n$ and Chapter~13 of \cite{Kumar02} for $H$ simple).
This means in particular that $\mathscr H$ is the direct limit of a system
$$\mathscr H_0\hookrightarrow\mathscr H_1\hookrightarrow\mathscr H_2
\hookrightarrow\cdots$$
of complex algebraic varieties and of closed embeddings. We endow
$\mathscr H$ with the direct limit of the Zariski topologies on the
varieties $\mathscr H_n$. A noetherian subspace $Z$ of $\mathscr H$
thus enjoys the specific topological properties of a subset of a complex
algebraic variety; for instance if $Z$ is locally closed, then
$\dim Z=\dim\overline Z$.

The affine Grassmannian of the groups $G$ and $T$ considered in
Section~\ref{ss:NotPinGps} will be denoted by $\mathscr G$ and $\mathscr
T$, respectively. The inclusion $T\subseteq G$ gives rise to a closed
embedding $\mathscr T\hookrightarrow\mathscr G$.

\subsection{Orbits}
\label{ss:Orbits}
We first look at the action of the group $G(\mathscr O)$ on $\mathscr G$
by left multiplication. The orbit $G(\mathscr O)[t^\lambda]$ depends only
on the $W$-orbit of $\lambda$ in $\Lambda$, and the Cartan decomposition
of $G(\mathscr K)$ says that
$$\mathscr G=\bigsqcup_{W\lambda\in\Lambda/W}G(\mathscr O)[t^\lambda].$$
For each coweight $\lambda\in\Lambda$, the orbit $\mathscr G_\lambda=
G(\mathscr O)[t^\lambda]$ is a noetherian subspace of $\mathscr G$.
If $\lambda$ is dominant, then the dimension of $\mathscr G_\lambda$ is
$\height(\lambda-w_0\lambda)$ and its closure is
\begin{equation}
\overline{\mathscr G_\lambda}=\bigsqcup_{\substack{\mu\in\Lambda_{++}\\[2pt]
\lambda\geqslant\mu}}\mathscr G_\mu.
\label{eq:ClosureOrbit}
\end{equation}

From this, one can quickly deduce that it is often possible to truncate
power series when dealing with the action of $G(\mathscr O)$ on
$\mathscr G$. Given an positive integer $s$, let $G_{(s)}$ denote the
$s$-th congruence subgroup of $G(\mathscr O)$, that is, the kernel of the
reduction map $G(\mathscr O)\to G(\mathscr O/t^s\mathscr O)$.
\begin{proposition}
\label{pr:CanTruncate}
For each noetherian subset $Z$ of $\mathscr G$, there exists a level
$s$ such that $G_{(s)}$ fixes $Z$ pointwise.
\end{proposition}
\begin{proof}
Let $\bigl(\Lambda_{++}^{(n)}\bigr)_{n\in\mathbb N}$ be an
increasing sequence of finite subsets of $\Lambda_{++}$ such that
$$\bigl\{\nu\in\Lambda_{++}\bigm|\nu\leqslant\mu\bigr\}\subseteq
\Lambda_{++}^{(n)}\quad\text{for each $\mu\in\Lambda_{++}^{(n)}$\quad
and that}\quad\bigcup_{n\in\mathbb N}\Lambda_{++}^{(n)}=\Lambda_{++}.$$
Set $\mathscr G_n=\bigsqcup_{\mu\in\Lambda_{++}^{(n)}}\mathscr G_\mu$.
The Cartan decomposition shows that $(\mathscr G_n)_{n\geqslant0}$ is
an increasing and exhaustive filtration of $\mathscr G$, and Equation
(\ref{eq:ClosureOrbit}) shows that each $\mathscr G_n$ is closed.
Therefore each noetherian subset $Z$ of $\mathscr G$ is contained in
$\mathscr G_n$ for $n$ sufficiently large. To prove the proposition,
it is thus enough to show that for each integer $n$, there is an
$s\geqslant1$ such that $G_{(s)}$ fixes $\mathscr G_n$ pointwise.

Let $\lambda\in\Lambda$, and choose $s\geqslant1$ larger
than $\langle\alpha,\lambda\rangle$ for all $\alpha\in\Phi$. Using that
$G_{(s)}$ is generated by elements $(1+t^sp)^\lambda$ and $x_\alpha(t^sp)$
with $\lambda\in\Lambda$, $\alpha\in\Phi$ and $p\in\mathscr O$, one
readily checks that $G_{(s)}$ fixes the point $[t^\lambda]$. Since
$G_{(s)}$ is normal in $G(\mathscr O)$, it pointwise fixes the orbit
$\mathscr G_\lambda$. The proposition then follows from the fact that
each $\mathscr G_n$ is a finite union of $G(\mathscr O)$-orbits.
\end{proof}

We now look at the action of the unipotent group $U^\pm(\mathscr K)$
on $\mathscr G$. It can be described by the Iwasawa decomposition
$$\mathscr G=\bigsqcup_{\lambda\in\Lambda}U^\pm(\mathscr K)[t^\lambda].$$
We will denote the orbit $U^\pm(\mathscr K)[t^\lambda]$ by $S_\lambda^\pm$.
Proposition~3.1~(a) in~\cite{MirkovicVilonen04} asserts that the closure
of a stratum $S_\lambda^\pm$ is the union
\begin{equation}
\overline{S_\lambda^\pm}=\bigsqcup_{\substack{\mu\in\Lambda\\[2pt]
\pm(\lambda-\mu)\geqslant0}}S_\mu^\pm.
\label{eq:ClosureStratum}
\end{equation}
This equation implies in particular
$$S_\lambda^\pm=\overline{S_\lambda^\pm}\setminus\left(\bigcup_{i\in I}
\overline{S_{\lambda\mp\alpha_i^\vee}^\pm}\right),$$
which shows that each stratum $S_\lambda^\pm$ is locally closed.

As pointed out by Mirkovi\'c and Vilonen (Equation~(3.5) in
\cite{MirkovicVilonen04}), these strata $S_\lambda^\pm$ can be understood
in terms of a Bia\l ynicki-Birula decomposition: indeed the choice of
a dominant and regular coweight $\xi\in\Lambda$ defines an action of
$\mathbb C^\times$ on $\mathscr G$, and
$$S_\lambda^\pm=\{x\in\mathscr G\mid\lim_{\substack{a\to0\\a\in\mathbb
C^\times}}a^{\pm\xi}\cdot x=[t^\lambda]\}$$
for each $\lambda\in\Lambda$. We will generalize this result in
Remark~\ref{rk:RetracCellsBB}. For now, we record the following two
(known and obvious) consequences:
\begin{itemize}
\item The set of points in $\mathscr G$ fixed by the action of $T$
is $\mathscr G^T=\bigl\{[t^\lambda]\bigm|\lambda\in\Lambda\bigr\}$;
in other words, $\mathscr G^T$ is the image of the embedding
$\mathscr T\hookrightarrow\mathscr G$.
\item If $Z$ is a closed and $T$-invariant subset of $\mathscr G$, then
$Z$ meets a stratum $S_\lambda^\pm$ if and only if $[t^\lambda]\in Z$.
\end{itemize}

The following proposition is in essence due to Kamnitzer (see Section~3.3
in~\cite{Kamnitzer05a}).
\begin{proposition}
\label{pr:DefMu}
Let $Z$ be an irreducible and noetherian subset of $\mathscr G$.
\begin{enumerate}
\item\label{it:PrDMa}
The set $\{\lambda\in\Lambda\mid Z\cap S_\lambda^+\neq\varnothing\}$
is finite and has a largest element. Denoting the latter by $\mu_+$, the
intersection $Z\cap S_{\mu_+}^+$ is open and dense in~$Z$.
\item\label{it:PrDMb}
The set $\{\lambda\in\Lambda\mid Z\cap S_\lambda^-\neq\varnothing\}$
is finite and has a smallest element. Denoting the latter by $\mu_-$, the
intersection $Z\cap S_{\mu_-}^-$ is open and dense in~$Z$.
\end{enumerate}
\end{proposition}
Given an irreducible and noetherian subset $Z$ in $\mathscr G$, we
indicate the coweights $\mu_\pm$ exhibited in Proposition~\ref{pr:DefMu}
by the notation $\mu_\pm(Z)$.

\trivlist
\item[\hskip\labelsep{\itshape Proof of Proposition~\ref{pr:DefMu}.}]
\upshape
The Cartan decomposition and the equality $\mathscr G^T=\{[t^\lambda]
\mid\lambda\in\Lambda\}$ imply that the obvious inclusion $(\mathscr
G_\nu)^T\supseteq\{[t^{w\nu}]\mid w\in W\}$ is indeed an equality for
each coweight $\nu\in\Lambda$. Therefore $X^T$ is finite for each subset
$X\subseteq\mathscr G$ that is a finite union of $G(\mathscr O)$-orbits.
This is in particular the case for each of the subsets $\mathscr G_n$
used in the proof of Proposition~\ref{pr:CanTruncate}. Since $\mathscr
G_n$ is moreover closed and $T$-invariant, this means that it meets only
finitely many strata $S_\lambda^+$. Thus a noetherian subset of
$\mathscr G$ meets only finitely many strata $S_\lambda^+$, for
it is contained in $\mathscr G_n$ for $n$ large enough.

Assume now that $Z$ is an irreducible and noetherian subset of $\mathscr
G$. Each intersection $Z\cap S_\lambda^+$ is locally closed in $Z$
and $Z$ is covered by finitely many such intersections, so there
exists a coweight $\mu_+$ for which the intersection $Z\cap S_{\mu_+}^+$
is dense in $Z$. Then $Z\subseteq\overline{S_{\mu_+}^+}$; by Equation
(\ref{eq:ClosureStratum}), this means that $\mu_+$ is the largest
element in $\{\lambda\in\Lambda\mid Z\cap S_\lambda^+\neq\varnothing\}$.
Moreover $Z\cap S_{\mu_+}^+$ is locally closed in $Z$; it is therefore
open in its closure in $Z$, which is $Z$.

The arguments above prove Assertion~\ref{it:PrDMa}. The proof of Assertion
\ref{it:PrDMb} is entirely similar.
\nobreak\noindent$\square$
\endtrivlist

\begin{other}{Examples}
\label{ex:DeterMu}
\begin{enumerate}
\item\label{it:ExDMa}
If $Z$ is an irreducible and noetherian subset of $\mathscr G$,
then $Z\cap S_{\mu_+(Z)}^+\cap S_{\mu_-(Z)}^-$ is dense in $Z$. Thus
$Z$ and $\overline Z$ are contained in $\overline{S_{\mu_+(Z)}^+\cap
S_{\mu_-(Z)}^-}$. One deduces from this the equality
$\mu_\pm(\overline Z)=\mu_\pm(Z)$.
\item\label{it:ExDMb}
For any coweight $\lambda\in\Lambda$, $\mu_+(\mathscr G_\lambda)
=\mu_+\bigl(\overline{\mathscr G_\lambda}\bigr)$ and
$\mu_-(\mathscr G_\lambda)=\mu_-\bigl(\overline{\mathscr
G_\lambda}\bigr)$ are the largest and the smallest element in the
orbit~$W\lambda$, respectively.
\end{enumerate}
\end{other}

We now present a method that allows to find the parameter $\lambda$ of
an orbit $\mathscr G_\lambda$ or $S_\lambda^\pm$ to which a given point
of $\mathscr G$ belongs. Given a $\mathbb C$-vector space $V$, we may
form the $\mathscr K$-vector space $V\otimes_{\mathbb C}\mathscr K$ by
extending the base field and regard $V$ as a subspace of it. In this
situation, we define the valuation $\val(v)$ of a non-zero vector
$v\in V\otimes_{\mathbb C}\mathscr K$ as the largest $n\in\mathbb Z$ such
that $v\in V\otimes t^n\mathscr O$; thus the valuation of a non-zero
element $v\in V$ is zero. We define the valuation $\val(f)$ of a non-zero
endomorphism $f\in\End_{\mathscr K}(V\otimes_{\mathbb C}\mathscr K)$ as
the largest $n\in\mathbb Z$ such that $f(V\otimes_{\mathbb C}\mathscr O)
\subseteq V\otimes t^n\mathscr O$; equivalently, $\val(f)$ is the
valuation of $f$ viewed as an element in
$\End_{\mathbb C}(V)\otimes_{\mathbb C}\mathscr K$.

For each weight $\eta\in X$, we denote by $V(\eta)$ the simple rational
representation of $G$ whose highest weight is the dominant weight in
the orbit $W\eta$, and we choose an extremal weight vector $v_\eta\in
V(\eta)$ of weight $\eta$. The structure map $g\mapsto g_{V(\eta)}$
from $G$ to $\End_{\mathbb C}(V(\eta))$ of this representation extends
to a map from $G(\mathscr K)$ to $\End_{\mathscr K}(V(\eta)
\otimes_{\mathbb C}\mathscr K)$; we denote this latter also by $g\mapsto
g_{V(\eta)}$, or simply by $g\mapsto g\cdot?$ if there is no risk of
confusion.

\begin{proposition}
\label{pr:KamnitzForm}
Let $g\in\mathscr G(\mathscr K)$.
\begin{enumerate}
\item\label{it:PrKFa}
The antidominant coweight $\lambda\in\Lambda$ such that $[g]\in\mathscr
G_\lambda$ is characterized by the equations
$$\forall\eta\in X_{++},\quad\langle\eta,\lambda\rangle=
\val\bigl(g_{V(\eta)}\bigr).$$
\item\label{it:PrKFb}
The coweight $\lambda\in\Lambda$ such that $[g]\in S_\lambda^\pm$ is
characterized by the equations
$$\forall\eta\in X_{++},\quad\pm\langle\eta,\lambda\rangle=
-\val(g^{-1}\cdot v_{\pm\eta}).$$
\end{enumerate}
\end{proposition}
\begin{proof}
Assertion~\ref{it:PrKFb} is due to Kamnitzer (this is Lemma~2.4
in~\cite{Kamnitzer05a}), so we only have to prove Assertion~\ref{it:PrKFa}.
Let $\lambda\in\Lambda$ be antidominant and let $\eta\in X_{++}$.
Then for each weight $\theta$ of $V(\eta)$, the element $t^\lambda$ acts
by $t^{\langle\lambda,\theta\rangle}$ on the $\theta$-weight subspace of
$V(\eta)$, with here $\langle\lambda,\theta\rangle\geqslant\langle\lambda,
\eta\rangle$ since $\theta\leqslant\eta$. It follows that
$\val\bigl((t^\lambda)_{V(\eta)}\bigr)=\langle\lambda,\eta\rangle$.
Thus the proposed formula holds for $g=t^\lambda$. To conclude the
proof, it suffices to observe that $\val\bigl(g_{V(\eta)}\bigr)$ depends
only of the double coset $G(\mathscr O)gG(\mathscr O)$, for the action of
$G(\mathscr O)$ leaves $V(\eta)\otimes_{\mathbb C}\mathscr O$ invariant.
\end{proof}

We end this section with a proposition that provides some information
concerning intersections of orbits. We agree to say that an assertion
$A(\lambda)$ depending on a coweight $\lambda\in\Lambda$ holds when
$\lambda$ is enough antidominant if
$$(\exists N\in\mathbb Z)\quad(\forall\lambda\in\Lambda)\quad
(\forall i\in I,\ \langle\alpha_i,\lambda\rangle\leqslant N)
\Longrightarrow A(\lambda).$$
\begin{proposition}
\label{pr:IncidPropStrat}
\begin{enumerate}
\item\label{it:PrIPSa}
Let $\lambda,\nu\in\Lambda$. If $S_\lambda^+\cap S_\nu^-\neq\varnothing$,
then $\lambda\geqslant\nu$.
\item\label{it:PrIPSb}
Let $\lambda\in\Lambda$. Then $S_\lambda^+\cap S_\lambda^-=
\bigl\{[t^\lambda]\bigr\}$.
\item\label{it:PrIPSc}
Let $\nu\in\Lambda$ such that $\nu\geqslant0$. If $\lambda\in\Lambda$ is
enough antidominant, then $S_{\lambda+\nu}^+\cap S_\lambda^-=
S_{\lambda+\nu}^+\cap\mathscr G_\lambda$.
\end{enumerate}
\end{proposition}
The proof of this proposition requires a lemma.

\begin{lemma}
\label{le:IncidPropStrat}
Let $\nu\in\Lambda$ such that $\nu\geqslant0$. If $\lambda\in\Lambda$ is
enough antidominant, then $S_{\lambda+\nu}^+\cap S_\lambda^-\subseteq
\mathscr G_\lambda$.
\end{lemma}
\begin{proof}
For the whole proof, we fix $\nu\in\Lambda$ such that $\nu\geqslant0$.

For each $\eta\in X_{++}$, we make the following construction. We form
the list $(\theta_1,\theta_2,\ldots,\theta_N)$ of all the weights of
$V(\eta)$, repeated according to their multiplicities and ordered in
such a way that $(\theta_i>\theta_j\Longrightarrow i<j)$ for all indices
$i,j$. Thus $N=\dim V(\eta)$, $\theta_1=\eta>\theta_i$ for all $i>1$,
and $\theta_1+\theta_2+\cdots+\theta_N$ is $W$-invariant hence
orthogonal to $\mathbb Z\Phi^\vee$. We say then that a coweight
$\lambda\in\Lambda$ satisfies Condition~$A_\eta(\lambda)$~if
$$\forall j\in\{1,\ldots,N\},\quad\langle\theta_1-\theta_j,\lambda
\rangle\leqslant\langle\theta_j+\theta_{j+1}+\cdots+\theta_N,\nu\rangle.$$
Certainly Condition~$A_\eta(\lambda)$ holds if $\lambda$ is enough
antidominant.

Now we choose a finite subset $Y\subseteq X_{++}$ that spans the
lattice $X$ up to torsion. To prove the lemma, it is enough to show that
$S_{\lambda+\nu}^+\cap S_\lambda^-\subseteq\mathscr G_\lambda$ for all
antidominant $\lambda$ satisfying Condition~$A_\eta(\lambda)$ for each
$\eta\in Y$.

Suppose that $\lambda$ satisfies these requirements and let $g\in
U^-(\mathscr K)t^\lambda$ be such that $[g]\in S_{\lambda+\nu}^+$.
We use Proposition~\ref{pr:KamnitzForm}~\ref{it:PrKFa} to show that
$[g]\in\mathscr G_\lambda$. Let $\eta\in Y$. Let $(v_1,v_2,\ldots,v_N)$
be a basis of $V(\eta)$ such that for each $i$, $v_i$ is a vector of
weight $\theta_i$. We denote the dual basis in $V(\eta)^*$ by $(v_1^*,
v_2^*,\ldots,v_N^*)$; thus $v_i^*$ is of weight $-\theta_i$. Then
$$\val\bigl(g_{V(\eta)}\bigr)=\min\bigl\{\val(\langle v_j^*,g\cdot
v_i\rangle)\bigm|1\leqslant i,j\leqslant N\bigl\}.$$
The choice $g\in U^-(\mathscr K)t^\lambda$ implies that the matrix of
$g_{V(\eta)}$ in the basis $(v_i)_{1\leqslant i\leqslant N}$ is lower
triangular, with diagonal entries $\bigl(t^{\langle\theta_i,\lambda
\rangle}\bigr)_{1\leqslant i\leqslant N}$. Let $i\leqslant j$ be two
indices. Then
$$g\cdot(v_i\wedge v_{j+1}\wedge v_{j+2}\wedge\cdots\wedge v_N)=
t^{\langle\theta_{j+1}+\theta_{j+2}+\cdots+\theta_N,\lambda\rangle}
(g\cdot v_i)\wedge v_{j+1}\wedge v_{j+2}\wedge\cdots\wedge v_N.$$
Therefore
\begin{align*}
\val(\langle v_j^*,g\cdot v_i\rangle)+\langle\theta_{j+1}
&+\theta_{j+2}+\cdots+\theta_N,\lambda\rangle\\[3pt]
&=\val(\langle v_j^*\wedge v_{j+1}^*\wedge v_{j+2}^*\wedge\cdots v_N^*,
g\cdot(v_i\wedge v_{j+1}\wedge v_{j+2}\wedge\cdots\wedge v_N)\rangle)\\[3pt]
&=\val(\langle g^{-1}\cdot(v_j^*\wedge v_{j+1}^*\wedge v_{j+2}^*\wedge
\cdots v_N^*),v_i\wedge v_{j+1}\wedge v_{j+2}\wedge\cdots\wedge
v_N\rangle)\\[3pt]
&\geqslant\val(g^{-1}\cdot(v_j^*\wedge v_{j+1}^*\wedge\cdots\wedge
v_N^*))\\[3pt]
&=\langle\theta_j+\theta_{j+1}+\cdots+\theta_N,\lambda+\nu\rangle;
\end{align*}
the last equality here comes from Proposition~\ref{pr:KamnitzForm}
\ref{it:PrKFb}, taking into account that $[g]\in S_{\lambda+\nu}^+$ and
that $v_j^*\wedge v_{j+1}^*\wedge\cdots\wedge v_N^*$ is a highest weight
vector of weight $-(\theta_j+\theta_{j+1}+\cdots+\theta_N)$ in
$\bigwedge^{N-j+1}V(\eta)^*$. By Condition~$A_\eta(\lambda)$, this
implies
$$\val(\langle v_j^*,g\cdot v_i\rangle)\geqslant\langle\theta_j,\lambda
\rangle+\langle\theta_j+\theta_{j+1}+\cdots+\theta_N,\nu\rangle\geqslant
\langle\eta,\lambda\rangle.$$
Therefore $\val\bigl(g_{V(\eta)}\bigr)\geqslant\langle\eta,\lambda\rangle$.
On the other hand, $\val\bigl(g_{V(\eta)}\bigr)\leqslant\val(\langle
v_1^*,g\cdot v_1\rangle)=\langle\eta,\lambda\rangle$. Thus the equality
$\val(g_{V(\eta)})=\langle\eta,\lambda\rangle$ holds for each $\eta\in Y$,
and we conclude by Proposition~\ref{pr:KamnitzForm}~\ref{it:PrKFa} that
$[g]\in\mathscr G_\lambda$.
\end{proof}

\trivlist
\item[\hskip\labelsep{\itshape Proof of Proposition~\ref{pr:IncidPropStrat}.}]
\upshape
We first prove Assertion~\ref{it:PrIPSa}. We let $\mathbb C^\times$ act
on $\mathscr G$ through a dominant and regular coweight $\xi\in\Lambda$.
Let $\lambda,\nu\in\Lambda$ and assume there exists an element $x\in
S_\lambda^+\cap S_\nu^-$. Then
$$[t^\nu]=\lim_{a\to0}a^{-\xi}\cdot x\quad\text{belongs to}\quad
\overline{S_\lambda^+}=\bigcup_{\substack{\mu\in\Lambda\\\lambda
\geqslant\mu}}S_\mu^+.$$
This shows that $\lambda\geqslant\nu$.

If $\mu\in\Lambda$ is enough antidominant, then
$$S_\mu^+\cap S_\mu^-\subseteq S_\mu^+\cap\mathscr G_\mu
=\bigl\{[t^\mu]\bigr\}$$
by Lemma~\ref{le:IncidPropStrat} and Formula~(3.6) in
\cite{MirkovicVilonen04}. Thus $S_\mu^+\cap S_\mu^-=
\bigl\{[t^\mu]\bigr\}$ if $\mu$ is enough antidominant. It
follows that for each $\lambda\in\Lambda$,
$$S_\lambda^+\cap S_\lambda^-=t^{\lambda-\mu}\cdot\bigl(S_\mu^+\cap
S_\mu^-\bigr)=t^{\lambda-\mu}\cdot\bigl\{[t^\mu]\bigr\}=
\bigl\{[t^\lambda]\bigr\}.$$
Assertion~\ref{it:PrIPSb} is proved.

Now let $\nu\in\Lambda$ such that $\nu\geqslant0$. By Lemma
\ref{le:IncidPropStrat}, the property
\begin{equation}
\forall\sigma,\tau\in\Lambda,\quad(0\leqslant\tau\leqslant\nu
\text{ and }\lambda\leqslant\sigma\leqslant\lambda+\nu)\Longrightarrow
(S_{\sigma+\tau}^+\cap S_\sigma^-\subseteq\mathscr G_\sigma)
\label{eq:PfPrIPS}
\end{equation}
holds if $\lambda$ is enough antidominant. We assume that this is the
case and that moreover
$$W\lambda\cap\{\sigma\in\Lambda\mid\sigma\leqslant\lambda+\nu\}=
\{\lambda\}.$$
We now show the equality $S_{\lambda+\nu}^+\cap S_\lambda^-=
S_{\lambda+\nu}^+\cap\mathscr G_\lambda$. Let us take
$x\in S_{\lambda+\nu}^+\cap\mathscr G_\lambda$. Calling $\sigma$ the
coweight such that $x\in S_\sigma^-$, we necessarily have $\lambda
\leqslant\sigma\leqslant\lambda+\nu$ (using Example~\ref{ex:DeterMu}
\ref{it:ExDMb} for the first inequality). Setting $\tau=\lambda+\nu-\sigma$,
we have $0\leqslant\tau\leqslant\nu$ and $x\in S_{\sigma+\tau}^+\cap
S_\sigma^-$, whence $x\in\mathscr G_\sigma$ by our assumption
(\ref{eq:PfPrIPS}). This entails $\sigma\in W\lambda$, then
$\sigma=\lambda$, and thus $x\in S_\lambda^-$. This reasoning shows
$S_{\lambda+\nu}^+\cap\mathscr G_\lambda\subseteq S_{\lambda+\nu}^+\cap
S_\lambda^-$. The converse inclusion also holds (set $\tau=\nu$ and
$\sigma=\lambda$ in~(\ref{eq:PfPrIPS})). Assertion~\ref{it:PrIPSc} is
proved.
\nobreak\noindent$\square$
\endtrivlist

\begin{other*}{Remark}
Assertion~\ref{it:PrIPSb} of Proposition~\ref{pr:IncidPropStrat} can also
be proved in the following way. Let $K$ be the maximal compact subgroup
of the torus $T$. The Lie algebra of $K$ is $\mathfrak k=i(\Lambda
\otimes_{\mathbb Z}\mathbb R)$. The affine Grassmannian $\mathscr G$ is a
K\"ahler manifold and the action of $K$ on $\mathscr G$ is hamiltonian.
Let $\mu:\mathscr G\to\mathfrak k^*$ be the moment map. Fix a dominant
and regular coweight $\xi\in\Lambda$. Then $\mathbb R_+^\times$ acts on
$\mathscr G$ through the map $\mathbb R_+^\times\hookrightarrow\mathbb
C^\times\xrightarrow\xi T$. The map $\langle\mu,i\xi\rangle$ from
$\mathscr G$ to $\mathbb R$ strictly increases along any non-constant
orbit for this $\mathbb R_+^\times$-action. Now take $\lambda\in\Lambda$
and $x\in S_\lambda^+\cap S_\lambda^-$. Then $\lim_{a\to0}a^\xi\cdot
x=\lim_{a\to\infty}a^\xi\cdot x=[t^\lambda]$. Thus $\langle\mu,i\xi\rangle$
cannot increases strictly along the orbit $\mathbb R_+^\times\cdot x$. This
implies that this orbit is constant; in other words, $x=[t^\lambda]$.
\end{other*}

\subsection{Mirkovi\'c-Vilonen cycles}
\label{ss:MVCycles}
Let $\lambda,\nu\in\Lambda$. In order that $S_\nu^+\cap\mathscr
G_\lambda\neq\varnothing$, it is necessary that $[t^\nu]\in
\overline{\mathscr G_\lambda}^T$, hence that $\nu-\lambda\in\mathbb
Z\Phi^\vee$ and that $\nu$ belongs to the convex hull of $W\lambda$
in $\Lambda\otimes_{\mathbb Z}\mathbb R$.

Assume that $\lambda$ is antidominant and denote by $L(w_0\lambda)$ the
irreducible rational representation of $G^\vee$ with lowest weight
$\lambda$. Mirkovi\'c and Vilonen proved that the intersection
$S_\nu^+\cap\mathscr G_\lambda$ is of pure dimension
$\height(\nu-\lambda)$ and has as many irreducible components as
the dimension of the $\nu$-weight subspace of $L(w_0\lambda)$ (Theorem
3.2 and Corollary~7.4 in \cite{MirkovicVilonen04}). From this
result and from Proposition~\ref{pr:IncidPropStrat}~\ref{it:PrIPSc},
one readily deduces the following fact.

\begin{proposition}
\label{pr:DimEstimate}
Let $\lambda,\nu\in\Lambda$ with $\nu\geqslant0$. Then the intersection
$S_{\lambda+\nu}^+\cap S_\lambda^-$ (viewed as a topological subspace of
$\mathscr G$) is noetherian of pure dimension $\height(\nu)$ and has as
many irreducible components as the dimension of the $\nu$-weight subspace
of~$U(\mathfrak n^{+,\vee})$.
\end{proposition}
\begin{proof}
As an abstract topological space, $S_{\lambda+\nu}^+\cap S_\lambda^-$ does
not depend on $\lambda$, because the action of $t^\mu$ on $\mathscr G$
maps $S_{\lambda+\nu}^+\cap S_\lambda^-$ onto $S_{\lambda+\mu+\nu}^+
\cap S_{\lambda+\mu}^-$, for any $\mu\in\Lambda$. We may therefore assume
that $\lambda$ is enough antidominant so that the conclusion of Proposition
\ref{pr:IncidPropStrat}~\ref{it:PrIPSc} holds and that the
$(\lambda+\nu)$-weight space of $L(w_0\lambda)$ has the same dimension
as the $\nu$-weight subspace of $U(\mathfrak n^{+,\vee})$.
The proposition then follows from Mirkovi\'c and Vilonen results.
\end{proof}

If $X$ is a topological space, we denote the set of irreducible
components of $X$ by $\Irr(X)$. For $\lambda,\nu\in\Lambda$, we set
$$\mathscr Z(\lambda)_\nu=\Irr\Bigl(\,\overline{S_\nu^+\cap\mathscr
G_\lambda}\,\Bigr).$$
An element $Z$ in a set $\mathscr Z(\lambda)_\nu$ is called an MV cycle.
Such a $Z$ is necessarily a closed, irreducible and noetherian subset
of $\mathscr G$. It is also $T$-invariant, for the action of the
connected group $T$ on $\overline{S_\nu^+\cap\mathscr G_\lambda}$ does
not permute the irreducible components of this intersection closure.
The coweight $\nu$ can be recovered from $Z$ by the rule $\mu_+(Z)=\nu$;
indeed $Z$ is the closure of an irreducible component $Y$ of $S_\nu^+\cap
\mathscr G_\lambda$, so that $\mu_+(Z)=\mu_+(Y)=\nu$. The union
$$\mathscr Z(\lambda)=\bigsqcup_{\nu\in\Lambda}\mathscr Z(\lambda)_\nu$$
is therefore disjoint.

We finally set
$$\mathscr Z=\bigsqcup_{\substack{\lambda,\nu\in\Lambda\\\lambda
\geqslant\nu}}\Irr\Bigl(\,\overline{S_\lambda^+\cap S_\nu^-}\,\Bigr).$$
Arguing as above, one sees that if $Z$ is an irreducible component of
$\overline{S_\lambda^+\cap S_\nu^-}$, then $\lambda$ and $\nu$
are determined by $Z$ through the equations $\mu_+(Z)=\lambda$ and
$\mu_-(Z)=\nu$. Using Example~\ref{ex:DeterMu}~\ref{it:ExDMa}, one
checks without difficulty that for any irreducible and noetherian
subset $Z$ of $\mathscr G$,
\begin{align}
\overline Z\in\mathscr Z&\Longleftrightarrow\overline Z\text{ is an
irreducible component of }\overline{S_{\mu_+(Z)}^+\cap S_{\mu_-(Z)}^-}
\notag\\[4pt]&\Longleftrightarrow\dim Z=\height(\mu_+(Z)-\mu_-(Z)).
\label{eq:CritMVCycle}
\end{align}

A result of Anderson (Proposition~3 in \cite{Anderson03}) asserts that
for any $\lambda,\nu\in\Lambda$ with $\lambda$ antidominant,
$$\mathscr Z(\lambda)_\nu=\bigl\{Z\in\mathscr Z\big|\mu_+(Z)=\nu,
\mu_-(Z)=\lambda\text{ and }Z\subseteq\overline{\mathscr G_\lambda}
\bigr\}.$$
This fact implies that if $\lambda$ and $\mu$ are two antidominant
coweights such that $\mu-\lambda\in\Lambda_{++}$ and if $Z\in\mathscr
Z(\mu)$, then $t^{\lambda-\mu}\cdot Z\in\mathscr Z(\lambda)$. The set
$\mathscr Z$ appears thus as the right way to stabilize the situation,
namely
$$\mathscr Z=\Biggl\{t^\nu\cdot Z\Biggm|\nu\in\Lambda,Z\in\bigsqcup_{\lambda
\in\Lambda_{++}}\mathscr Z(\lambda)\Biggr\}.$$
It seems therefore legitimate to call MV cycles the elements of $\mathscr Z$.

From now on, our main aim will be to describe MV cycles as precisely as
possible. We treat here the case where $G$ has semisimple rank~$1$.
We set $\mathbb C[t^{-1}]_0^+=\mathbb C[t^{-1}]_0^*=\{0\}$.
For each positive integer $n$, we consider the subsets
$$\mathbb C[t^{-1}]_n^+=\bigl\{a_{-n}t^{-n}+\cdots+a_{-1}t^{-1}
\bigm|(a_{-n},\ldots,a_{-1})\in\mathbb C^n\bigr\}$$
and
$$\mathbb C[t^{-1}]_n^*=\bigl\{a_{-n}t^{-n}+\cdots+a_{-1}t^{-1}
\bigm|(a_{-n},\ldots,a_{-1})\in\mathbb C^n,\ a_{-n}\neq0\bigr\}$$
of $\mathscr K$; these are affine complex varieties. Finally we set
$\mathbb C[t^{-1}]^+=t^{-1}\mathbb C[t^{-1}]=\bigcup_{n\in\mathbb N}
\mathbb C[t^{-1}]_n^+$ and endow it with the inductive limit of the
Zariski topologies on the subspaces $\mathbb C[t^{-1}]_n^+$.

\begin{proposition}
\label{pr:GpsOfRankOne}
Assume that $G$ has semisimple rank~$1$. Let $\nu\in\Lambda$ and denote
the unique simple root by $\alpha$. Then the map $f:p\mapsto
x_{-\alpha}\bigl(pt^{-\langle\alpha,\nu\rangle}\bigr)[t^\nu]$ from
$\mathbb C[t^{-1}]^+$ onto $S_\nu^-$ is a homeomorphism. Moreover
for each $n\in\mathbb N$, the map $f$ induces homeomorphisms
$$\mathbb C[t^{-1}]_n^+\stackrel\simeq\longrightarrow
\overline{S_{\nu+n\alpha^\vee}^+}\cap S_\nu^-\quad\text{and}\quad
\mathbb C[t^{-1}]_n^*\stackrel\simeq\longrightarrow
S_{\nu+n\alpha^\vee}^+\cap S_\nu^-.$$
\end{proposition}
This proposition implies that if $G$ has semisimple rank~$1$, then each
intersection $S_\lambda^+\cap S_\nu^-$ is either empty or irreducible.
In this case thus, the map $Z\mapsto(\mu_+(Z),\mu_-(Z))$ is a bijection
from $\mathscr Z$ onto $\{(\lambda,\nu)\mid\lambda\geqslant\nu\}$,
with inverse bijection $(\lambda,\nu)\mapsto\overline{S_\lambda^+\cap
S_\nu^-}$.

\trivlist
\item[\hskip\labelsep{\itshape Proof of Proposition~\ref{pr:GpsOfRankOne}.}]
\upshape
Let $G$, $\alpha$, $\nu$ and $f$ be as in the statement of the proposition.
The additive group $\mathscr K$ acts transitively on $S_\nu^-$ through the
map $(p,z)\mapsto x_{-\alpha}\bigl(pt^{-\langle\alpha,\nu\rangle}\bigr)z$,
where $p\in\mathscr K$ and $z\in S_\nu^-$. The stabilizer in $\mathscr K$
of $[t^\nu]$ is $\mathscr O$. Since $\mathscr K/\mathscr O\cong\mathbb
C[t^{-1}]^+$, the map $f$ is bijective. It is also continuous.

Now let $n\in\mathbb N$. Set $\lambda=\nu+n\alpha^\vee$; then
$n=\langle\alpha,\lambda-\nu\rangle/2$. Specializing the equality
$$x_{-\alpha}(-a^{-1})=x_\alpha(-a)\,a^{\alpha^\vee}\overline{s_\alpha}
\,x_\alpha(-a)$$
to the value $a=-qt^n$, where $q\in\mathscr O^\times$, multiplying
it on the left by $t^\nu$ and noticing that $(-q)^{\alpha^\vee}
\overline{s_\alpha}\,x_\alpha(qt^n)\in G(\mathscr O)$, we get
$$\bigl[x_{-\alpha}\bigl(q^{-1}t^{-\langle\alpha,\lambda+\nu\rangle/2}
\bigr)\,t^\nu\bigr]=\bigl[x_\alpha\bigl(qt^{\langle\alpha,\lambda+
\nu\rangle/2}\bigr)\,t^{\lambda}\bigr].$$
This equality immediately implies that $f\bigl(\mathbb C[t^{-1}]_n^*
\bigr)\subseteq S_{\nu+n\alpha^\vee}^+\cap S_\nu^-$.
Since
$$\mathbb C[t^{-1}]^+=\bigsqcup_{n\in\mathbb N}\mathbb C[t^{-1}
]_n^*\quad\text{and}\quad S_\nu^-=\bigsqcup_{n\in\mathbb N}
\bigl(S_{\nu+n\alpha^\vee}^+\cap S_\nu^-\bigr),$$
we deduce that $f\bigl(\mathbb C[t^{-1}]_n^*\bigr)=
S_{\nu+n\alpha^\vee}^+\cap S_\nu^-$, and then, using
(\ref{eq:ClosureStratum}), that $f\bigl(\mathbb C[t^{-1}]_n^+\bigr)=
\overline{S_{\nu+n\alpha^\vee}^+}\cap S_\nu^-$. The map $f$ yields thus
a continuous bijection from $\mathbb C[t^{-1}]_n^+$ onto
$\overline{S_{\nu+n\alpha^\vee}^+}\cap S_\nu^-$.

It remains to show the continuity of $f^{-1}$. We may assume without
loss of generality that $\nu=0$. We first look at the particular case
$G=\SL_2$ with its usual pinning. Given an element $p\in\mathscr K$, we
write $p=\{p\}_{<0}+\{p\}_{\geqslant0}$ according to the decomposition
$\mathscr K=\mathbb C[t^{-1}]^+\oplus\mathscr O$, and we denote by $p_0$
the coefficient of $t^0$ in $p$. We consider the subsets
$$\Omega'=\biggl\{\begin{pmatrix}a&b\\c&d\end{pmatrix}\biggm|
a_0\neq0\biggr\}\quad\text{and}\quad
\Omega''=\biggl\{\begin{pmatrix}a&b\\c&d\end{pmatrix}\biggm|
b_0\neq0\biggr\}$$
of $G(\mathscr K)$, and we define maps
$$h':\begin{pmatrix}a&b\\c&d\end{pmatrix}\mapsto\bigl\{c/\{a
\}_{\geqslant0}\bigr\}_{<0}\quad\text{and}\quad
h'':\begin{pmatrix}a&b\\c&d\end{pmatrix}\mapsto\bigl\{d/\{b
\}_{\geqslant0}\bigr\}_{<0}$$
from $\Omega'$ and $\Omega''$, respectively, to $\mathbb C[t^{-1}]^+$.
Certainly, $\Omega'$ and $\Omega''$ are open subsets of $G(\mathscr K)$,
and $h'$ and $h''$ are continuous (see Proposition~1.2 in
\cite{BeauvilleLaszlo94} for details on the inductive system that
defines the topology on $G(\mathscr K)$). We now observe that
$U^-(\mathscr K)G(\mathscr O)\subseteq\Omega'\cup\Omega''$ and that the
map $h:g\mapsto f^{-1}([g])$ from $U^-(\mathscr K)G(\mathscr O)$ to
$\mathbb C[t^{-1}]^+$ is given on $\Omega'\cap U^-(\mathscr K)G(\mathscr
O)$ by the restriction of $h'$, and on $\Omega''\cap U^-(\mathscr K)
G(\mathscr O)$ by the restriction of $h''$. This map $h$ is thus
continuous, and we conclude that $f^{-1}$ is continuous in our
particular case $G=\SL_2$.

The continuity of $f^{-1}$ is then guaranteed whenever $G$ is the product
of $\SL_2$ with a torus. Now any connected reductive group of semisimple
rank~$1$ is isogenous to such a product; the general case follows, because
an isogeny between two connected reductive groups induces a homeomorphism
between the neutral connected components of their respective Grassmannians
(see for instance Section~2 of~\cite{GaussentLittelmann05}).
\nobreak\noindent$\square$
\endtrivlist

\subsection{Parabolic retractions}
\label{ss:ParabRet}
In Section~(5.3.28) of~\cite{BeilinsonDrinfeld??}, Beilinson and Drinfeld
describe a way to relate $\mathscr G$ with the affine Grassmannians of
Levi subgroups of $G$. We rephrase their construction in a slightly less
general context.

Let $P$ be a parabolic subgroup of $G$ which contains $T$, let $M$
be the Levi factor of $P$ that contains $T$, and let $\mathscr P$
and $\mathscr M$ be the affine Grassmannians of $P$ and $M$. The diagram
$G\hookleftarrow P\twoheadrightarrow M$ yields similar diagrams
$G(\mathscr K)\hookleftarrow P(\mathscr K)\twoheadrightarrow M(\mathscr K)$
and $\mathscr G\xleftarrow i\mathscr P\xrightarrow\pi\mathscr M$.
The continuous map $i$ is bijective but is not a homeomorphism in general
($\mathscr P$ has usually more connected components than $\mathscr G$).
We may however define the (non-continuous) map $r_P=\pi\circ i^{-1}$ from
$\mathscr G$ to $\mathscr M$.

The group $P(\mathscr K)$ acts on $\mathscr M$ via the projection
$P(\mathscr K)\twoheadrightarrow M(\mathscr K)$ and acts on $\mathscr G$
via the embedding $P(\mathscr K)\hookrightarrow G(\mathscr K)$. The map
$r_P$ can then be characterized as the unique $P(\mathscr K)$-equivariant
section of the embedding $\mathscr M\hookrightarrow\mathscr G$ that
arises from the inclusion $M\subseteq G$.

For instance, consider the case where $P$ is the Borel subgroup $B^\pm$;
then the Levi factor $M$ is the torus $T$ and the group $P(\mathscr K)$
contains the group $U^\pm (\mathscr K)$. The map $r_{B^\pm}:\mathscr
G\to\mathscr T$, being a $U^\pm(\mathscr K)$-equivariant section of the
embedding $\mathscr T\hookrightarrow\mathscr G$, sends the whole stratum
$S_\lambda^\pm$ to the point $[t^\lambda]$, for each $\lambda\in\Lambda$.

\begin{other}{Remark}
\label{rk:RetracCellsBB}
The map $r_P$ can also be understood in terms of a Bia\l ynicki-Birula
decomposition. Indeed let $\mathfrak g$, $\mathfrak p$ and $\mathfrak t$
be the Lie algebras of $G$, $P$ and $T$. We write $\mathfrak g=\mathfrak
t\oplus\bigoplus_{\alpha\in\Phi}\mathfrak g^\alpha$ for the root
decomposition of $\mathfrak g$ and put $\Phi_P=\{\alpha\in\Phi\mid
\mathfrak g^\alpha\subseteq\mathfrak p\}$. Choosing now $\xi\in\Lambda$
such that
$$\forall\alpha\in\Phi_P,\ \langle\alpha,\lambda\rangle\geqslant0
\quad\text{and}\quad\forall\alpha\in\Phi\setminus\Phi_P,\
\langle\alpha,\lambda\rangle<0,$$
one may check that $r_P(x)=\lim_{\substack{a\to0\\a\in\mathbb C^\times}}
a^\xi\cdot x$ for each $x\in\mathscr G$. This construction justifies the
name of parabolic retraction we give to the map $r_P$.
\end{other}

As noted by Beilinson and Drinfeld (see the proof of Proposition~5.3.29
in~\cite{BeilinsonDrinfeld??}), parabolic retractions enjoy a transitivity
property. Namely considering a pair $(P,M)$ inside $G$ as above and a
pair $(Q,N)$ inside $M$, we get maps $\mathscr G\stackrel{r_P}
\longrightarrow\mathscr M\stackrel{r_Q}\longrightarrow\mathscr N$.
The preimage $R$ of $Q$ by the quotient map $P\twoheadrightarrow M$
is a parabolic subgroup of $G$, and $N$ is the Levi factor of $R$ that
contains $T$. The composition $r_Q\circ r_P$ is a $R(\mathscr
K)$-equivariant section of the embedding $\mathscr N\hookrightarrow
\mathscr G$; it thus coincides with $r_R$.

We will mainly apply these constructions to the case of standard
parabolic subgroups. Let us fix the relevant terminology. For each
subset $J\subseteq I$, we denote by $U^\pm_J$ the subgroup of $G$
generated by the images of the morphisms $x_{\pm\alpha_j}$ for
$j\in J$. We denote the subgroup generated by $T\cup U^+_J\cup U^-_J$
by $M_J$ and we denote the subgroup generated by $B^+\cup M_J$ by
$P_J$. Thus $M_J$ is the Levi factor of $P_J$ that contains $T$.
We shorten the notation and denote the parabolic retraction $r_{P_J}$
simply by $r_J$. The Weyl group of $M_J$ can be identified with the
parabolic subgroup $W_J$ of $W$ generated by the simple reflections
$s_j$ with $j\in J$; we denote the longest element of $W_J$ by $w_{0,J}$.

The Iwasawa decomposition for $M_J$ gives
$$\mathscr M_J=\bigsqcup_{\lambda\in\Lambda}U^\pm_J(\mathscr K)
[t^\lambda].$$
For $\lambda\in\Lambda$, we denote the $U^\pm_J(\mathscr K)$-orbit
of $[t^\lambda]$ by $S_{\lambda,J}^\pm$.
\begin{lemma}
\label{le:RetractStrata}
For each $\lambda\in\Lambda$,\quad$S_\lambda^+=(r_J)^{-1}\bigl(
S_{\lambda,J}^+\bigr)$ \ and\quad$\overline{w_{0,J}}S_{w_{0,J}^{-1}
\lambda}^+=(r_J)^{-1}\bigl(S_{\lambda,J}^-\bigr)$.
\end{lemma}
\begin{proof}
Consider the transitivity property $r_R=r_Q\circ r_P$ of parabolic
retractions written above for $P=P_J$, $M=M_J$ and $N=T$. For the first
formula, one chooses moreover $Q=TU^+_J$, so that $R=B^+$. Recalling
the equality $(r_{B^+})^{-1}\bigl([t^\lambda]\bigr)=S_\lambda^+$ and
its analogue $(r_Q)^{-1}\bigl([t^\lambda]\bigr)=S_{\lambda,J}^+$ for
$\mathscr M_J$, we see that the desired formula simply computes the
preimage of $[t^\lambda]$ by the map $r_R=r_Q\circ r_P$.

For the second formula, one chooses $Q=TU^-_J$, whence
$R=\overline{w_{0,J}}\,B^+\,\overline{w_{0,J}}^{-1}$. Here we have
$$(r_R)^{-1}\Bigl(\bigl[t^\lambda\bigr]\Bigr)=\overline{w_{0,J}}\
(r_{B^+})^{-1}\Bigl(\bigl[t^{w_{0,J}^{-1}\lambda}\bigr]\Bigr)=
\overline{w_{0,J}}\,S_{w_{0,J}^{-1}\lambda}^+$$
and $(r_Q)^{-1}\bigl([t^\lambda]\bigr)=S_{\lambda,J}^-$. Again the
desired formula simply computes the preimage of $[t^\lambda]$ by the
map $r_R=r_Q\circ r_P$.
\end{proof}

To conclude this section, we note that for any $\mathscr K$-point $h$
of the unipotent radical of $P_J$, any $g\in P_J(\mathscr K)$ and any
$x\in\mathscr G$,
\begin{equation}
\label{eq:CommutInRetr}
r_J(gh\cdot x)=(ghg^{-1})\cdot r_J(gx)=r_J(gx),
\end{equation}
because $ghg^{-1}$ is a $\mathscr K$-point of the unipotent radical of
$P_J$ and thus acts trivially on $\mathscr M_J$.

\section{Crystal structure and string parametrizations}
\label{se:CrystStruct}
For each dominant coweight $\lambda$, the set $\mathscr Z(\lambda)$
yields a basis of the rational $G^\vee$-module $L(\lambda)$. One may
therefore expect that $\mathscr Z(\lambda)$ can be turned in a natural
way into a crystal isomorphic to $\mathbf B(\lambda)$. Braverman and
Gaitsgory made this idea precise in~\cite{BravermanGaitsgory01}.
Later in~\cite{BravFinkGait03}, these two authors and Finkelberg
extended this result by endowing $\mathscr Z$ with the structure of
a crystal isomorphic to $\widetilde{\mathbf B(-\infty)}$. We recall
this crucial result in Section~\ref{ss:BFG} and characterize the
crystal operations on $\mathscr Z$ in a suitable way for comparisons
(Proposition~\ref{pr:CritOpCrys}).

We begin Section~\ref{ss:MVCycString} by translating the geometric
definition of Braverman, Finkelberg and Gaitsgory's crystal structure
on $\mathscr Z$ in more algebraic terms (Proposition~\ref{pr:CrystOpAlg}).
From there, we deduce a quite explicit description of MV cycles. More
precisely, let $b\in\mathbf B(-\infty)$ and let $\Xi(t_0\otimes b)$
be the MV cycle that corresponds to $t_0\otimes b\in\widetilde{\mathbf
B(-\infty)}$. Theorem~\ref{th:ParamOpen} exhibits a parametrization of
an open and dense subset of $\Xi(t_0\otimes b)$ by a variety of the form
$(\mathbb C^\times)^m\times\mathbb C^n$; this parametrization generalizes
the description in semisimple rank~1 given in
Proposition~\ref{pr:GpsOfRankOne}.

The next Section~\ref{ss:TildeYic} introduces subsets $\tilde Y_{\mathbf
i,\mathbf c}$ of the affine Grassmannian $\mathscr G$, where $\mathbf
i\in I^l$ and $\mathbf c\in\mathbb Z^l$. When $\mathbf c$ is the string
parameter in direction $\mathbf i$ of an element $b\in\mathbf B(-\infty)$,
the definition of $\tilde Y_{\mathbf i,\mathbf c}$ reflects the
construction in the statement of Theorem~\ref{th:ParamOpen}, so that
$\Xi(t_0\otimes b)=\overline{\tilde Y_{\mathbf i,\mathbf c}}$.
It turns out that the closure $\overline{\tilde Y_{\mathbf i,\mathbf c}}$
is always an MV cycle, even when $\mathbf c$ does not belong to the string
cone in direction $\mathbf i$. Proposition~\ref{pr:TildeYic} presents
a necessary and sufficient condition on
$\overline{\tilde Y_{\mathbf i,\mathbf c}}$ in order that $\mathbf c$
may belong to the string cone; its proof relies on Berenstein and
Zelevinsky's characterization of the string cone in terms
of $\mathbf i$-trails \cite{BerensteinZelevinsky01}.

The introduction of the subsets $\tilde Y_{\mathbf i,\mathbf c}$ finds
its justification in Section~\ref{ss:LuszAlgGeomParam}. Here we use
them to explain how the algebraic-geometric parametrization of
$\mathbf B(-\infty)$ devised by Lusztig in~\cite{Lusztig96} is related
to MV cycles.

In the course of his work on MV polytopes~\cite{Kamnitzer05a,Kamnitzer05b},
Kamnitzer was led to a description of MV cycles similar to the equality
$\Xi(t_0\otimes b)=\overline{\tilde Y_{\mathbf i,\mathbf c}}$, but
starting from the Lusztig parameter of $b$ instead of the string parameter.
In Section~\ref{ss:LinkKamnCons}, we show that the equality and
Kamnitzer's description are in fact equivalent results.

\subsection{Braverman, Finkelberg and Gaitsgory's crystal structure}
\label{ss:BFG}
In Section~13 of~\cite{BravFinkGait03}, Braverman, Finkelberg and
Gaitsgory endow $\mathscr Z$ with the structure of a crystal with an
involution $*$. The main step of their construction is an analysis of
the behavior of MV cycles with respect to the standard parabolic
retractions. For a subset $J\subseteq I$, we denote the analogues
of the maps $\mu_\pm$ for the affine Grassmannian $\mathscr M_J$ by
$\mu_{\pm,J}$. The following theorem is due to Braverman, Finkelberg
and Gaitsgory; we nevertheless recall quickly its proof since we
ground the proof of the forthcoming Propositions~\ref{pr:CritOpCrys}
and~\ref{pr:CrystOpAlg} on it.

\begin{theorem}
\label{th:MVCycFibr}
Let $J$ be a subset of $I$ and let $Z\in\mathscr Z$ be an MV cycle. Set
$$Z_J=\overline{r_J\bigl(Z\cap S_\nu^-\bigr)\cap S_{\rho,J}^-}\quad
\text{and}\quad Z^J=\overline{Z\cap S_\nu^-\cap(r_J)^{-1}([t^\rho])},$$
where $\nu=\mu_-(Z)$ and $\rho=w_{0,J}\,\mu_+(\overline{w_{0,J}}^{-1}Z)$.
Then the map $Z\mapsto(Z_J,Z^J)$ is a bijection from $\mathscr Z$ onto
the set of all pairs $(Z',Z'')$, where $Z'$ is an MV cycle in
$\mathscr M_J$ and $Z''$ is an MV cycle in $\mathscr G$ which satisfy
\begin{equation}
\label{eq:CondFactor}
\mu_{-,J}(Z')=\mu_+(Z'')=w_{0,J}\,\mu_+(\overline{w_{0,J}}^{-1}Z'').
\end{equation}
Under this correspondence, one has
\begin{align*}
\mu_+(Z)&=\mu_{+,J}(Z_J),\\
\mu_-(Z)&=\mu_-(Z^J),\\
w_{0,J}\,\mu_+(\overline{w_{0,J}}^{-1}Z)&=\mu_{-,J}(Z_J)=\mu_+(Z^J)
=w_{0,J}\,\mu_+(\overline{w_{0,J}}^{-1}Z^J).
\end{align*}
\end{theorem}
\begin{proof}
Let us consider three coweights $\lambda,\nu,\rho\in\Lambda$, in the
same coset modulo $\mathbb Z\Phi^\vee$, and unrelated to the
MV cycle $Z$ for the moment. The group $H=U^-_J(\mathscr K)$ acts
on $\mathscr G$, leaving $S_\nu^-$ stable. On the other hand,
$S_{\rho,J}^-$ is the $H$-orbit of $[t^\rho]$; we denote by $K$ the
stabilizer of $[t^\rho]$ in $H$, so that $S_{\rho,J}^-\cong H/K$.
Since the map $r_J$ is $H$-equivariant, the action of $H$ leaves stable
the intersection $S_\nu^-\cap(r_J)^{-1}\bigl(S_{\rho,J}^-\bigr)$, the
action of $K$ leaves stable the intersection $F=S_\nu^-\cap(r_J)^{-1}
([t^\rho])$, and we have a commutative diagram
$$\xymatrix{F\ar@{^{(}->}[r]&H\times_KF\ar[d]\ar[r]^<>(.5)\simeq&
S_\nu^-\cap(r_J)^{-1}\bigl(S_{\rho,J}^-\bigr)\ar[d]^{r_J}\\&H/K
\ar[r]^\simeq&S_{\rho,J}^-.}$$
In this diagram, the two leftmost arrows define a fiber bundle.

By Lemma~\ref{le:RetractStrata}, $F\subseteq S_\rho^+\cap S_\nu^-$;
therefore the dimension of $F$ is at most $\height(\rho-\nu)$.
The group $K$ is connected --- indeed $K=U^-_J(\mathscr K)\cap t^\rho
G(\mathscr O)t^{-\rho}$, so it leaves invariant each irreducible
component of $F$. We thus have a canonical bijection $C\mapsto
\tilde C=H\times_KC$ from $\Irr(F)$ onto $\Irr(H\times_KF)$.
If moreover $X$ is a subspace of $H/K=S_{\rho,J}^-$, then the assignment
$(C,D)\mapsto\tilde C\cap(r_J)^{-1}(D)$ is a bijection from
$\Irr(F)\times\Irr(X)$ onto $\Irr(S_\nu^-\cap(r_J)^{-1}(X))$.
We will apply this fact to $X=S_{\rho,J}^-\cap\overline{S_{\lambda,J}^+}$;
using (\ref{eq:ClosureStratum}) and Proposition~\ref{pr:DimEstimate}, one
sees easily that $X$ has then dimension at most $\height(\lambda-\rho)$.
Since $\tilde C\cap(r_J)^{-1}(D)$ is a fiber bundle with fiber $C$ and
base $D$, its dimension is
$$\dim C+\dim D\leqslant\height(\rho-\nu)+\height(\lambda-\rho)=
\height(\lambda-\nu).$$

Now let $Z$ be an MV cycle and set $\lambda=\mu_+(Z)$, $\nu=\mu_-(Z)$
and $\rho=w_{0,J}\,\mu_+(\overline{w_{0,J}}^{-1}Z)$ in the previous
setting. By Proposition~\ref{pr:DefMu} and Lemma~\ref{le:RetractStrata},
$$Z\cap S_\nu^-\quad\text{and}\quad\overline{w_{0,J}}
\Bigl(\overline{w_{0,J}}^{-1}Z\cap S_{w_{0,J}^{-1}\rho}^+\Bigr)=
Z\cap(r_J)^{-1}\bigl(S_{\rho,J}^-\bigr)$$
are open and dense subsets in $Z$. Thus $\dot Z=Z\cap S_\nu^-\cap
(r_J)^{-1}\bigl(S_{\rho,J}^-\bigr)$ is a closed irreducible subset
of $S_\nu^-\cap(r_J)^{-1}\bigl(S_{\rho,J}^-\bigr)$ of dimension
$\dim Z=\height(\lambda-\nu)$; this subset $\dot Z$ is actually contained
in $S_\nu^-\cap(r_J)^{-1}(X)$, because $\dot Z\subseteq Z\subseteq
\overline{S_\lambda^+}$. It is therefore an irreducible component
$\tilde C\cap(r_J)^{-1}(D)$, with moreover $\dim C=\height(\rho-\nu)$
and $\dim D=\height(\lambda-\rho)$.

One observes then that $[t^\rho]\in D$, because $D$ is a closed and
$T$-invariant subset of $S_{\rho,J}^-$. Then
$$C=\tilde C\cap(r_J)^{-1}([t^\rho])=\dot Z\cap(r_J)^{-1}([t^\rho])=
Z\cap S_\nu^-\cap(r_J)^{-1}([t^\rho]),$$
and thus, by Lemma~\ref{le:RetractStrata}, $C\subseteq S_\nu^-\cap
S_\rho^+\cap\overline{w_{0,J}}\,S_{w_{0,J}^{-1}\rho}^+$. Therefore
$\mu_-(C)=\nu$ and $\mu_+(C)=w_{0,J}\,\mu_+(\overline{w_{0,J}}^{-1}C)
=\rho$; Equivalence~(\ref{eq:CritMVCycle}) and the estimate $\dim
C=\height(\rho-\nu)$ imply then that $\overline C$ is an MV cycle.
On the other hand, the relations $\mu_{+,J}(D)\leqslant\lambda$,
$\mu_{-,J}(D)=\rho$ and $\dim D=\height(\lambda-\rho)$ imply altogether
that $\overline D$ is an MV cycle in $\mathscr M_J$ and that
$\mu_{+,J}(D)=\lambda$. Moreover
$$D=r_J(\dot Z)=r_J\bigl(Z\cap S_\nu^-\bigr)\cap S_{\rho,J}^-.$$
Thus $Z_J=\overline D$ and $Z^J=\overline C$ satisfy the conditions
stated in the theorem.

Conversely, given $Z'$ and $Z''$ as in the statement of the theorem, we
take $\lambda=\mu_{+,J}(Z')$, $\nu=\mu_-(Z'')$ and $\rho=\mu_{-,J}(Z')$
in the construction above, and we set $C=Z''\cap F$, $D=Z'\cap
S_{\rho,J}^-$ and $\dot Z=\tilde C\cap(r_J)^{-1}(D)$. Then $C$ is an
open and dense subset in $Z''$; it is therefore irreducible with the
same dimension as $Z''$, namely $\height(\rho-\nu)$. Since it is a
closed subset of $F$, $C$ is an irreducible component of $F$. Likewise
$D$ has dimension $\height(\lambda-\rho)$ and is an irreducible
component of $X=S_{\rho,J}^-\cap\overline{S_{\lambda,J}^+}$. The first
part of the reasoning above implies thus that $\dot Z$ is irreducible
of dimension $\dim C+\dim D=\height(\lambda-\nu)$. Since
$\mu_+(\dot Z)=\lambda$ and $\mu_-(\dot Z)=\nu$, it follows from
Equivalence~(\ref{eq:CritMVCycle}) that $Z=\overline{\dot Z}$ is an
MV cycle.

It is then routine to check that the two maps $Z\mapsto(Z_J,Z^J)$ and
$(Z',Z'')\mapsto Z$ are mutually inverse bijections.
\end{proof}

We are now ready to define Braverman, Finkelberg and Gaitsgory's crystal
structure on $\mathscr Z$. Let $Z$ be an MV cycle. We set
$$\wt(Z)=\mu_+(Z).$$
Given $i\in I$, we apply Theorem~\ref{th:MVCycFibr} to $Z$ and $J=\{i\}$.
We set $\rho=s_i\,\mu_+(\overline{s_i}^{-1}Z)$ and get
a decomposition $\bigl(Z_{\{i\}},Z^{\{i\}}\bigr)$ of $Z$. Then we set
$$\varepsilon_i(Z)=\biggl\langle\alpha_i,\frac{-\mu_+(Z)-\rho}2\biggr
\rangle\quad\text{and}\quad\varphi_i(Z)=\biggl\langle\alpha_i,
\frac{\mu_+(Z)-\rho}2\biggr\rangle.$$
Since $\mu_+(Z)-\rho=\mu_{+,\{i\}}\bigl(Z_{\{i\}}\bigr)-
\mu_{-,\{i\}}\bigl(Z_{\{i\}}\bigr)$ belongs to $\mathbb N\alpha_i^\vee$,
the definition for $\varphi_i(Z)$ is equivalent to the equation
\begin{equation}
\label{eq:DefPhiBFG}
\mu_+(Z)-\rho=\varphi_i(Z)\,\alpha_i^\vee.
\end{equation}

The MV cycles $\tilde e_iZ$ and $\tilde f_iZ$ are defined by the
following requirements:
\begin{align*}
\mu_+(\tilde e_iZ)&=\mu_+(Z)+\alpha_i^\vee,\\
\mu_+(\tilde f_iZ)&=\mu_+(Z)-\alpha_i^\vee,\\
(\tilde e_iZ)^{\{i\}}&=(\tilde f_iZ)^{\{i\}}=Z^{\{i\}};
\end{align*}
if $\mu_+(Z)=\rho$, that is, if $\varphi_i(Z)=0$, then we set
$\tilde f_iZ=0$.

These conditions do define the MV cycles $\tilde e_iZ$ and $\tilde f_iZ$.
Indeed they prescribe the components $(\tilde e_iZ)^{\{i\}}$ and
$(\tilde f_iZ)^{\{i\}}$ and require
\begin{align*}
\mu_{+,\{i\}}\bigl((\tilde e_iZ)_{\{i\}}\bigr)&=\mu_+(\tilde e_iZ)
=\mu_+(Z)+\alpha^\vee_i=\mu_{+,\{i\}}\bigl(Z_{\{i\}}\bigr)+
\alpha_i^\vee\\[3pt]\mu_{-,\{i\}}\bigl((\tilde e_iZ)_{\{i\}}\bigr)&=
\mu_+\bigl((\tilde e_iZ)^{\{i\}}\bigr)=\mu_+\bigl(Z^{\{i\}}\bigr)=
\mu_{-,\{i\}}\bigl(Z_{\{i\}}\bigr)
\end{align*}
and
\begin{align*}
\mu_{+,\{i\}}\bigl((\tilde f_iZ)_{\{i\}}\bigr)&=\mu_+(\tilde f_iZ)
=\mu_+(Z)-\alpha^\vee_i=\mu_{+,\{i\}}\bigl(Z_{\{i\}}\bigr)-
\alpha_i^\vee\\[3pt]\mu_{-,\{i\}}\bigl((\tilde f_iZ)_{\{i\}}\bigr)&=
\mu_+\bigl((\tilde f_iZ)^{\{i\}}\bigr)=\mu_+\bigl(Z^{\{i\}}\bigr)=
\mu_{-,\{i\}}\bigl(Z_{\{i\}}\bigr).
\end{align*}
These latter equations fully determine the components $(\tilde
e_iZ)_{\{i\}}$ and $(\tilde f_iZ)_{\{i\}}$ because $M_{\{i\}}$
has semisimple rank~1 (see the comment after the statement of
Proposition~\ref{pr:GpsOfRankOne}).

One checks without difficulty that $\mathscr Z$, endowed with these
maps $\wt$, $\varepsilon_i$, $\varphi_i$, $\tilde e_i$ and $\tilde f_i$,
satisfies Kashiwara's axioms of a crystal. On the other hand, let
$g\mapsto g^t$ be the antiautomorphism of $G$ that fixes $T$ pointwise
and that maps $x_{\pm\alpha}(a)$ to $x_{\mp\alpha}(a)$ for each simple
root $\alpha$ and each $a\in\mathbb C$. Then the involutive automorphism
$g\mapsto(g^t)^{-1}$ of $G$ extends to $G(\mathscr K)$ and induces an
involution on $\mathscr G$, which we denote by $x\mapsto x^*$. The
image of an MV cycle $Z$ under this involution is an MV cycle $Z^*$.
The properties of this involution $Z\mapsto Z^*$ with respect to the
crystal operations allow Braverman, Finkelberg and
Gaitsgory~\cite{BravFinkGait03} to establish the existence of an
isomorphism of crystals $\Xi:\widetilde{\mathbf B(-\infty)}
\stackrel\simeq\longrightarrow\mathscr Z$. This isomorphism is unique
and is compatible with the involutions $*$ on $\widetilde{\mathbf
B(-\infty)}$ and $\mathscr Z$. One checks that
\begin{xalignat}2
\Xi(t_\lambda\otimes1)&=\bigl\{[t^\lambda]\bigr\},&
\mu_-\bigl(\Xi(t_\lambda\otimes b)\bigr)&=\lambda,
\label{eq:PptyXi}\\[3pt]
\Xi(t_\lambda\otimes b)&=t^\lambda\cdot\Xi(t_0\otimes b),&
\dim\,\Xi(t_\lambda\otimes b)&=\height(\wt(b)),\notag
\end{xalignat}
for all $\lambda\in\Lambda$ and $b\in\mathbf B(-\infty)$.

The following proposition gives a useful criterion which says when two
MV cycles are related by an operator $\tilde e_i$.
\begin{proposition}
\label{pr:CritOpCrys}
Let $Z$ and $Z'$ be two MV cycles in $\mathscr G$ and let $i\in I$. Then
$Z'=\tilde e_iZ$ if and only if the four following conditions hold:
\begin{align*}
\mu_-(Z')&=\mu_-(Z),\\
s_i\,\mu_+(\overline{s_i}^{-1}Z')&=s_i\,\mu_+(\overline{s_i}^{-1}Z),\\
\mu_+(Z')&=\mu_+(Z)+\alpha_i^\vee,\\
Z'&\supseteq Z.
\end{align*}
\end{proposition}
\begin{proof}
We first prove that the conditions in the statement of the proposition
are sufficient to ensure that $Z'=\tilde e_iZ$. We assume that the two
MV cycles $Z$ and $Z'$ enjoy the conditions above and we set
\begin{align*}
\rho&=s_i\,\mu_+(\overline{s_i}^{-1}Z)=s_i\,\mu_+(\overline{s_i}^{-1}Z'),\\
\nu&=\mu_-(Z)=\mu_-(Z'),\\
F&=S_\nu^-\cap(r_{\{i\}})^{-1}([t^\rho]).
\end{align*}
The proof of Theorem~\ref{th:MVCycFibr} tells us that $C=Z\cap F$
and $C'=Z'\cap F$ are two irreducible components of $F$. The condition
$Z'\supseteq Z$ entails then $C'\supseteq C$, and thus $C'=C$. It
follows that
$$Z^{\{i\}}=\overline{C^{\vphantom\prime}}=\overline{C'}=Z^{\prime\{i\}}.$$
This being known, the assumption $\mu_+(Z')=\mu_+(Z)+\alpha_i^\vee$
implies $Z'=\tilde e_iZ$.

Conversely, assume that $Z'=\tilde e_iZ$. Routine arguments show then
that the three first conditions in the statement of the proposition hold.
Setting $\rho$, $\nu$, $F$, $C$ and $C'$ as in the first part of the proof,
we get
$$C=\overline{C^{\vphantom\prime}}\cap F=Z^{\{i\}}\cap F=Z^{\prime\{i\}}
\cap F=\overline{C'}\cap F=C'.$$
On the other hand, set $D=Z_{\{i\}}\cap S_{\rho,\{i\}}^-$ and
$D'=Z'_{\{i\}}\cap S_{\rho,\{i\}}^-$. Using
Proposition~\ref{pr:GpsOfRankOne}, we see that
$$D=\overline{S_{\mu_+(Z),\{i\}}^+\cap S_{\rho,\{i\}}^-}\cap
S_{\rho,\{i\}}^-=\overline{S_{\mu_+(Z),\{i\}}^+}\cap S_{\rho,\{i\}}^-$$
is contained in
$$D'=\overline{S_{\mu_+(Z'),\{i\}}^+\cap S_{\rho,\{i\}}^-}\cap
S_{\rho,\{i\}}^-=\overline{S_{\mu_+(Z'),\{i\}}^+}\cap S_{\rho,\{i\}}^-.$$
Adopting the notation $\tilde C$ from the proof of Theorem
\ref{th:MVCycFibr}, we deduce that $\tilde C\cap(r_{\{i\}})^{-1}(D)$
is contained in $\tilde C\cap(r_{\{i\}})^{-1}(D')$. The closure $Z$ of
the first set is thus contained in the closure $Z'$ of the second set.
\end{proof}

For each dominant coweight $\lambda\in\Lambda_{++}$, the two sets $\mathbf
B(\lambda)$ and $\mathscr Z(\lambda)$ have the same cardinality; indeed
they both index bases of two isomorphic vector spaces, namely the
rational irreducible $G^\vee$-module with highest weight $\lambda$ and
the intersection cohomology of $\overline{\mathscr G_\lambda}$,
respectively. More is true: in~\cite{BravermanGaitsgory01}, Braverman
and Gaitsgory endow $\mathscr Z(\lambda)$ with the structure of a
crystal and show the existence of an isomorphism of crystals
$\Xi(\lambda):\mathbf B(\lambda)\stackrel\simeq\longrightarrow\mathscr
Z(\lambda)$ (see \cite{BravermanGaitsgory01}, p.~569).

\begin{proposition}
\label{pr:CompBGBFG}
The following diagram commutes:
$$\xymatrix{\mathbf B(\lambda)\ar@{_{(}->}[d]_{\iota_{w_0\lambda}}
\ar[r]^{\Xi(\lambda)}&\mathscr Z(\lambda)\ar@{_{(}->}[d]\\
\mathbf T_{w_0\lambda}\otimes\mathbf B(-\infty)\ar[r]_<>(.5)
\Xi&\mathscr Z.}$$
\end{proposition}
\begin{proof}
Let $Z,Z'\in\mathscr Z(\lambda)$ and assume that $Z'$ is the image
of $Z$ by the crystal operator defined in Section~3.3 of
\cite{BravermanGaitsgory01}. The definition of this operator is so
similar to the definition of our (in fact, Braverman, Finkelberg and
Gaitsgory's) crystal operator $\tilde e_i$ that a slight modification
of the proof of Proposition~\ref{pr:CritOpCrys} yields
\begin{align*}
\mu_-(Z')&=\mu_-(Z),\\
s_i\,\mu_+(\overline{s_i}^{-1}Z')&=s_i\,\mu_+(\overline{s_i}^{-1}Z),\\
\mu_+(Z')&=\mu_+(Z)+\alpha_i^\vee,\\
Z'&\supseteq Z.
\end{align*}
By Proposition~\ref{pr:CritOpCrys}, this implies that $Z'$ is the image of
$Z$ by our crystal operator $\tilde e_i$. In other words, the inclusion
$\mathscr Z(\lambda)\hookrightarrow\mathscr Z$ is an embedding of
crystals when $\mathscr Z(\lambda)$ is endowed with the crystal structure
from \cite{BravermanGaitsgory01}.

Thus both maps $\Xi\circ\iota_{w_0\lambda}$ and $\Xi(\lambda)$ are crystal
embeddings of $\mathbf B(\lambda)$ into $\mathscr Z$. Also both maps send
the lowest weight element $b_\low$ of $\mathbf B(\lambda)$ onto the MV cycle
$\bigl\{[t^{w_0\lambda}]\bigr\}$. The proposition then follows from the
fact that each element of $\mathbf B(\lambda)$ can be obtained by
applying a sequence of crystal operators to $b_\low$.
\end{proof}

\begin{other*}{Remark}
One can establish the equality $\Xi\circ\iota_{w_0\lambda}(\mathbf
B(\lambda))=\mathscr Z(\lambda)$ without using Braverman and Gaitsgory's
isomorphism $\Xi(\lambda)$ by the following direct argument.
Let $Z\in\mathscr Z(\lambda)$. Certainly $\mu_-(Z)=w_0\lambda$, so
by Equation~(\ref{eq:PptyXi}), $\Xi^{-1}(Z)$ may be written
$t_{w_0\lambda}\otimes b$ with $b\in\mathbf B(-\infty)$.
Take $i\in I$ and set $\rho=s_i\,\mu_-(\overline{s_i}^{-1}Z)$.
Then $\overline{s_i}^{-1}Z$ meets $S_{s_i^{-1}\rho}^-$, and thus
$\Bigl[t^{s_i^{-1}\rho}\Bigr]$ belongs to $\overline{s_i}^{-1}Z$,
for $\overline{s_i}^{-1}Z$ is closed and $T$-stable. From the
inclusion $Z\subseteq\overline{\mathscr G_\lambda}$, we then deduce
that $[t^\rho]\in\overline{\mathscr G_\lambda}$. Using Equation
(\ref{eq:ClosureOrbit}) and the description $(\mathscr G_\mu)^T
=\bigl\{[t^{w\mu}]\bigm|w\in W\bigr\}$ (see the proof of Proposition
\ref{pr:DefMu}), this yields
$$\rho\in\bigl\{w\mu\bigm|w\in W,\ \mu\in\Lambda_{++}\text{ such
that }\lambda\geqslant\mu\bigr\}.$$
On the other side,
$$\rho-w_0\lambda=s_i\,\mu_-\bigl(\overline{s_i}^{-1}Z\bigr)-\mu_-(Z)=
\mu_+(Z^*)-s_i\,\mu_+\bigl(\overline{s_i}^{-1}Z^*\bigr)=\varphi_i(Z^*)
\alpha_i^\vee.$$
These two facts together entail $\varphi_i(Z^*)\leqslant\langle
\alpha_i,-w_0\lambda\rangle$. Since
$$\varphi_i(Z^*)=\varphi_i(\Xi^{-1}(Z^*))=\varphi_i(\Xi^{-1}(Z)^*)
=\varphi_i((t_{w_0\lambda}\otimes b)^*)=\varphi_i(t_{-w_0\lambda-\wt(b)}
\otimes b^*)=\varphi_i(b^*),$$
we obtain $\varphi_i(b^*)\leqslant\langle\alpha_i,-w_0\lambda
\rangle$. This inequality holds for each $i\in I$, therefore the
element $t_{w_0\lambda}\otimes b$ belongs to $\iota_{w_0\lambda}
(\mathbf B(\lambda))$. We have thus established the inclusion
$\Xi^{-1}(\mathscr Z(\lambda))\subseteq\iota_{w_0\lambda}(\mathbf
B(\lambda))$. Since $\mathbf B(\lambda)$ and $\mathscr Z(\lambda)$
have the same cardinality, this inclusion is an equality.
\end{other*}

\subsection{Description of an MV cycle from the string parameter}
\label{ss:MVCycString}
We begin this section with a proposition that translates Braverman,
Finkelberg and Gaitsgory's geometrical definition for the crystal
operation $\tilde e_i$ into a more algebraic language. This proposition
comes in to flavors: Statement~\ref{it:PrCOAa} is terse, whereas
Statement~\ref{it:PrCOAb} is verbose but yields more refined information.
We recall that the notations $\mathbb C[t^{-1}]_k^+$ and $\mathbb
C[t^{-1}]_k^*$ have been defined in Section~\ref{ss:MVCycles}.

\begin{proposition}
\label{pr:CrystOpAlg}
Let $Z$ be an MV cycle, let $i\in I$, let $k\in\mathbb N$, and set
$Z'=\tilde e_i^k(Z)$.
\begin{enumerate}
\item\label{it:PrCOAa}
For each $p\in\mathscr O$, the action of $y_i\bigl(pt^{\varepsilon_i(Z)}
\bigr)$ stabilizes $Z$. The MV cycle $Z'$ is the closure of
$$\bigl\{y_i(p)\,z\bigm|z\in Z\text{ and }p\in\mathscr K^\times\text{
such that }\val(p)=-k+\varepsilon_i(Z)\bigr\}.$$
\item\label{it:PrCOAb}
Set $\nu=\mu^-(Z)$, $\rho=s_i\mu_+(\overline{s_i}^{-1}Z)$,
$\dot Z=Z\cap S_\nu^-\cap\Bigl(\overline{s_i}\,S_{s_i^{-1}\rho}^+\Bigr)$
and $\dot Z'=Z'\cap S_\nu^-\cap\Bigl(\overline{s_i}\,S_{s_i^{-1}\rho}^+
\Bigr)$. Then the map $f:(p,z)\mapsto y_i\bigl(pt^{\varepsilon_i(Z)}\bigr)
z$ is a homeomorphism from $\mathbb C[t^{-1}]_k^+\times\dot Z$ onto
$\dot Z'$. If moreover $\rho=\mu_+(Z)$, then $\dot Z=Z\cap
S_\nu^-\cap S_{\mu_+(Z)}^+$ and $f$ induces a homeomorphism from
$\mathbb C[t^{-1}]_k^*\times\dot Z$ onto an open and dense subset of
$Z'\cap S_\nu^-\cap S_{\mu_+(Z')}^+$.
\end{enumerate}
\end{proposition}
\begin{proof}
We begin with the proof of Statement~\ref{it:PrCOAb}. Let $Z$ be an
MV cycle and let $i\in I$. We adopt the notation used in the proof of
Theorem~\ref{th:MVCycFibr}, with here $J=\{i\}$. We set
$\lambda=\mu_+(Z)$, $\nu=\mu_-(Z)$, $\rho=s_i\mu_+(\overline{s_i}^{-1}Z)$,
$n=\varphi_i(Z)=\langle\alpha_i,\lambda-\rho\rangle/2$,
$F=S_\nu^-\cap(r_{\{i\}})^{-1}([t^\rho])$ and $X=S_{\rho,\{i\}}^-
\cap\overline{S_{\lambda,\{i\}}^+}$. Then
$$C=Z\cap S_\nu^-\cap(r_{\{i\}})^{-1}([t^\rho])\quad\text{and}\quad
D=r_{\{i\}}\bigl(Z\cap S_\nu^-\bigr)\cap S_{\rho,\{i\}}^-$$
are irreducible components of $F$ and $X$, respectively. Proposition
\ref{pr:GpsOfRankOne} implies then that $D=X$ and that the map
$h:p\mapsto y_i\bigl(pt^{-\langle\alpha_i,\rho\rangle}\bigr)[t^\rho]$
from $\mathscr K$ to $\mathscr M_{\{i\}}$ induces a homeomorphism from
$\mathbb C[t^{-1}]_n^+$ onto $D$.

Let $k\in\mathbb N$ and set $D'=S_{\rho,\{i\}}^-\cap
\overline{S_{\lambda+k\alpha_i^\vee,\{i\}}^+}$. Then $h$ induces a
homeomorphism from $\mathbb C[t^{-1}]_{n+k}^+$ onto $D'$. Since
$-\langle\alpha_i,\rho\rangle=\varepsilon_i(Z)+n$, it follows that
the map $g:(p,x)\mapsto y_i\bigl(pt^{\varepsilon_i(Z)}\bigr)x$
from $\mathscr K\times\mathscr M_{\{i\}}$ to $\mathscr M_{\{i\}}$
induces a homeomorphism from $\mathbb C[t^{-1}]_k^+\times D$ onto $D'$.
Now set
$$Z'=\tilde e_i^k(Z),\quad\dot Z=Z\cap S_\nu^-\cap\Bigl(\overline{s_i}
\,S_{s_i^{-1}\rho}^+\Bigr)\quad\text{and}\quad\dot Z'=Z'\cap S_\nu^-
\cap\Bigl(\overline{s_i}\,S_{s_i^{-1}\rho}^+\Bigr).$$
The proof of Theorem~\ref{th:MVCycFibr} gives us $\dot Z=\tilde C\cap
(r_{\{i\}})^{-1}(D)$ and $\dot Z'=\tilde C\cap(r_{\{i\}})^{-1}(D')$.
Consider the map $f:(p,z)\mapsto y_i\bigl(pt^{\varepsilon_i(Z)}\bigr)z$
from $\mathscr K\times\mathscr G$ to $\mathscr G$. Using that the
action of the group $y_i(\mathscr K)$ stabilizes $\tilde C$ and
commutes with the parabolic retraction $r_{\{i\}}$, we conclude that
$f$ induces a homeomorphism from $\mathbb C[t^{-1}]_k^+\times\dot Z$
onto $\dot Z'$. The first assertion in Statement~\ref{it:PrCOAb} is thus
shown.

Suppose now that $\lambda=\rho$, and denote by $\mathscr N$ the
connected component of $\mathscr M_{\{i\}}$ that contains $[t^\rho]$.
By Lemma~\ref{le:RetractStrata}, $r_{\{i\}}(Z)\cap\mathscr N$
is contained in both $\overline{S_{\lambda,\{i\}}^+}$ and
$\overline{S_{\rho,\{i\}}^-}$, hence in their intersection $\bigl\{
[t^\rho]\bigr\}$. This shows that $r_{\{i\}}(Z)\cap S_{\lambda,\{i\}}^+
=\bigl\{[t^\rho]\bigr\}=r_{\{i\}}(Z)\cap S_{\rho,\{i\}}^-$, and thus
that $Z\cap S_\lambda^+=Z\cap\Bigl(\overline{s_i}\,S_{s_i^{-1}\rho,
\{i\}}^+\Bigr)$, again by Lemma~\ref{le:RetractStrata}. Therefore
$\dot Z=Z\cap S_\nu^-\cap S_\lambda^+$. Now if $k=0$, then
$$f\bigl(\mathbb C[t^{-1}]_k^*\times\dot Z\bigr)=\dot Z=Z\cap
S_\nu^-\cap S_\lambda^+=Z'\cap S_\nu^-\cap S_{\mu_+(Z')}^+.$$
And if $k>0$, then by Proposition~\ref{pr:GpsOfRankOne}
$$g(\mathbb C[t^{-1}]_k^*\times D)=h(\mathbb C[t^{-1}]_{n+k}^*)=
S_{\rho,\{i\}}^-\cap S_{\lambda+k\alpha_i^\vee,\{i\}}^+=D'\cap
S_{\lambda+k\alpha_i^\vee,\{i\}}^+,$$
and thus by Lemma~\ref{le:RetractStrata}
$$f(\mathbb C[t^{-1}]_k^*\times\dot Z)=\dot Z'\cap S_{\mu_+(Z')}^+,$$
which is an open subset of $Z'\cap S_\nu^-\cap S_{\mu_+(Z')}^+$. This
concludes the proof of Statement~\ref{it:PrCOAb}.

We now turn to the proof of Statement~\ref{it:PrCOAa}. We first observe
that $h(\mathscr O)=\bigl\{[t^\rho]\bigr\}$. Let $p\in\mathscr O$ and
write $pt^{-n}=q+r$, with $q\in\mathbb C[t^{-1}]_n^+$ and $r\in\mathscr
O$. For each $x\in D$, we can find $s\in\mathbb C[t^{-1}]_n^+$ such that
$x=h(s)$, and then
$$y_i\bigl(pt^{\varepsilon_i(Z)}\bigr)\cdot x=y_i\bigl((q+r)t^{-\langle
\alpha_i,\rho\rangle}\bigr)\cdot h(s)=h(q+r+s)=h(q+s)\cdot h(r)=h(q+s)$$
belongs to $D$. The action of $y_i\bigl(pt^{\varepsilon_i(Z)}\bigr)$
therefore stabilizes $D$. Since it stabilizes also $\tilde C$ and commutes
with $r_{\{i\}}$, it stabilizes $\dot Z$. We conclude that it stabilizes
$\overline{\dot Z}=Z$.

Using this, we see that
$$\bigl\{y_i(p)\,z\bigm|z\in Z\text{ and }p\in\mathscr K^\times\text{
such that }\val(p)=-k+\varepsilon_i(Z)\bigr\}=f\bigl(\mathbb
C[t^{-1}]_k^*\times Z\bigr).$$
This set has the same closure as $f(\mathbb C[t^{-1}]_k^*\times\dot Z)$,
namely $Z'$. This completes the proof of Statement~\ref{it:PrCOAa}.
\end{proof}

We now recall the definition of the string parameter of an element in
$\mathbf B(-\infty)$. To each sequence $\mathbf i=(i_1,\ldots,i_l)$ of
elements of $I$, we associate an injective map $\Psi_{\mathbf i}$ from
$\mathbf B(-\infty)$ to $\mathbb N^l\times\mathbf B(-\infty)$ by the
following recursive definition:
\begin{itemize}
\item If $l=0$, then $\Psi_{()}:\mathbf B(-\infty)\to\mathbf B(-\infty)$
is the identity map.
\item If $l>1$ and $b\in\mathbf B(-\infty)$, then $\Psi_{\mathbf i}(b)=
\bigl(c_1,\Psi_{\mathbf j}(\tilde f_{i_1}^{c_1}b)\bigr)$, where
$c_1=\varphi_{i_1}(b)$ and $\mathbf j=(i_2,\ldots,i_l)$.
\end{itemize}
To the sequence $\mathbf i$, one also associates recursively an element
$w_{\mathbf i}\in W$ by setting $w_{()}=1$ and asking that $w_{\mathbf i}$
is the longest of the two elements $w_{\mathbf j}$ and
$s_{i_1}w_{\mathbf j}$, where $\mathbf j=(i_2,\ldots,i_l)$ as above.
Finally, one defines the subset
$$\mathbf B(-\infty)_{\mathbf i}=\bigl\{b\in\mathbf B(-\infty)\bigm|
\exists(k_1,\ldots,k_l)\in\mathbb N^l,\ b=\tilde e_{i_1}^{k_1}
\cdots\tilde e_{i_l}^{k_l}1\bigr\}.$$

\noindent From Kashiwara's work on Demazure modules~\cite{Kashiwara93b}
(see also Section~12.4 in~\cite{Kashiwara95}), one deduces that:
\begin{itemize}
\item $\mathbf B(-\infty)_{\mathbf i}$ depends only on $w_{\mathbf i}$
and not on $\mathbf i$.
\item If $\mathbf i$ is a reduced decomposition of the longest element
$w_0$ of $W$, then $\mathbf B(-\infty)_{\mathbf i}=\mathbf B(-\infty)$.
\item $\mathbf B(-\infty)_{\mathbf i}$ is the set of all $b\in
\mathbf B(-\infty)$ such that $\Psi_{\mathbf i}(b)$ has the form
$\bigl(\mathbf c_{\mathbf i}(b),1\bigr)$ for a certain
$\mathbf c_{\mathbf i}(b)\in\mathbb N^l$.
\end{itemize}
The map $\mathbf c_{\mathbf i}:\mathbf B(-\infty)_{\mathbf i}\to\mathbb
N^l$ implicitly defined in the third item above is called the string
parametrization in the direction $\mathbf i$. Its image is called
the string cone and is denoted by $\mathcal C_{\mathbf i}$.

The next theorem affords an explicit description of the MV cycle
$\Xi(t_0\otimes b)$ from the string parameter of $b$. It shows in
particular that MV cycles are rational varieties, a fact already
known from Gaussent and Littelmann's work (see for instance Theorem~4
in~\cite{GaussentLittelmann05}).

\begin{theorem}
\label{th:ParamOpen}
Let $\mathbf i\in I^l$ and $b\in\mathbf B(-\infty)_{\mathbf i}$.
Write $\mathbf c_{\mathbf i}(b)=(c_1,\ldots,c_l)$, set
$$e_j=-\sum_{k=j+1}^lc_k\langle\alpha_{i_j},\alpha_{i_k}^\vee\rangle$$
for each $j\in\{1,\ldots,l\}$, and set $Z=\Xi(t_0\otimes b)$.
Then the map
$$(p_1,\ldots,p_l)\mapsto\bigl[y_{i_1}(p_1t^{e_1})\cdots y_{i_l}
(p_lt^{e_l})\bigr]$$
is an embedding of $\mathbb C[t^{-1}]_{c_1}^*\times\cdots\times
\mathbb C[t^{-1}]_{c_l}^*$ as an open and dense subset of
$Z\cap S_{\mu_+(Z)}^+\cap S_{\mu_-(Z)}^-$.
\end{theorem}
\begin{proof}
We use induction on the length $l$ of the sequence $\mathbf i$.
The assertion certainly holds when $l=0$, for in this case $b=1$
and thus $\tilde Y_{\mathbf i,\mathbf c}=\bigl\{[t^0]\bigr\}$.

Now let $\mathbf i\in I^l$ and $b\in\mathbf B(-\infty)_{\mathbf i}$.
We write $\mathbf i=(i_1,\ldots,i_l)$ and $\mathbf c_{\mathbf i}(b)
=(c_1,\ldots,c_l)$. We set $\mathbf i'=(i_2,\ldots,i_l)$ and
$b'=\tilde f_{i_1}^{c_1}b$. We will apply the induction hypothesis
to $\mathbf i'$ and $b'$.

We note that $\varphi_{i_1}(b')=0$ and that $\mathbf c_{\mathbf i'}
(b')=(c_2,\ldots,c_l)$. For $j\in\{1,\ldots,l\}$, we set
$e_j=-\sum_{k=j+1}^lc_k\langle\alpha_{i_j},\alpha_{i_k}^\vee\rangle$.
We set $Z=\Xi(t_0\otimes b)$ and $Z'=\Xi(t_0\otimes b')$; then
$Z=\tilde e_{i_1}^{c_1}(Z')$, for $\Xi$ is an isomorphism of crystals.
The equality $\varphi_{i_1}(b')=0$ implies that
$$\varepsilon_{i_1}(Z')=\varepsilon_{i_1}(t_0\otimes b')=
\varepsilon_{i_1}(b')=-\langle\alpha_{i_1},\wt(b')\rangle=e_1.$$
Thanks to~(\ref{eq:DefPhiBFG}), the equality $\varphi_{i_1}(b')=0$
also leads to
$$\mu_+(Z')=s_i\mu_+\bigl(\overline{s_i}^{-1}Z'\bigr).$$

Proposition~\ref{pr:CrystOpAlg}~\ref{it:PrCOAb} thus asserts that the
map $(p,z)\mapsto y_{i_1}(pt^{e_1})z$ is a homeomorphism from
$\mathbb C[t^{-1}]_{c_1}^*\times\bigl(Z'\cap S_{\mu_+(Z')}^+
\cap S_{\mu_-(Z')}^-\bigr)$ onto an open and dense subset of
$Z\cap S_{\mu_+(Z)}^+\cap S_{\mu_-(Z)}^-$. Theorem~\ref{th:ParamOpen}
then follows immediately by induction.
\end{proof}

\subsection{The subsets $\tilde Y_{\mathbf i,\mathbf c}$}
\label{ss:TildeYic}
Given a sequence $\mathbf i=(i_1,\ldots,i_l)$ of elements of $I$
and a sequence $\mathbf p=(p_1,\ldots,p_l)$ of elements of $\mathscr K$,
we form the element
$$y_{\mathbf i}(\mathbf p)=y_{i_1}(p_1)\cdots y_{i_l}(p_l).$$
Given the sequence $\mathbf i$ as above and a sequence
$\mathbf c=(c_1,\ldots,c_l)$ of integers, we set
$$\tilde Y_{\mathbf i,\mathbf c}=\bigl\{[y_{\mathbf i}(\mathbf p)]\bigm|
\mathbf p\in(\mathscr K^\times)^l\text{ such that }\val(p_j)=\tilde c_j
\bigr\},$$
where $\tilde c_j=-c_j-\sum_{k=j+1}^lc_k\langle\alpha_{i_j},
\alpha_{i_k}^\vee\rangle$.

\begin{proposition}
\label{pr:TildeYic}
\begin{enumerate}
\item\label{it:PrTYa}
Let $\mathbf i\in I^l$, let $b\in\mathbf B(-\infty)_{\mathbf i}$ and
set $\mathbf c=\mathbf c_{\mathbf i}(b)$. Then the MV cycle
$\Xi(t_0\otimes b)$ is the closure of $\tilde Y_{\mathbf i,\mathbf c}$.
\item\label{it:PrTYb}
Let $\mathbf i=(i_1,\ldots,i_N)$ be a reduced decomposition of $w_0$
and let $\mathbf c=(c_1,\ldots,c_N)$ be an element in $\mathbb Z^N$. Let
$Z$ be the closure of $\tilde Y_{\mathbf i,\mathbf c}$ and let $\lambda$
be the coweight $c_1\alpha_{i_1}^\vee+\cdots+c_N\alpha_{i_N}^\vee$.
Then $Z$ is an MV cycle, $\mu_-(Z)=0$ and $\mu_+(Z)\geqslant\lambda$.
Moreover $\mu_+(Z)=\lambda$ if and only if $\mathbf c\in\mathcal
C_{\mathbf i}$.
\end{enumerate}
\end{proposition}

Many assertions of this proposition follow easily from Proposition
\ref{pr:CrystOpAlg} and Theorem~\ref{th:ParamOpen}. The truly new points
are the inequality $\mu_+(Z)\geqslant\lambda$ in Statement~\ref{it:PrTYb}
and the fact that the equality $\mu_+(Z)=\lambda$ holds only
if $\mathbf c\in\mathcal C_{\mathbf i}$. We will ground our proof on
the notion of $\mathbf i$-trail in Berenstein and Zelevinsky's
work~\cite{BerensteinZelevinsky01}. We first recall what it is about.

We denote the differential at $0$ of the one-parameter subgroups
$x_{\alpha_i}$ and $x_{-\alpha_i}$ by $E_i$ and $F_i$, respectively;
they are elements of the Lie algebra of $G$. Let $\mathbf i=(i_1,\ldots,
i_N)$ be a reduced decomposition of $w_0$, let $\gamma$ and $\delta$ two
weights in $X$, let $V$ be a rational $G$-module, and write
$V=\bigoplus_{\eta\in X}V_\eta$ for its decomposition in weight subspaces.
According to Definition~2.1 in \cite{BerensteinZelevinsky01}, an
$\mathbf i$-trail from $\gamma$ to $\delta$ in $V$ is a sequence of
weights $\pi=(\gamma=\gamma_0,\gamma_1,\ldots,\gamma_N=\delta)$ such that
each difference $\gamma_{j-1}-\gamma_j$ has the form $n_j\alpha_{i_j}$
for some non-negative integer $n_j$, and such that $E_{i_1}^{n_1}\cdots
E_{i_N}^{n_N}$ defines a non-zero map from $V_\delta$ to $V_\gamma$. To
such an $\mathbf i$-trail $\pi$, Berenstein and Zelevinsky associate the
sequence of integers $d_j(\pi)=\langle\gamma_{j-1}+\gamma_j,
\alpha_{i_j}^\vee\rangle/2$.

Assume moreover that $G$ is semisimple and simply connected. In that case,
$X$ is the free $\mathbb Z$-module with basis the set $\{\omega_i\mid
i\in I\}$ of fundamental weights. For each $i\in I$, we can thus speak of
the simple rational $G$-module with highest weight $\omega_i$, which we
denote by $V(\omega_i)$. Then by Theorem~3.10 in
\cite{BerensteinZelevinsky01}, the string cone $\mathcal C_{\mathbf i}$
is the set of all $(c_1,\ldots,c_N)\in\mathbb Z^N$ such that
$\sum_jd_j(\pi)c_j\geqslant0$ for any $i\in I$ and any $\mathbf i$-trail
$\pi$ from $\omega_i$ to $w_0s_i\omega_i$ in $V(\omega_i)$.

The following lemma explains why $\mathbf i$-trails are relevant to
our problem.
\begin{lemma}
\label{le:Mu&Trails}
Let $\mathbf i$, $\mathbf c$, $Z$ and $\lambda$ be as in the statement
of Proposition~\ref{pr:TildeYic}~\ref{it:PrTYb}, let $i\in I$, and assume
that $G$ is semisimple and simply connected. Then $\langle\omega_i,
\lambda-\mu_+(Z)\rangle$ is the minimum of the numbers $\sum_jd_j(\pi)c_j$
for all weights $\delta\in X$ and all $\mathbf i$-trails $\pi$ from
$\omega_i$ to $\delta$ in $V(\omega_i)$.
\end{lemma}
\begin{proof}
Let us consider an $\mathbf i$-trail $\pi=(\gamma_0,\gamma_1,\ldots,
\gamma_N)$ in $V(\omega_i)$ which starts from $\gamma_0=\omega_i$.
Introducing the integers $n_j$ such that $\gamma_{j-1}-\gamma_j=n_j
\alpha_{i_j}$, we obtain $\gamma_j=\omega_i-\sum_{k=1}^jn_k\alpha_{i_k}$
for each $j\in\{1,\ldots,N\}$ and so
$$d_j(\pi)=\langle\omega_i,\alpha_{i_j}^\vee\rangle-\sum_{k=1}^{j-1}
n_k\langle\alpha_{i_k},\alpha_{i_j}^\vee\rangle-n_j.$$
We then compute
$$\sum_{j=1}^Nd_j(\pi)c_j-\langle\omega_i,\lambda\rangle=\sum_{j=1}^N
\left(-n_j-\sum_{k=1}^{j-1}\langle\alpha_{i_k},\alpha_{i_j}^\vee\rangle
n_k\right)c_j=n_1\tilde c_1+\cdots+n_N\tilde c_N,$$
where we set as usual $\tilde c_j=-c_j-\sum_{k=j+1}^Nc_k\langle
\alpha_{i_j},\alpha_{i_k}^\vee\rangle$ for each $j\in\{1,\ldots,N\}$.

We adopt the notational conventions set up before Proposition
\ref{pr:KamnitzForm}. In particular, we embed $V(\omega_i)$ inside
$V(\omega_i)\otimes_{\mathbb C}\mathscr K$ and we view this latter
as a representation of the group $G(\mathscr K)$. We also consider
a non-degenerate contravariant bilinear form $(?,?)$ on $V(\omega_i)$;
it is compatible with the decomposition of $V(\omega_i)$ as the sum
of its weight subspaces and it satisfies $(v,E_iv')=(F_iv,v')$ for any
$i\in I$ and any vectors $v$ and $v'$ in $V(\omega_i)$. We extend the
contravariant bilinear form to $V(\omega_i)\otimes_{\mathbb C}\mathscr K$
by multilinearity.

By Proposition~\ref{pr:DefMu}, $\langle\omega_i,\mu_+(Z)\rangle$ is
the maximum of $\langle\omega_i,\nu\rangle$ for those $\nu\in\Lambda$
such that $S_\nu^+$ meets $\tilde Y_{\mathbf c,\mathbf i}$. Using
Proposition~\ref{pr:KamnitzForm}~\ref{it:PrKFb}, we deduce that
\begin{align*}
\langle\omega_i,\mu_+(Z)\rangle
&=\max\Bigl\{-\val\bigl(g^{-1}\cdot v_{\omega_i}\bigr)\Bigm|
g\in G(\mathscr K)\text{ such that }[g]\in\tilde Y_{\mathbf
c,\mathbf i}\Bigr\}\\[3pt]
&=\max\Biggl\{-\val\bigl(\bigl(v,y_{\mathbf i}(\mathbf p)^{-1}\cdot
v_{\omega_i}\bigr)\bigr)\Biggm|\begin{aligned}v\in V(\omega_i),&\
\mathbf p\in(\mathscr K^\times)^N\\\text{such that }&\val(p_j)=\tilde
c_j\end{aligned}\Biggr\},
\end{align*}
where we wrote $\mathbf p=(p_1,\ldots,p_N)$ as usual. Moreover we may
ask that the vector $v$ in the last maximum is a weight vector.

Let us denote by $M$ the minimum of the numbers $\sum_jd_j(\pi)c_j$ for
all $\mathbf i$-trails $\pi$ in $V(\omega_i)$ which start from $\omega_i$.
We expand the product
$$y_{\mathbf i}(\mathbf p)^{-1}=\exp(-p_NF_{i_N})\cdots\exp(-p_1F_{i_1})=
\sum_{n_1,\ldots,n_N\geqslant0}\frac{(-1)^{n_1+\cdots+n_N}\;p_1^{n_1}
\cdots p_N^{n_N}}{n_1!\cdots n_N!}\;F_{i_N}^{n_N}\cdots F_{i_1}^{n_1}$$
and we substitute in $\bigl(v,y_{\mathbf i}(\mathbf p)^{-1}\cdot
v_{\omega_i}\bigr)$: we get a sum of terms of the form
$$\frac{(-1)^{n_1+\cdots+n_N}\;p_1^{n_1}\cdots p_N^{n_N}}{n_1!\cdots n_N!}
\;\Bigl(v,F_{i_N}^{n_N}\cdots F_{i_1}^{n_1}\cdot v_{\omega_i}\Bigr).$$
If such a term is not zero, then the sequence
$$\pi=(\omega_i,\ \omega_i-n_1\alpha_{i_1},\ \omega_i-n_1\alpha_{i_1}-n_2
\alpha_{i_2},\ \ldots,\ \omega_i-n_1\alpha_{i_1}-\cdots-n_N\alpha_{i_N})$$
is an $\mathbf i$-trail and the term has valuation
$$n_1\tilde c_1+\cdots+n_N\tilde c_N=\sum_{j=1}^Nd_j(\pi)c_j-\langle
\omega_i,\lambda\rangle\geqslant M-\langle\omega_i,\lambda\rangle.$$
Therefore the valuation of $(v,y_{\mathbf i}(\mathbf p)^{-1}\cdot
v_{\omega_i})$ is greater or equal to $M-\langle\omega_i,\lambda\rangle$
for any $v\in V(\omega_i)$; we conclude that $\langle\omega_i,\mu_+(Z)
\rangle\leqslant\langle\omega_i,\lambda\rangle-M$.

Conversely, let $\pi$ be an $\mathbf i$-trail in $V(\omega_i)$ which
starts from $\omega_i$ and which is such that $\sum_jd_j(\pi)c_j=M$. With
this $\mathbf i$-trail come the numbers $n_1$, \dots, $n_N$ as before. By
definition of an $\mathbf i$-trail, there is then a weight vector
$v\in V(\omega_i)$ such that
$$\Bigl(v,F_{i_N}^{n_N}\cdots F_{i_1}^{n_1}\cdot v_{\omega_i}\Bigr)\neq0.$$
Given $(a_1,\ldots,a_N)\in(\mathbb C^\times)^N$, we set $\mathbf p=
(a_1t^{\tilde c_1},\ldots,a_Nt^{\tilde c_N})$ and look at the coefficient
$f$ of $t^{M-\langle\omega_i,\lambda\rangle}$ in $\bigl(v,y_{\mathbf i}
(\mathbf p)^{-1}\cdot v_{\omega_i}\bigr)$. The computation above shows
that $f$ is a polynomial in $(a_1,\ldots,a_N)$; it is not zero since the
coefficient of $a_1^{n_1}\cdots a_N^{n_N}$ in $f$ is
$$\frac{(-1)^{n_1+\cdots+n_N}}{n_1!\cdots n_N!}\Bigl(v,F_{i_N}^{n_N}
\cdots F_{i_1}^{n_1}\cdot v_{\omega_i}\Bigr)\neq0.$$
Therefore there exists $\mathbf p\in(\mathscr K^\times)^N$ with $\val(p_j)
=\tilde c_j$ such that $\bigl(v,y_{\mathbf i}(\mathbf p)^{-1}\cdot
v_{\omega_i}\bigr)$ has valuation $\leqslant M-\langle\omega_i,\lambda
\rangle$. It follows that $\langle\omega_i,\mu_+(Z)\rangle\geqslant
\langle\omega_i,\lambda\rangle-M$, which completes the proof.
\end{proof}

\trivlist
\item[\hskip\labelsep{\itshape Proof of Proposition~\ref{pr:TildeYic}.}]
\upshape
Statement~\ref{it:PrTYa} is established in the same fashion as Theorem
\ref{th:ParamOpen}, using Proposition~\ref{pr:CrystOpAlg}~\ref{it:PrCOAa}
instead of Proposition~\ref{pr:CrystOpAlg}~\ref{it:PrCOAb}.

Now let $\mathbf i$, $\mathbf c$, $Z$ and $\lambda$ as in the statement
of Statement~\ref{it:PrTYb}. Applying repeatedly Proposition
\ref{pr:CrystOpAlg}~\ref{it:PrTYa}, one shows easily that $Z$ is an MV
cycle. Furthermore by its very definition, $\tilde Y_{\mathbf i,\mathbf c}$
is contained in $S_0^-$; this entails that $\mu_-(Z)=0$.

If $\mathbf c$ is the string in direction $\mathbf i$ of an element
$b\in\mathbf B(-\infty)$, then $Z=\Xi(t_0\otimes b)$, and thus
$$\mu_+(Z)=\wt(Z)=\wt(t_0\otimes b)=\wt(b)=\wt\bigl(\tilde
e_{i_1}^{c_1}\cdots\tilde e_{i_N}^{c_N}1\bigr)=\lambda.$$
The equality $\mu_+(Z)=\lambda$ holds therefore for each
$\mathbf c\in\mathcal C_{\mathbf i}$.

It remains to show that $\mu_+(Z)\geqslant\lambda$ with equality
only if $\mathbf c\in\mathcal C_{\mathbf i}$. Let us first consider
the case where $G$ is semisimple and simply connected. Then
$\Lambda=\mathbb Z\Phi^\vee$ and we can speak of the fundamental weights
$\omega_i$ and of the $G$-modules $V(\omega_i)$.

Let $i\in I$. The sequence
$$\pi=(\omega_i,\ s_{i_1}\omega_i,\ s_{i_2}s_{i_1}\omega_i,\ \ldots,\
w_0\omega_i)$$
is an $\mathbf i$-trail in $V(\omega_i)$ for which $d_j(\pi)=0$ for
each $j$. By Lemma~\ref{le:Mu&Trails}, we deduce
$$\langle\omega_i,\lambda-\mu_+(Z)\rangle\leqslant\sum_jd_j(\pi)c_j=0.$$
This is enough to guarantee that $\mu_+(Z)\geqslant\lambda$.

Suppose now that $\mu_+(Z)=\lambda$. Lemma~\ref{le:Mu&Trails} implies
then that $\sum_jd_j(\pi)c_j\geqslant0$ for all $i\in I$, all weights
$\delta\in X$, and all $\mathbf i$-trails $\pi$ from $\omega_i$ to
$\delta$ in $V(\omega_i)$. In particular, this holds for all $i\in I$
and all $\mathbf i$-trails $\pi$ from $\omega_i$ to $w_0s_i\omega_i$ in
$V(\omega_i)$. By Theorem~3.10 in~\cite{BerensteinZelevinsky01}, this
implies $\mathbf c\in\mathcal C_{\mathbf i}$. The proof is thus complete
in the case where $G$ is semisimple and simply connected.

In the general case, we note that the inclusion $\mathbb Z\Phi^\vee
\hookrightarrow\Lambda$ defines an epimorphism from a semisimple simply
connected group $\dot G$ onto $G$, such that $\mathbb Z\Phi^\vee$
is the cocharacter group of a maximal torus $\dot T$ of $\dot G$ and
$\Phi$ is the root system of $(\dot G,\dot T)$. The morphism from $\dot
G$ to $G$ then induces a homeomorphism from $\dot{\mathscr G}$ onto the
neutral connected component of $\mathscr G$. The subsets $\tilde
Y_{\mathbf i,\mathbf c}$ of $\dot{\mathscr G}$ and $\mathscr G$ match
under this homeomorphism, as do the functions $\mu_\pm$. Since
Proposition~\ref{pr:TildeYic} holds for $\dot G$, it holds for $G$ as well.
\nobreak\noindent$\square$
\endtrivlist

\subsection{Lusztig's algebraic-geometric parametrization of $\mathbf B$}
\label{ss:LuszAlgGeomParam}
As we have seen in Section~\ref{ss:MVCycString}, the choice of a
reduced decomposition $\mathbf i$ of $w_0$ determines a bijection
$\mathbf c_{\mathbf i}:\mathbf B(-\infty)\to\mathcal C_{\mathbf i}$,
called the ``string parametrization''. The decomposition $\mathbf i$
also determines a bijection $b_{\mathbf i}:\mathbb N^N\to\mathbf
B(-\infty)$, called the ``Lusztig parametrization'', which reflects
Lusztig's original construction~\cite{Lusztig90} of the canonical
basis on a combinatorial level. We refer the reader to~\cite{Lusztig92},
\cite{Saito94} and Section~3.1 in~\cite{BerensteinZelevinsky01} for
additional information on the map $b_{\mathbf i}$ and its construction.

The Lusztig parametrizations $b_{\mathbf i}$ are convenient because
they permit a study of $\mathbf B(-\infty)$ by way of numerical data
in a fixed domain $\mathbb N^N$, but they are not intrinsic, for they
depend on the choice of $\mathbf i$. To avoid this drawback, Lusztig
introduces in~\cite{Lusztig96} a parametrization of $\mathbf B(-\infty)$
in terms of closed subvarieties in arc spaces on $U^-$. We will first
recall briefly his construction and then we will explain a relationship
with MV cycles. For simplicity, Lusztig restricts himself to the case
where $G$ is simply laced, but he explains in the introduction of
\cite{Lusztig96} that his results hold in the general case as well.

Lusztig starts by recalling a general construction. To a complex
algebraic variety $X$ and a non-negative integer $s$, one can associate
the space $X_s$ of all jets of curves drawn on $X$, of order $s$. In
formulas, one looks at the algebra $\mathbb C_s=\mathbb C[[t]]/(t^{s+1})$
and defines $X_s$ as the set of morphisms from $\Spec\mathbb C_s$ to $X$.
If $X$ is smooth of dimension $n$, then $X_s$ is smooth of dimension
$(s+1)n$. There exist morphisms of truncation
$$\cdots\to X_{s+1}\to X_s\to\cdots\to X_1\to X_0=X;$$
the projective limit of this inverse system of maps is the space
$X(\mathscr O)$. Finally the assignment $X\rightsquigarrow X_s$ is
functorial, hence $X_s$ is a group as soon as $X$ is one.

Now let $\mathbf i$ be a reduced decomposition of $w_0$. The morphism
$$y_{\mathbf i}:(a_1,\ldots,a_N)\mapsto y_{i_1}(a_1)\cdots y_{i_N}(a_N)$$
from $(\mathbb C)^N$ to $U^-$ gives by functoriality a morphism
$(y_{\mathbf i})_s:(\mathbb C_s)^N\to(U^-)_s$. Given an element
$\mathbf d=(d_1,\ldots,d_N)$ in $\mathbb N^N$, we may look at the image
of the subset
$$(t^{d_1}\mathbb C_s)\times\cdots\times(t^{d_N}\mathbb C_s)\subseteq
(\mathbb C_s)^N$$
by $(y_{\mathbf i})_s$. This is a constructible, irreducible subset of
$(U^-)_s$. If $s$ is big enough, then the closure of this subset depends
only on $b=b_{\mathbf i}(\mathbf d)$ and not on $\mathbf i$ and $\mathbf
d$ individually. (This is Lemma~5.2 of~\cite{Lusztig96}; the precise
condition is that $s$ must be strictly larger than $\height(\wt b)$.)
One may therefore denote this closure by $V_{b,s}$; it is a closed
irreducible subset of $(U^-)_s$ of codimension $\height(\wt b)$.
Proposition~7.5 in \cite{Lusztig96} asserts that moreover the assignment
$b\mapsto V_{b,s}$ is injective for $s$ big enough: there is a constant
$M$ depending only on the root system $\Phi$ such that
$$\Bigl(V_{b,s}=V_{b',s}\quad\text{and}\quad s>M\height(\wt b)\Bigr)
\quad\Longrightarrow\quad b=b'$$
for any $b$, $b'\in\mathbf B(-\infty)$. Thus $b\mapsto V_{b,s}$ may be
seen as a parametrization of $\mathbf B(-\infty)$ by closed irreducible
subvarieties of $(U^-)_s$.

Our next result shows that Lusztig's construction is related to MV cycles
and to Braverman, Finkelberg and Gaitsgory's theorem. We fix a dominant
coweight $\lambda\in\Lambda_{++}$. By Proposition~\ref{pr:CanTruncate},
the map $x\mapsto x\cdot[t^{w_0\lambda}]$ from $G(\mathscr O)$ to
$\mathscr G$ factorizes through the reduction map $G(\mathscr O)\to
G_s$ when $s$ is big enough, defining thus a map
$$\Upsilon_s:\ G_s\to\mathscr G,\ x\mapsto x\cdot[t^{w_0\lambda}].$$
On the other hand, we may consider the two embeddings of crystals
$\kappa_\lambda:\mathbf B(\lambda)\hookrightarrow\mathbf B(\infty)
\otimes\mathbf T_\lambda$ and $\iota_{w_0\lambda}:\mathbf B(\lambda)
\hookrightarrow\mathbf T_{w_0\lambda}\otimes\mathbf B(-\infty)$,
as in Section~\ref{ss:Crystals}. Finally, the isomorphism $\mathbf
B(\infty)^\vee\cong\mathbf B(-\infty)$ yields a bijection $b\mapsto
b^\vee$ from $\mathbf B(\infty)$ onto $\mathbf B(-\infty)$.
\begin{proposition}
\label{pr:CompLuszBFG}
We adopt the notations above and assume that $s$ is big enough so that
the map $\Upsilon_s$ exists and that the closed subsets $V_{b^\vee,s}$ are
defined for each $b\otimes t_\lambda$ in the image of $\kappa_\lambda$.
Then the diagram
$$\xymatrix{\mathbf B(\lambda)\ar@{_{(}->}[d]_{\iota_{w_0\lambda}}
\ar@{^{(}->}[r]^<>(.5){\kappa_\lambda}&\im(\kappa_\lambda)
\ar[d]^{b\otimes t_\lambda\mapsto\overline{\Upsilon_s(V_{b^\vee,s})}}\\
\mathbf T_{w_0\lambda}\otimes\mathbf B(-\infty)\ar[r]_<>(.5)\Xi
&\mathscr Z}$$
commutes.
\end{proposition}
\begin{proof}
This is a consequence of Proposition~\ref{pr:TildeYic}~\ref{it:PrTYa},
combined with a result of Morier-Genoud \cite{Morier-Genoud03}. We
first look at the commutative diagram that defines the embedding
$\iota_{w_0\lambda}$, namely
$$\xymatrix{&\mathbf B(\lambda)\ar@{_{(}->}[dl]_{\kappa_\lambda}
\ar[r]^<>(.5)\simeq\ar@{_{(}->}[rd]_<>(.5){\iota_{w_0\lambda}}&\mathbf
B(-w_0\lambda)^\vee\ar@{_{(}->}[d]&\mathbf B(-w_0\lambda)\ar@{-->}[l]
\ar@{_{(}->}[d]^{\kappa_{-w_0\lambda}}\\\mathbf B(\infty)\otimes\mathbf
T_\lambda&&\mathbf T_{w_0\lambda}\otimes\mathbf B(-\infty)&\mathbf
B(\infty)\otimes\mathbf T_{-w_0\lambda}.\ar@{-->}[l]}$$
The two arrows in broken line on this diagram are the maps
$b\mapsto b^\vee$; they are not morphisms of crystals. The map from
$\mathbf B(-w_0\lambda)$ to $\mathbf B(\lambda)$ obtained by composing
the two arrows on the top line intertwines the raising operators
$\tilde e_i$ with their lowering counterparts $\tilde f_i$ and sends
the highest weight element of $\mathbf B(-w_0\lambda)$ to the lowest
weight element of $\mathbf B (\lambda)$; it therefore coincides with
the map denoted by $\Phi_{-w_0\lambda}$ in~\cite{Morier-Genoud03}.

Now let $b\in\mathbf B(\lambda)$. We write $\kappa_\lambda(b)=b'\otimes
t_\lambda$ and $\kappa_{-w_0\lambda}(\Phi_{-w_0\lambda}^{-1}(b))=
b''\otimes t_{-w_0\lambda}$; thus $\iota_{w_0\lambda}(b)=
t_{w_0\lambda}\otimes(b'')^\vee$. We choose a reduced decomposition
$\mathbf i$ of $w_0$ and we set $\mathbf c=(c_1,\ldots,c_N)=
\mathbf c_{\mathbf i}^{}((b'')^\vee)$ and $(d_1,\ldots,d_N)=
b_{\mathbf i}^{-1}((b')^\vee)$. We additionally set $\tilde c_j=-c_j
-\sum_{k=j+1}^lc_k\langle\alpha_{i_j},\alpha_{i_k}^\vee\rangle$ for
each $j\in\{1,\ldots,N\}$.\vspace{2pt} Corollary~3.5 in~\cite{Morier-Genoud03}
then asserts that $d_j=\langle\alpha_{i_j},-w_0\lambda\rangle+\tilde c_j$
for all $j$. Now comparing the definition of Lusztig's subset
$V_{(b')^\vee,s}$ with the definition of $\tilde Y_{\mathbf i,\mathbf c}$
and using Proposition~\ref{pr:TildeYic}~\ref{it:PrTYa}, we compute
$$\overline{V_{(b')^\vee,s}\cdot[t^{w_0\lambda}]}=
\overline{t^{w_0\lambda}\cdot\tilde Y_{\mathbf i,\mathbf c}}=
t^{w_0\lambda}\cdot\Xi\bigl(t_0\otimes(b'')^\vee\bigr)=\Xi\bigl(
t_{w_0\lambda}\otimes(b'')^\vee\bigr)=(\Xi\circ\iota_{w_0\lambda})(b).$$
\end{proof}

\subsection{Link with Kamnitzer's construction}
\label{ss:LinkKamnCons}
Let $b\in\mathbf B(-\infty)$ and let $\mathbf i$ be a reduced
decomposition of $w_0$. Theorem~\ref{th:ParamOpen} explains how to
construct an open and dense subset in the MV cycle $\Xi(t_0\otimes b)$
when one knows the string parameter $\mathbf c_{\mathbf i}(b)$.
In his work on MV polytopes, Kamnitzer~\cite{Kamnitzer05a} presents
a similar result, which provides a dense subset of $\Xi(t_0\otimes
b)$ from the datum of the Lusztig parameter $b_{\mathbf i}^{-1}(b)$.
These two results are twin; indeed Kamnitzer's result and Proposition
\ref{pr:TildeYic}~\ref{it:PrTYa} can be quickly derived one from the
other. This section, which does not contain any formalized statement,
aims at explaining how.

Our main tool here is Berenstein, Fomin and Zelevinsky's work. In a
series of papers (among which~\cite{BerensteinFominZelevinsky96,
BerensteinZelevinsky97,BerensteinZelevinsky01}), these three authors
devise an elegant method that yields all transition maps between the
different parametrizations of $\mathbf B(-\infty)$ we have met, namely
the maps
$$b_{\mathbf j}^{-1}\circ b_{\mathbf i}^{}:
\mathbb N^N\to\mathbb N^N,\quad
\mathbf c_{\mathbf j}^{}\circ b_{\mathbf i}^{}:
\mathbb N^N\to\mathcal C_{\mathbf j},\quad
b_{\mathbf j}^{-1}\circ\mathbf c_{\mathbf i}^{-1}:
\mathcal C_{\mathbf i}\to\mathbb N^N,\quad
\mathbf c_{\mathbf j}^{}\circ\mathbf c_{\mathbf i}^{-1}:
\mathcal C_{\mathbf i}\to\mathcal C_{\mathbf j},$$
where $\mathbf i$ and $\mathbf j$ are two reduced decomposition of $w_0$.
In recalling their results hereafter, we will slightly modify their
notation; our modifications simplify the presentation, perhaps at the
price of the loss of positivity results.

We first alter the string parameter $\mathbf c_{\mathbf i}$ by defining
a map $\tilde{\mathbf c}_{\mathbf i}$ from $\mathbf B(-\infty)$ to
$\mathbb Z^N$ as follows: an element $b\in\mathbf B(-\infty)$ with string
parameter $\mathbf c_{\mathbf i}(b)=(c_1,\ldots,c_N)$ in direction
$\mathbf i$ is sent to the $N$-tuple $(\tilde c_1,\ldots,\tilde
c_N)$, where $\tilde c_j=-c_j-\sum_{k=j+1}^Nc_k\langle\alpha_{i_j},
\alpha_{i_k}^\vee\rangle$. We denote the image of this map
$\tilde{\mathbf c}_{\mathbf i}$ by $\tilde{\mathcal C}_{\mathbf i}$.

Let $\mathbf i=(i_1,\ldots,i_l)$ be a sequence of elements of $I$ and
let $\mathbf a=(a_1,\ldots,a_l)$ be a sequence of elements of $\mathbb
C^\times$. Assuming that the product $s_{i_1}\cdots s_{i_l}$ is a
reduced decomposition of an element $w\in W$, Theorem~1.2 in
\cite{BerensteinZelevinsky97} implies there is a unique element
in $U^-\cap B^+y_{\mathbf i}(\mathbf a)\,\overline w^{-1}$; we
denote it by $z_{\mathbf i}(\mathbf a)$. Theorem~1.2 in
\cite{BerensteinZelevinsky97} also asserts that if $\mathbf i$ is
a reduced decomposition of $w_0$, then the map $z_{\mathbf i}$ is a
birational morphism from $(\mathbb C^\times)^N$ to $U^-$. Now under the
same assumption, the map $y_{\mathbf i}$ is a birational morphism from
$\mathbb C^N$ to $U^-$. If $\mathbf i$ and $\mathbf j$ are both reduced
decompositions of $w_0$, we therefore get birational maps
\begin{equation}
\label{eq:GeomLiftings}
z_{\mathbf j}^{-1}\circ z_{\mathbf i}^{},\quad
y_{\mathbf j}^{-1}\circ z_{\mathbf i}^{},\quad
z_{\mathbf j}^{-1}\circ y_{\mathbf i}^{}\quad\text{and}\quad
y_{\mathbf j}^{-1}\circ y_{\mathbf i}^{}
\end{equation}
from $\mathbb C^N$ to itself. After extension of the base field, we may
view them as birational maps from $\mathscr K^N$ to itself.

We need now to define the process of tropicalization. Here we depart from
Berenstein, Fomin and Zelevinsky's purely algebraic method based on total
positivity and semifields and adopt a more pedestrian approach.

Let $k$ and $l$ be two positive integers and let $\mathbf f:\mathscr K^k
\to\mathscr K^l$ be a rational map, represented as a sequence $(f_1,\ldots,
f_l)$ of rational functions in $k$ indeterminates. These indeterminates
are collectively denoted as a sequence $\mathbf p=(p_1,\ldots,p_k)$. We
suppose that no component $f_j$ vanishes identically. Now choose
$j\in\{1,\ldots,l\}$ and $\mathbf m=(m_1,\ldots,m_k)\in\mathbb Z^k$.
There exists a non-empty (Zariski) open subset $\Omega\subseteq(\mathbb
C^\times)^k$ such that the valuation of $f_j(a_1t^{m_1},\ldots,a_kt^{m_k})$
is a constant $\hat f_j$, independent on the point
$\mathbf a=(a_1,\ldots,a_k)$ in $\Omega$. (It is here implicitely
understood that if $\mathbf a\in\Omega$, then neither the numerator nor
the denominator of the rational function $f_j$ vanishes after
substitution.) The term of lowest degree in $f_j(a_1t^{m_1},\ldots,
a_kt^{m_k})$ may then be written $\bar f_j(\mathbf a)t^{\hat f_j}$,
where $\bar f_j$ is a rational function with complex coefficients in
the indeterminates $a_1$, \dots, $a_k$. Of course, $\hat f_j$ and
$\bar f_j$ depend on the choice of $\mathbf m\in\mathbb Z^k$, but
the open subset $\Omega$ may be chosen to meet the demand simultaneously
for all $\mathbf m$. Indeed, as we make the substitution $p_i=a_it^{m_i}$,
each monomial in the indeterminates $p_1$, \dots, $p_k$ in the numerator
or in the denominator of $f_j$ becomes a non-zero element of $\mathscr K$.
To find the term $\bar f_j(\mathbf a)t^{\hat f_j}$ of lowest degree in
$f_j(a_1t^{m_1},\ldots,a_kt^{m_k})$, we collect the monomials in the
numerator of $f_j$ that get minimal valuation, and likewise in the
denominator. The r\^ole of the condition $\mathbf a\in\Omega$ is to ensure
that no accidental cancellation occurs when we make the sum of these
monomials, in the numerator as well as in the denominator. Since there are
only finitely many monomials, there are only finitely many possibilities
for accidental cancellations, hence finitely many conditions on $\mathbf
a$ to be prescribed by $\Omega$. Moreover monomials in the numerator or
in the denominator of $f_j$ are selected or discarded according to their
valuation, and we can divide $\mathbb R^k$ into finitely many regions,
say $\mathbb R^k=D^{(1)}\sqcup\cdots\sqcup D^{(t)}$, so that the set of
selected monomials depends only on the domain $D^{(r)}$ to which $\mathbf
m$ belongs. Since the valuation of each monomial depends affinely on
$\mathbf m$, the regions $D^{(1)}$, \dots, $D^{(t)}$ are indeed
intersections of affine hyperplanes and open affine half-spaces, hence
are locally closed, convex and polyhedral. For the same reason, $\hat f_j$
depends affinely on $\mathbf m$ in each region $D^{(r)}$; for its part,
$\bar f_j$ remains constant when $\mathbf m$ varies inside a region
$D^{(r)}$. Finally we note that the choice of the domain $\Omega\subseteq
(\mathbb C^\times)^k$, the decomposition $\mathbb R^k=D^{(1)}\sqcup\cdots
\sqcup D^{(t)}$ and the reduction $f_j\mapsto(\hat f_j,\bar f_j)$ may be
carried out for all $j\in\{1,\ldots,l\}$ at the same time. In particular
each $\mathbf m\in\mathbb Z^k$ yields a tuple $\hat{\mathbf f}=(\hat f_1,
\ldots,\hat f_l)$ of integers and a rational map $\bar{\mathbf f}=
(\bar f_1,\ldots,\bar f_l)$ from $\mathbb C^k$ to $\mathbb C^l$.
We summarize these observations in a formalized statement:\\[5pt]
\textit{Let $\mathbf f:\mathscr K^k\to\mathscr K^l$ be a rational map,
without identically vanishing component. Then there exists a partition
$\mathbb R^k=D^{(1)}\sqcup\cdots\sqcup D^{(t)}$ of $\mathbb R^k$ into
finitely many locally closed polyhedral convex subsets, there exist
affine maps $\hat{\mathbf f}^{(1)},\ldots,\hat{\mathbf f}^{(t)}:\mathbb
R^k\to\mathbb R^l$, there exist rational maps $\bar{\mathbf f}^{(1)},
\ldots,\bar{\mathbf f}^{(t)}:\mathbb C^k\to\mathbb C^l$, and there exists
an open subset $\Omega\subseteq(\mathbb C^\times)^k$ with the following
property: for each $r\in\{1,\ldots,t\}$, each lattice point $\mathbf m$
in $D^{(r)}\cap\mathbb Z^k$, each point $\mathbf a\in\Omega$, and each
sequence $\mathbf p\in(\mathscr K^\times)^k$ such that the lower degree
term of $p_i$ is $a_it^{m_i}$, the map $\mathbf f$ has a well-defined
value in $(\mathscr K^\times)^l$ at $\mathbf p$, the map $\bar{\mathbf
f}^{(r)}$ has a well-defined value in $(\mathbb C^\times)^l$ at
$\mathbf a$, and the term of lower degree of $f_j(\mathbf p)$ has
valuation $\hat f^{(r)}_j(\mathbf m)$ and coefficient
$\bar f^{(r)}_j(\mathbf a)$.}
\vspace{5pt}

We define the tropicalization of $\mathbf f$ as the map $\mathbf f^\trop:
\mathbb R^k\to\mathbb R^l$ whose restriction to each $D^{(r)}$ coincides
with the restriction of the corresponding $\hat{\mathbf f}^{(r)}$; this
is a continuous piecewise affine map. If the rational map $\mathbf f$
we started with has complex coefficients (that is, if it comes from
a rational map from $\mathbb C^k$ to $\mathbb C^l$ by extension of the
base field), then the convex subsets $D^{(r)}$ are cones and the affine
maps $\hat{\mathbf f}^{(r)}$ are linear.

With this notation and this terminology, Theorems~5.2 and 5.7 in
\cite{BerensteinZelevinsky01} implies that the maps
$$b_{\mathbf j}^{-1}\circ b_{\mathbf i}^{}:
\mathbb N^N\to\mathbb N^N,\quad
\tilde{\mathbf c}_{\mathbf j}^{}\circ b_{\mathbf i}^{}:
\mathbb N^N\to\tilde{\mathcal C}_{\mathbf j},\quad
b_{\mathbf j}^{-1}\circ\tilde{\mathbf c}_{\mathbf i}^{-1}:
\tilde{\mathcal C}_{\mathbf i}\to\mathbb N^N,\quad
\tilde{\mathbf c}_{\mathbf j}^{}\circ\tilde{\mathbf c}_{\mathbf i}^{-1}:
\tilde{\mathcal C}_{\mathbf i}\to\tilde{\mathcal C}_{\mathbf j}$$
are restrictions of the tropicalizations of the maps in
(\ref{eq:GeomLiftings}).

One may here observe a hidden symmetry. Using the equality
$\overline{w_0}^2=(-1)^{2\rho^\vee}$, where $2\rho^\vee$ is the
sum of all positive coroots in $\Phi_+^\vee$, one checks that the
birational maps $y_{\mathbf j}^{-1}\circ z_{\mathbf i}^{}$ and
$z_{\mathbf j}^{-1}\circ y_{\mathbf i}^{}$ are equal. These maps
have therefore the same tropicalization. In other words,
$\tilde{\mathbf c}_{\mathbf j}^{}\circ b_{\mathbf i}^{}$ and
$b_{\mathbf j}^{-1}\circ\tilde{\mathbf c}_{\mathbf i}^{-1}$ are
given by the same piecewise affine formulas. The sentence following
Theorem~3.8 in \cite{BerensteinZelevinsky01} seems to indicate that
this fact has escaped observation up to now.

In~\cite{Kamnitzer05a}, Kamnitzer introduces subsets $A^{\mathbf i}
(n_\bullet)$ in $\mathscr G$, where $\mathbf i$ is a reduced decomposition
of $w_0$ and $n_\bullet\in\mathbb N^N$. Combining Theorem~4.7 in
\cite{Kamnitzer05b} with the proof of Theorem~3.1 in~\cite{Kamnitzer05a},
one can see that $\Xi(t_0\otimes b_{\mathbf i}(n_\bullet))$ is the
closure of $A^{\mathbf i}(n_\bullet)$. On the other hand, Theorem~4.5
in~\cite{Kamnitzer05a} says that
$$A^{\mathbf i}(n_\bullet)=\bigl\{[z_{\mathbf i}(\mathbf q)]\bigm|
\mathbf q=(q_1,\ldots,q_N)\in(\mathscr K^\times)^N\text{ such that }
\val(q_j)=n_j\bigr\}.$$

Now fix $b\in\mathbf B(-\infty)$ and a reduced decomposition $\mathbf i$
of $w_0$. Call $\tilde{\mathbf c}=(\tilde c_1,\ldots,\tilde c_N)$ the
modified string parameter $\tilde{\mathbf c}_{\mathbf i}(b)$ of $b$ in
direction $\mathbf i$ and call $n_\bullet=(n_1,\ldots,n_N)$ the Lusztig
parameter $b_{\mathbf i}^{-1}(b)$ of $b$ with respect to $\mathbf i$. The
rational maps $\mathbf f=z_{\mathbf i}^{-1}\circ y_{\mathbf i}^{}$ and
$\mathbf g=y_{\mathbf i}^{-1}\circ z_{\mathbf i}^{}$ are mutually
inverse birational maps from $\mathscr K^N$ to itself, and by
Berenstein and Zelevinsky's theorem,
$$\mathbf f^\trop(\tilde{\mathbf c})=n_\bullet\quad\text{and}\quad
\mathbf g^\trop(n_\bullet)=\tilde{\mathbf c}.$$
The analysis that we made to define the tropicalization of $\mathbf f$
and $\mathbf g$ shows the existence of open subsets $\Omega$ and
$\Omega'$ of $(\mathbb C^\times)^N$ and of rational maps $\bar{\mathbf f}$
and $\bar{\mathbf g}$ from $\mathbb C^N$ to itself such that:
\begin{itemize}
\item For each $\mathbf a\in\Omega$ and $\mathbf b\in\Omega'$,
$\bar{\mathbf f}(\mathbf a)$ and $\bar{\mathbf g}(\mathbf b)$
have well-defined values in $(\mathbb C^\times)^N$.
\item For any $N$-tuple $\mathbf p$ of Laurent series whose terms
of lower degree are $a_1t^{\tilde c_1}$, \dots, $a_Nt^{\tilde c_N}$
with $(a_1,\ldots,a_N)\in\Omega$, the evaluation $\mathbf f(\mathbf p)$
is a well-defined element $\mathbf q$ of $(\mathscr K^\times)^N$;
moreover the lower degree terms of the components of $\mathbf q$ are
$\bar f_1(\mathbf a)t^{n_1}$, \dots, $\bar f_N(\mathbf a)t^{n_N}$.
\item For any $N$-tuple $\mathbf q$ of Laurent series whose terms
of lower degree are $b_1t^{n_1}$, \dots, $b_Nt^{n_N}$ with $(b_1,\ldots,
b_N)\in\Omega'$, the evaluation $\mathbf g(\mathbf q)$ is a well-defined
element $\mathbf p$ of $(\mathscr K^\times)^N$; moreover the lower
degree terms of the components of $\mathbf p$ are $\bar g_1(\mathbf b)
t^{\tilde c_1}$, \dots, $\bar g_N(\mathbf b)t^{\tilde c_N}$.
\end{itemize}
Because $\mathbf f$ and $\mathbf g$ are mutually inverse birational maps,
so are $\bar{\mathbf f}$ and $\bar{\mathbf g}$. One can then assume that
these two latter maps are mutually inverse isomorphisms between $\Omega$
and $\Omega'$, by shrinking these open subsets if necessary. Thus
$\mathbf f$ and $\mathbf g$ set up a bijective correspondence between
$$\widehat\Omega=\left\{\mathbf p\in(\mathscr K^\times)^N\Biggm|
\begin{aligned}\text{each $p_j$ has }&\text{lower degree term}\\
a_jt^{\tilde c_j}\text{ with }&(a_1,\ldots,a_N)\in\Omega\end{aligned}
\right\}$$
and
$$\widehat\Omega'=\left\{\mathbf q\in(\mathscr K^\times)^N\Biggm|
\begin{aligned}\text{each $q_j$ has }&\text{lower degree term}\\
b_jt^{n_j}\text{ with }&(b_1,\ldots,b_N)\in\Omega'\end{aligned}
\right\}.$$
In other words, to each $\mathbf p\in\widehat\Omega$ corresponds
a $\mathbf q\in\widehat\Omega'$ such that $y_{\mathbf i}(\mathbf p)
=z_{\mathbf i}(\mathbf q)$, and conversely. This shows the equality
$$\bigl\{[y_{\mathbf i}(\mathbf p)]\bigm|\mathbf p\in\widehat\Omega\bigr\}
=\bigl\{[z_{\mathbf i}(\mathbf q)]\bigm|\mathbf q\in\widehat\Omega'\bigr
\}.$$
By Kamnitzer's theorem, the right-hand side is dense in $A^{\mathbf i}
(n_\bullet)$ hence in $\Xi(t_0\otimes b)$. We thus get another proof
of our Proposition~\ref{pr:TildeYic}~\ref{it:PrTYa}, which claims that
$\Xi(t_0\otimes b)$ is the closure of the left-hand side.

\begin{other*}{Remark}
We fix a reduced decomposition $\mathbf i$ of $w_0$. Each MV cycle $Z$
such that $\mu_-(Z)=0$ is the closure of a set $\tilde Y_{\mathbf i,
\mathbf c}$ for a certain $\mathbf c\in\mathcal C_{\mathbf i}$; indeed
there exists $b\in\mathbf B(-\infty)$ such that $Z=\Xi(t_0\otimes b)$,
and one takes then $\mathbf c=\mathbf c_{\mathbf i}(b)$. It follows that
$S_0^-$ is contained in the union $\bigcup_{\mathbf c\in\mathcal
C_{\mathbf i}}\overline{\tilde Y_{\mathbf i,\mathbf c}}$. On the other
side, each $\tilde Y_{\mathbf i,\mathbf c}$ is contained in $S_0^-$.
One could then hope that $S_0^-$ is the disjoint union of the $\tilde
Y_{\mathbf i,\mathbf c}$ for $\mathbf c\in\mathcal C_{\mathbf i}$, because
the analogous property $S_0^-=\bigsqcup_{n_\bullet\in\mathbb N^N}A^{\mathbf
i}(n_\bullet)$ for Kamnitzer's subsets holds (see Proposition~4.1 in
\cite{Kamnitzer05a}).

This is alas not the case in general, as the following counter-example
shows. We take $G=\SL_4$ with its usual pinning and enumerate the simple
roots in the usual way $(\alpha_1,\alpha_2,\alpha_3)$. We choose the
reduced decomposition $\mathbf i=(2,1,3,2,1,3)$ and consider
$$g=y_2(-1)\;y_1(1/t)\;y_3(1/t)\;y_2(t)\;y_1(-1/t)\;y_3(-1/t)=
\begin{pmatrix}1&0&0&0\\0&1&0&0\\-1&t-1&1&0\\-1/t&1&0&1\end{pmatrix}.$$
If one tries to factorize an element in $gG(\mathscr O)\cap
U^-(\mathscr K)$ as a product
$$y_2(p_1)\;y_1(p_2)\;y_3(p_3)\;y_2(p_4)\;y_1(p_5)\;y_3(p_6)$$
using Berenstein, Fomin and Zelevinsky's method
\cite{BerensteinFominZelevinsky96}, and if after that one adjusts
$\mathbf c=(c_1,\ldots,c_6)$ so that $(\val(p_1),\ldots,\val(p_6))=
(\tilde c_1,\ldots,\tilde c_6)$, then one finds
$$c_1\leqslant0,\quad c_2\leqslant0,\quad c_3\leqslant0,\quad
c_4\geqslant1,\quad c_5\geqslant1,\quad c_6\geqslant1.$$
These conditions on $\mathbf c$ must be satisfied in order that $[g]$
can belong to $\tilde Y_{\mathbf i,\mathbf c}$. However the equations
that define the cone $\mathcal C_{\mathbf i}$ are
$$c_1\geqslant0,\quad c_2\geqslant c_6\geqslant0,\quad c_3\geqslant
c_5\geqslant0,\quad c_2+c_3\geqslant c_4\geqslant c_5+c_6.$$
We conclude that $[g]\not\in\bigcup_{\mathbf c\in\mathcal C_{\mathbf i}}
\tilde Y_{\mathbf i,\mathbf c}$.
\end{other*}

\section{BFG crystal operations on MV cycles and root operators on
LS galleries}
\label{se:BFG&GalMod}
Let $\lambda\in\Lambda_{++}$ be a dominant coweight. Littelmann's path
model~\cite{Littelmann95} affords a concrete realization of the crystal
$\mathbf B(\lambda)$ in terms of piecewise linear paths drawn on
$\Lambda\otimes_{\mathbb Z}\mathbb R$; it depends on the choice of a
path joining $0$ to $\lambda$ and contained in the dominant Weyl chamber.
In~\cite{GaussentLittelmann05}, Gaussent and Littelmann present a variation
of the path model, replacing piecewise linear paths by galleries in the
Coxeter complex of the affine Weyl group $W^\aff$. They define a set
$\Gamma^+_\LS(\gamma_\lambda)$ of ``LS~galleries'', which depends on the
choice of a minimal gallery $\gamma_\lambda$ joining $0$ to $\lambda$
and contained in the dominant Weyl chamber. Defining ``root operators''
$e_\alpha$ and $f_\alpha$ for each simple root $\alpha$ in $\Phi$, they
endow $\Gamma^+_\LS(\gamma_\lambda)$ with the structure of a crystal,
which happens to be isomorphic to $\mathbf B(\lambda)$. Using a
Bott-Samelson resolution\vspace{2pt} $\pi:\hat\Sigma(\gamma_\lambda)\to
\overline{\mathscr G_\lambda}$ and a Bia\l ynicki-Birula decomposition of
$\hat\Sigma(\gamma_\lambda)$ into a disjoint union of cells $C(\delta)$,
Gaussent and Littelmann associate a closed subvariety $Z(\delta)=
\overline{\pi(C(\delta))}$ of $\mathscr G$ to each LS gallery $\delta$ and
show that the map $Z$ is a bijection from $\Gamma^+_\LS(\gamma_\lambda)$
onto $\mathscr Z(\lambda)$.

The main result of this section is Theorem~\ref{th:CompLSBG}, which says
that $Z$ is an isomorphism of crystals. In other words, the root operators
on LS galleries match Braverman and Gaitsgory's crystal operations on MV
cycles under the bijection $Z$.

Strictly speaking, our proof for this comparison result is valid only
when $\lambda$ is regular. The advantage of this situation is that
elements in $\Gamma^+_\LS(\gamma_\lambda)$ are then galleries of alcoves.
In the case where $\lambda$ is singular, Gaussent and Littelmann's
constructions involve a more general class of galleries (see Section~4
in~\cite{GaussentLittelmann05}). Such a sophistication is however not
needed: our presentation of Gaussent and Littelmann's results in
Section~\ref{ss:GaCeMVCy} below makes sense even if $\lambda$ is
singular. Within this framework, our comparison theorem is valid for any
$\lambda$, regular or singular.

A key idea of Gaussent and Littelmann is to view the affine Grassmannian
as a subset of the set of vertices of the (affine) Bruhat-Tits building of
$G(\mathscr K)$. In Section~\ref{ss:CoxCo&BTBui}, we review quickly
basic facts about the latter and study the stabilizer in $U^+(\mathscr K)$
of certain of its faces. We warn here the reader that we use our own
convention pertaining the Bruhat-Tits building: indeed our Iwahori
subgroup is the preimage of $B^-$ by the specialization map at $t=0$
from $G(\mathscr O)$ to $G$, whereas Gaussent and Littelmann use the
preimage of $B^+$. Our convention is unusual, but it makes the statement
of our comparison result more natural. Section~\ref{ss:GaCeMVCy}
recalls the main steps in Gaussent and Littelmann's construction, in a
way that encompasses the peculiarities of the case where $\lambda$ is
singular. The final Section~\ref{ss:CompTh} contains the proof of
our comparison theorem. To prove the equality $\tilde e_iZ(\delta)=
Z(e_{\alpha_i}\delta)$ for each LS gallery $\delta$ and each $i\in I$,
we use the criterion of Proposition~\ref{pr:CritOpCrys}. The first three
conditions are easily checked, while the inclusion $Z(\delta)\subseteq
Z(e_{\alpha_i}\delta)$ is established in Proposition~\ref{pr:InclMVCyc}.

\subsection{Affine roots, the Coxeter complex and the Bruhat-Tits building}
\label{ss:CoxCo&BTBui}
We consider the vector space $\Lambda_{\mathbb R}=\Lambda
\otimes_{\mathbb Z}\mathbb R$. We define a real root of the affine root
system (for short, an affine root) as a pair $(\alpha,n)\in\Phi\times
\mathbb Z$. To an affine root $(\alpha,n)$, we associate:
\begin{itemize}
\item the reflection $s_{\alpha,n}:x\mapsto x-\bigl(\langle\alpha,
x\rangle-n\bigr)\,\alpha^\vee$ of $\Lambda_{\mathbb R}$;
\item the affine hyperplane $H_{\alpha,n}=\{x\in\Lambda_{\mathbb R}
\mid\langle\alpha,x\rangle=n\}$ of fixed points of $s_{\alpha,n}$;
\item the closed half-space $H^-_{\alpha,n}=\{x\in\Lambda_{\mathbb R}
\mid\langle\alpha,x\rangle\leqslant n\}$;
\item the one-parameter additive subgroup $x_{\alpha,n}:b\mapsto
x_\alpha(bt^n)$ of $G(\mathscr K)$; here $b$ belongs to either
$\mathbb C$ or $\mathscr K$.
\end{itemize}
We denote the set of all affine roots by $\Phi^\aff$. We embed $\Phi$
in $\Phi^\aff$ by identifying a root $\alpha\in\Phi$ with the affine
root $(\alpha,0)$. We choose an element $0$ that does not belong to
$I$; we set $I^\aff=I\sqcup\{0\}$ and $\alpha_0=(-\theta,-1)$, where
$\theta$ is the highest root of $\Phi$. The elements $\alpha_i$ with
$i\in I^\aff$ are called simple affine roots.

The group of affine transformations of $\Lambda_{\mathbb R}$ generated by
all reflections $s_{\alpha,n}$ is called the affine Weyl group and is
denoted by $W^\aff$. For each $i\in I^\aff$, we set $s_i=s_{\alpha_i}$.
Then $W^\aff$ is a Coxeter system when equipped with the set of generators
$\{s_i\mid i\in I^\aff\}$. The parabolic subgroup of $W^\aff$
generated by the simple reflections $s_i$ with $i\in I$ is
isomorphic to $W$. For each $\lambda\in\mathbb Z\Phi^\vee$, the
translation $\tau_\lambda:x\mapsto x+\lambda$ belongs to $W^\aff$. All
these translations form a normal subgroup in $W^\aff$, isomorphic to the
coroot lattice $\mathbb Z\Phi^\vee$, and $W^\aff$ is the semidirect
product $W^\aff=\mathbb Z\Phi^\vee\rtimes W$.

The group $W^\aff$ acts on the set $\Phi^\aff$ of affine roots: one
demands that $w(H^-_\beta)=H^-_{w\beta}$ for each element
$w\in W^\aff$ and each affine root $\beta\in\Phi^\aff$. The action
of an element $w\in W$ or a translation $\tau_\lambda$ on an
affine root $(\alpha,n)\in\Phi\times\mathbb Z$ is given by
$w\,(\alpha,n)=(w\alpha,n)$ or $\tau_\lambda\,(\alpha,n)=
\bigl(\alpha,n+\langle\alpha,\lambda\rangle\bigr)$. One checks
that $ws_\alpha w^{-1}=s_{w\alpha}$ for all $w\in W^\aff$ and
$\alpha\in\Phi^\aff$. Using Equation~(\ref{eq:EqCommut1}), one
checks that
\begin{equation}
\label{eq:EqCommut6}
(t^\lambda\,\overline w\,)\,x_\alpha(a)\,(t^\lambda\,\overline
w\,)^{-1}=x_{\tau_\lambda w(\alpha)}(\pm a)
\end{equation}
in $G(\mathscr K)$, for all $\lambda\in\mathbb Z\Phi^\vee$, $w\in W$,
$\alpha\in\Phi^\aff$ and $a\in\mathscr K$.

We denote by $\mathfrak H$ the arrangement formed by the hyperplanes
$H_\beta$, where $\beta\in\Phi^\aff$. It divides the vector space
$\Lambda_{\mathbb R}$ into faces. Faces with maximal dimension are
called alcoves; they are the connected components of $\Lambda_{\mathbb R}
\setminus\bigcup_{H\in\mathfrak H}H$. Faces of codimension~$1$ are called
facets; faces of dimension~$0$ are called vertices. The closure of a face
is the disjoint union of faces of smaller dimension. Endowed with the set
of all faces, $\Lambda_{\mathbb R}$ becomes a polysimplicial complex,
called the Coxeter complex $\mathscr A^\aff$; it is endowed with an
action of $W^\aff$.

The dominant open Weyl chamber is the subset
$$C_\dom=\{x\in\Lambda_{\mathbb R}\mid\forall i\in I,\ \langle
\alpha_i,x\rangle>0\}.$$
The fundamental alcove
$$A_\fund=\{x\in C_\dom\mid\langle\theta,x\rangle<1\}$$
is the complement\vspace{2pt} of $\bigcup_{i\in I^\aff}H^-_{\alpha_i}$.
We label the faces contained in $\overline{A_\fund}$ by proper subsets
of $I^\aff$ by setting
$$\phi_J=\left(\bigcap_{i\in J\vphantom\setminus}H_{\alpha_i}\right)
\setminus\left(\bigcup_{i\in I^\aff\setminus J}H^-_{\alpha_i}\right)$$
for each $J\subset I^\aff$. For instance $\phi_\varnothing$ is the
alcove $A_\fund$ and $\phi_I$ is the vertex $\{0\}$. Any face of our
arrangement $\mathfrak H$ is conjugated under the action of $W^\aff$
to exactly one face contained in $\overline{A_\fund}$, because this
latter is a fundamental domain for the action of $W^\aff$ on
$\Lambda_{\mathbb R}$. We say that a subset $J\subset I^\aff$ is
the type of a face $F$ if $F$ is conjugated to $\phi_J$ under $W^\aff$.

We denote by $\hat B$ the (Iwahori) subgroup of $G(\mathscr K)$
generated by the torus $T(\mathscr O)$ and by the elements $x_\alpha(ta)$
and $x_{-\alpha}(a)$, where $\alpha\in\Phi_+$ and $a\in\mathscr O$. In
other words, $\hat B$ is the preimage of the Borel subgroup $B^-$
under the specialization map at $t=0$ from $G(\mathscr O)$ to $G$. We
lift the simple reflections $s_i$ to the group $G(\mathscr K)$ by setting
$$\overline{s_i}=x_{\alpha_i}(1)x_{-\alpha_i}(-1)x_{\alpha_i}(1)
=x_{-\alpha_i}(-1)x_{\alpha_i}(1)x_{-\alpha_i}(-1)$$
for each $i\in I^\aff$. We lift any element $w\in W^\aff$ to an
element $\overline w\in G(\mathscr K)$ so that $\overline w=
\overline{s_{i_1}}\cdots\overline{s_{i_l}}$ for each reduced decomposition
$s_{i_1}\cdots s_{i_l}$ of $w$. This notation does not conflict with
our earlier notation $\overline{s_i}$ for $i\in I$ and $\overline w$
for $w\in W$. For each $\lambda\in\mathbb Z\Phi^\vee$, the lift
$\overline{\tau_\lambda}$ of the translation $\tau_\lambda$ coincides
with $t^\lambda$ up to a sign (that is, up to the multiplication by an
element of the form $(-1)^\mu$ with $\mu\in\mathbb Z\Phi^\vee$).

The affine Bruhat-Tits building $\mathscr I^\aff$ is a polysimplicial
complex endowed with an action of $G(\mathscr K)$. The affine Coxeter
complex $\mathscr A^\aff$ can be embedded in $\mathscr I^\aff$ as the
subcomplex formed by the faces fixed by $T$; in this identification,
the action of an element $w\in W^\aff$ on $\mathscr A^\aff$ matches
the action of $\overline w$ on $(\mathscr I^\aff)^T$. Each face of
$\mathscr I^\aff$ is conjugated under the action of $G(\mathscr K)$
to exactly one face contained in $\overline{A_\fund}$; we say that a
subset $J\subset I^\aff$ is the type of a face $F$ if $F$ is
conjugated to $\phi_J$. Finally there is a $G(\mathscr K)$-equivariant
map of the affine Grassmannian $\mathscr G$ into $\mathscr I^\aff$,
which extends the map $[t^\lambda]\mapsto\{\lambda\}$ from $\mathscr G^T$
into $\mathscr A^\aff\cong(\mathscr I^\aff)^T$.

Given a subset $J\subseteq I^\aff$, we denote by $\hat P_J$ the
subgroup of $G(\mathscr K)$ generated by $\hat B$ and the elements
$\overline{s_i}$ for $i\in J$; thus $\hat B=\hat P_\varnothing$ and
$G(\mathscr O)=\hat P_I$. (The subgroup $\hat P_J$ is the stabilizer in
$G(\mathscr K)$ of the face $\phi_J$. For each $g\in G(\mathscr K)$,
the stabilizer of the face $g\phi_J$ is thus the parahoric subgroup
$g\hat P_Jg^{-1}$. This bijection between the set of faces in the
affine building and the set of parahoric subgroups in $G(\mathscr K)$
is indeed the starting point for the definition of the building,
see~\S2.1 in~\cite{BruhatTits72}.) To shorten the notation, we will
write $\hat P_i$ instead of $\hat P_{\{i\}}$ for each $i\in I^\aff$.
Similarly, for each $i\in I^\aff$, we will write $W_i$ to indicate
the subgroup $\{1,s_i\}$ of $W^\aff$.

We denote the stabilizer in $U^+(\mathscr K)$ of a face $F$ of the
affine building by $\Stab_+(F)$. Our last task in this section is
to determine as precisely as possible the group $\Stab_+(F)$ and the
set $\Stab_+(F')/\Stab_+(F)$ when $F$ and $F'$ are faces of the Coxeter
complex such that $F'\subseteq\overline F$. We need additional notation
for that. Given a real number $a$, we denote the smallest integer greater
than $a$ by $\lceil a\rceil$. To a face $F$ of the Coxeter complex,
Bruhat and Tits (see (7.1.1) in \cite{BruhatTits72}) associate the
function $f_F:\alpha\mapsto\sup_{x\in F}\langle\alpha,x\rangle$ on
the dual space of $\Lambda_{\mathbb R}$. If $\alpha\in\Phi$, then $\lceil
f_F(\alpha)\rceil$ is the smallest integer $n$ such that $F$ lies in
the closed half-space $H^-_{\alpha,n}$. The function $f_F$ is convex
and positively homogeneous of degree $1$; in particular,
$f_F(i\alpha+j\beta)\leqslant if_F(\alpha)+jf_F(\beta)$ for all roots
$\alpha,\beta\in\Phi$ and all positive integers $i,j$. When $F$ and
$F'$ are two faces of the Coxeter complex such that
$F'\subseteq\overline F$, we denote by $\Phi_+^\aff(F',F)$ the set of
all affine roots $\beta\in\Phi_+\times\mathbb Z$ such that
$F'\subseteq H_\beta$ and $F\not\subseteq H^-_\beta$; in other words,
$(\alpha,n)\in\Phi_+^\aff(F',F)$ if and only if $\alpha\in\Phi_+$,
$n=f_{F'}(\alpha)$ and $n+1=\lceil f_F(\alpha)\rceil$. We denote by
$\Stab_+(F',F)$ the subgroup of $U^+(\mathscr K)$ generated by the
elements of the form $x_\beta(a)$ with $\beta\in\Phi_+^\aff(F',F)$
and $a\in\mathbb C$.

\begin{proposition}
\label{pr:DescStab+}
\begin{enumerate}
\item\label{it:PrDS+a}
The stabilizer $\Stab_+(F)$ of a face $F$ of the Coxeter
complex is generated by the elements $x_\alpha(p)$, where
$\alpha\in\Phi_+$ and $p\in\mathscr K$ satisfy $\val(p)\geqslant
f_F(\alpha)$.
\item\label{it:PrDS+b}
Let $F$ and $F'$ be two faces of the Coxeter complex such that
$F'\subseteq\overline F$. Then $\Stab_+(F',F)$ is a set of
representatives for the right cosets of $\Stab_+(F)$ in $\Stab_+(F')$.
For any total order on the set $\Phi_+^\aff(F',F)$, the map
$$(a_\beta)_{\beta\in\Phi_+^\aff(F',F)}\mapsto\prod_{\beta\in
\Phi_+^\aff(F',F)}x_\beta(a_\beta)$$
is a bijection from $\mathbb C^{\Phi_+^\aff(F',F)}$ onto $\Stab_+(F',F)$.
\end{enumerate}
\end{proposition}
\begin{proof}
Item~\ref{it:PrDS+a} is proved in Bruhat and Tits's paper
\cite{BruhatTits72}, see in particular Sections~(7.4.4) and
Equation~(1) in Section~(7.1.8). We note here that this
result implies that for any total order on $\Phi_+$, the map
$$(p_\alpha)_{\alpha\in\Phi_+}\mapsto\prod_{\alpha\in\Phi_+}
x_\alpha\Bigl(p_\alpha t^{\lceil f_F(\alpha)\rceil}\Bigr)$$
is a bijection from $\mathscr O^{\Phi_+}$ onto $\Stab_+(F)$.

We now turn to Item~\ref{it:PrDS+b}. We first observe the
following property of $\Phi_+^\aff(F',F)$: for each pair $i,j$ of
positive integers and each pair $(\alpha,m),(\beta,n)$ of affine
roots in $\Phi_+^\aff(F',F)$ such that $i\alpha+j\beta\in\Phi$, the
affine root $(i\alpha+j\beta,im+jn)$ belongs to $\Phi_+^\aff(F',F)$.
Indeed $F'\subseteq H_{\alpha,m}\cap H_{\beta,n}$ implies
$F'\subseteq H_{i\alpha+j\beta,im+jn}$, and the inequality
$$f_F(i\alpha+j\beta)\geqslant if_F(\alpha)-jf_F(-\beta)=
if_F(\alpha)+jn>im+jn$$
shows that $F\not\subseteq H^-_{i\alpha+j\beta,im+jn}$.
Standard arguments based on Chevalley's commutator
formula~(\ref{eq:EqCommut5}) show then the second assertion in
Item~\ref{it:PrDS+b}.

Now the map $(\alpha,m)\mapsto\alpha$ from $\Phi^\aff$ to $\Phi$ restricts
to a bijection from $\Phi_+^\aff(F',F)$ onto a subset $\Phi'_+$ of
$\Phi_+$. We set $\Phi''_+=\Phi_+\setminus\Phi'_+$. We endow $\Phi_+$
with a total order, chosen so that every element in $\Phi'_+$ is smaller
than every element in $\Phi''_+$, and we transport the order induced
on $\Phi'_+$ to $\Phi_+^\aff(F',F)$. By Item~\ref{it:PrDS+a}, each
element in $\Stab_+(F')$ may be uniquely written as a product
\begin{equation}
\label{eq:PfPrDS+}
\prod_{\alpha\in\Phi_+}x_\alpha\Bigl(p_\alpha t^{\lceil f_{F'}(\alpha)
\rceil}\Bigr)
\end{equation}
with $(p_\alpha)_{\alpha\in\Phi_+}$ in $\mathscr O^{\Phi_+}$. We write
$p_\alpha=a_\alpha+tq_\alpha$ for each $\alpha\in\Phi'_+$, with
$a_\alpha\in\mathbb C$ and $q_\alpha\in\mathscr O$. Thus for each
$(\alpha,m)\in\Phi_+^\aff(F',F)$, we have $p_\alpha t^{\lceil f_{F'}
(\alpha)\rceil}=a_\alpha t^m+q_\alpha t^{\lceil f_F(\alpha)\rceil}$.
On the other hand, $\lceil f_{F'}(\alpha)\rceil=\lceil f_F(\alpha)\rceil$
for each $\alpha\in\Phi''_+$. We may therefore rewrite the product in
(\ref{eq:PfPrDS+}) as
$$\left(\prod_{(\alpha,m)\in\smash{\Phi_+^\aff}(F',F)}x_\alpha
\Bigl(a_\alpha t^m\Bigr)\;x_\alpha\Bigl(q_\alpha t^{\lceil f_F(\alpha)
\rceil}\Bigr)\right)\left(\prod_{\alpha\in\Phi''_+}x_\alpha
\Bigl(p_\alpha t^{\lceil f_F(\alpha)\rceil}\Bigr)\right).$$
We rearrange the first product above using again Chevalley's commutator
formula: there exists a family $(r_\alpha)_{\alpha\in\Phi'_+}$ of power
series such that this product is
$$\left(\prod_{(\alpha,m)\in\smash{\Phi_+^\aff}(F',F)}
x_\alpha\Bigl(a_\alpha t^m\Bigr)\right)\left(\prod_{\alpha\in\Phi'_+}
x_\alpha\Bigl(r_\alpha t^{\lceil f_F(\alpha)\rceil}\Bigr)\right),$$
and for fixed numbers $a_\alpha$, the map $(q_\alpha)\mapsto
(r_\alpha)$ is a bijection from $\mathcal O^{\Phi'_+}$ onto itself. We
conclude that the map
$$((a_\beta),(p_\alpha))\mapsto\left(\prod_{\beta\in\smash{\Phi_+^\aff}
(F',F)}x_\beta\Bigl(a_\beta\Bigr)\right)\left(\prod_{\alpha\in\Phi_+}
x_\alpha\Bigl(p_\alpha t^{\lceil f_F(\alpha)\rceil}\Bigr)\right)$$
is a bijection from $\mathbb C^{\Phi_+^\aff(F',F)}\times\mathscr
O^{\Phi_+}$ onto $\Stab_+(F')$. This means exactly that the map
$(g,h)\mapsto gh$ is a bijection from $\Stab_+(F',F)\times\Stab_+(F)$
onto $\Stab_+(F')$. The proof of Item~\ref{it:PrDS+b} is now complete.
\end{proof}

Things are more easy to grasp when $F$ is an alcove and $F'$ is a
facet of $\overline F$, because then $\Phi_+^\aff(F',F)$ has at most
one element. In this particular case, certain commutators involving
elements of $\Stab_+(F')$ and $\Stab_+(F)$ automatically belong to
$\Stab_+(F)$.
\begin{lemma}
\label{le:CommStab+}
Let $F$ be an alcove of the Coxeter complex and let $F'$ be a facet of
$\overline F$. Let $(\alpha,m)\in\Phi_+\times\mathbb Z$ be the affine
root such that $F'$ lies in the wall $H_{\alpha,m}$ and let $(\beta,n)
\in\Phi^\aff$ be such that $F\subseteq H^-_{\beta,n}$. We assume that
$\beta$ is either positive or is the opposite of a simple root, and that
$\beta\neq-\alpha$. Then for each $q\in\mathscr O$ and each $v\in
\Stab_+(F',F)$, the commutator $x_{\beta,n}(q)\,v\,x_{\beta,n}(q)^{-1}
\,v^{-1}$ belongs to $\Stab_+(F)$.
\end{lemma}
\begin{proof}
There is nothing to show if $F\subseteq H^-_{\alpha,m}$ since $v=1$
in this case. We may thus assume that $\Stab_+(F',F)=\bigl\{(\alpha,m)
\bigr\}$; then there is an $a\in\mathbb C$ such that
$v=x_{\alpha,m}(a)$.

Suppose first that $\beta=\alpha$. Then
$$x_{\beta,n}(q)\,v\,x_{\beta,n}(q)^{-1}\,v^{-1}=
x_{\beta,n}(q)\,x_{\alpha,m}(a)\,x_{\beta,n}(-q)\,x_{\alpha,m}(-a)=
x_\alpha(qt^n+at^m-qt^n-at^m)=1.$$
Therefore the assertion holds in this case.

Suppose now that $\beta\neq\alpha$. The facet $F'$ is contained in
the closure of exactly two alcoves, $F$ and say $F^*$, the
latter lying in $H^-_{\alpha,m}$. Then $f_{F^*}(\alpha)=m$. We observe
that no wall other than $H_{\alpha,m}$ separates $F^*$ and $F$.
In particular, $H_{\beta,n}$ does not separate $F^*$ and $F$, because
$\beta\neq\pm\alpha$. Since $F$ lies in $H^-_{\beta,n}$, so does
$F^*$, and thus $f_{F^*}(\beta)\leqslant n$. Therefore for any pair
of positive integers $i,j$ such that $i\alpha+j\beta$ is a root,
$f_{F^*}(i\alpha+j\beta)\leqslant im+jn$. This means that $F^*$ lies
in the half-space $H^-_{i\alpha+j\beta,im+jn}$. Again, the wall
$H_{i\alpha+j\beta,im+jn}$ does not separate $F^*$ and $F$,
and we conclude that $F$ lies in the half-space
$H^-_{i\alpha+j\beta,im+jn}$. Chevalley's commutator formula
(\ref{eq:EqCommut5}) implies that
\begin{align*}
x_{\beta,n}(q)\,v\,x_{\beta,n}(q)^{-1}\,v^{-1}&=x_{\beta,n}(q)\,
x_{\alpha,m}(a)\,x_{\beta,n}(-q)\,x_{\alpha,m}(-a)\\[5pt]
&=\prod_{i,j>0}x_{i\alpha+j\beta,im+jn}\bigl(C_{i,j,\alpha,\beta}
\,a^i(-q)^j\bigr).
\end{align*}
Here the product is taken over all pairs of positive integers $i,j$
such that $i\alpha+j\beta$ is a root. The assumption about $\beta$ in
the statement of the lemma implies that such a root $i\alpha+j\beta$ is
necessarily positive. By Proposition~\ref{pr:DescStab+}~\ref{it:PrDS+a},
each factor $x_{i\alpha+j\beta,im+jn}\bigl(C_{i,j,\alpha,\beta}
a^i(-q)^j\bigr)$ belongs to $\Stab_+(F)$. Thus the commutator
$x_{\beta,n}(q)\,v\,x_{\beta,n}(q)^{-1}\,v^{-1}$ belongs to $\Stab_+(F)$.
\end{proof}

\begin{other*}{Remark}
The first assertion in Proposition~\ref{pr:DescStab+}~\ref{it:PrDS+b}
means that $\Stab_+(F')$ has the structure of a bicrossed product
$\Stab_+(F',F)\Join\Stab_+(F)$ (see~\cite{Takeuchi81}) whenever
$F$ and $F'$ are two faces in the Coxeter complex such that
$F'\subseteq\overline F$. Suppose now that $F$ is an alcove and that
$F'$ is a facet of $\overline F$. Then Proposition~\ref{pr:DescStab+}
\ref{it:PrDS+a} and Lemma~\ref{le:CommStab+} imply that each element
$v\in\Stab_+(F',F)$ normalizes the group $\Stab_+(F)$. Thus $\Stab_+(F)$
is a normal subgroup of $\Stab_+(F')$ and $\Stab_+(F')$ is the
semidirect product $\Stab_+(F',F)\ltimes\Stab_+(F)$.
\end{other*}

\subsection{Galleries, cells and MV cycles}
\label{ss:GaCeMVCy}
We fix a dominant coweight $\lambda\in\Lambda_{++}$. As usual, we denote
by $P_\lambda$ the standard parabolic subgroup $P_J$ of $G$, where
$J=\{j\in I\mid\langle\alpha_j,\lambda\rangle=0\}$. Besides, we denote by
$\{\lambda_\fund\}$ the vertex in $\overline{A_\fund}$ with the same type
as $\{\lambda\}$. Finally, there is a unique element $w_\lambda$ in
$W^\aff$ with minimal length such that $\lambda=w_\lambda(\lambda_\fund)$.
Thus among all alcoves in $\mathscr A^\aff$ having $\{\lambda\}$ as
vertex, $w_\lambda(A_\fund)$ is the one closest to $A_\fund$.

We denote the length of $w_\lambda$ by $p$ and we choose a reduced
decomposition $s_{i_1}\cdots s_{i_p}$ of it, with $(i_1,\ldots,i_p)
\in(I^\aff)^p$. The geometric translation of this choice is the datum
of the sequence
$$\gamma_\lambda=\bigl(\{0\}\subset\overline{\Gamma^{}_0}\supset
\Gamma'_1\subset\overline{\Gamma^{}_1}\supset\cdots\supset\Gamma'_p
\subset\overline{\Gamma^{}_p}\supset\{\lambda\}\bigr)$$
of alcoves and facets (also known as a gallery) in $\mathscr A^\aff$, where
$$\Gamma^{}_j=s_{i_1}\cdots s_{i_j}(A_\fund)\quad\text{and}
\quad\Gamma'_j=s_{i_1}\cdots s_{i_{j-1}}\bigl(\phi_{\{i_j\}}\bigr).$$
By Proposition~2.19~(iv) in~\cite{Tits74}, these alcoves and facets are
all contained in the dominant Weyl chamber $C_\dom$. The choice of the
reduced decomposition $s_{i_1}\cdots s_{i_p}$ of $w_\lambda$ and the
notations $P_\lambda$, $\lambda_\fund$, $\gamma_\lambda$ will be kept
for the rest of Section~\ref{se:BFG&GalMod}.

We define the Bott-Samelson variety as the smooth projective variety
$$\hat\Sigma(\gamma_\lambda)=G(\mathscr O)\underset{\hat B}\times
\hat P_{i_1}\underset{\hat B}\times\cdots\underset{\hat B}\times
\hat P_{i_p}/\hat B.$$
We will denote the image in $\hat\Sigma(\gamma_\lambda)$ of an element
$(g_0,g_1,\ldots,g_p)\in G(\mathscr O)\times\hat P_{i_1}\times\cdots
\times\hat P_{i_p}$ by the usual notation $[g_0,g_1,\ldots,g_p]$.
The group $G(\mathscr O)$ acts on $\hat\Sigma(\gamma_\lambda)$ by
left multiplication on the first factor. There is a
$G(\mathscr O)$-equivariant map $\pi:[g_0,g_1,\ldots,g_p]\mapsto
g_0g_1\cdots g_p\bigl[t^{\lambda_\fund}\bigr]$ from
$\hat\Sigma(\gamma_\lambda)$ onto $\overline{\mathscr G_\lambda}$.

The geometric language of buildings is of great convenience in the
study of the Bott-Samelson variety. Indeed each element
$d=[g_0,g_1,\ldots,g_p]$ in $\hat\Sigma(\gamma_\lambda)$ may be
viewed as a gallery
\begin{equation}
\label{eq:Gallery}
\delta=\bigl(\{0\}=\Delta'_0\subset\overline{\Delta^{}_0}\supset
\Delta'_1\subset\overline{\Delta^{}_1}\supset\cdots\supset\Delta'_p
\subset\overline{\Delta^{}_p}\supset\Delta'_{p+1}\bigr)
\end{equation}
in $\mathscr I^\aff$, where
\begin{gather*}
\Delta^{}_j=g_0\cdots g_j(A_\fund)\quad\text{for $0\leqslant j\leqslant
p$,}\\\Delta'_j=g_0\cdots g_{j-1}\bigl(\phi_{\{i_j\}}\bigr)\quad\text{for
$1\leqslant j\leqslant p$,}\\\text{and}\quad\Delta'_{p+1}=g_0\cdots
g_p\{\lambda_\fund\}.
\end{gather*}
(This gallery has the same type as $\gamma_\lambda$, that is, each
facet $\Delta'_j$ of $\delta$ has the same type as the corresponding
element $\Gamma'_j$ in $\gamma_\lambda$. We also observe that the
vertex $\Delta'_{p+1}$ of the affine building corresponds to the
element $\pi(d)$ of the affine Grassmannian.) Thus for instance the
point $\bigl[1,\overline{s_{i_1}},\overline{s_{i_2}},\ldots,
\overline{s_{i_p}}\bigr]$ in $\hat\Sigma(\gamma_\lambda)$ is viewed as
the gallery $\gamma_\lambda$. With this picture in mind, one proves
easily the following proposition.

\begin{proposition}
\label{pr:PiAboveOrbit}
The restriction of $\pi$ to $\pi^{-1}(\mathscr G_\lambda)$ is a fiber
bundle with fiber isomorphic to~$P_\lambda/B^+$.
\end{proposition}
\begin{proof}
Let $J=\{j\in I\mid\langle\alpha_j,\lambda\rangle=0\}$ and let $P_J^-$
be the parabolic subgroup of $G$ generated by $B^-\cup M_J$. The set $S$
of alcoves whose closure contains $\phi_J$ is in canonical bijection with
the set of all Iwahori subgroups of $G(\mathscr K)$ contained in
$\hat P_J$, hence with $\hat P_J/\hat B\cong P_J^-/B^-$. In particular,
$P_J^-$ acts transitively on $S$ and $S$ is isomorphic to $P_\lambda/B^+$.

Now let $F=\pi^{-1}([t^\lambda])$ and let $H$ be the stabilizer of
$[t^\lambda]$ in $G(\mathscr O)$; thus $H\supseteq P_J^-$. Since
$\pi$ is $G(\mathscr O)$-equivariant, $H$ acts on $F$ and there is a
commutative diagram
$$\xymatrix@C+=30pt{G(\mathscr O)\times_HF\ar[d]\ar[r]^<>(.5)\simeq&
\pi^{-1}(\mathscr G_\lambda)\ar[d]_\pi\\G(\mathscr O)/H\ar[r]^<>(.5)
\simeq&\mathscr G_\lambda.}$$
It thus suffices to prove that $F$ is isomorphic to $S$.

Each element $d\in F$ can be viewed as a gallery
$$\delta=\bigl(\{0\}\subset\overline{\Delta^{}_0}\supset\Delta'_1
\subset\overline{\Delta^{}_1}\supset\cdots\supset\Delta'_p\subset
\overline{\Delta^{}_p}\supset\{\lambda\}\bigr)$$
in $\mathscr I^\aff$ stretching from $\{0\}$ to $\{\lambda\}$.
We claim that $\overline{\Delta_0}$ always contains $\phi_J$. When
all faces of $\delta$ belong to $\mathscr A^\aff$, this claim follows
from the proof of Proposition~2.29 in~\cite{Tits74} (with
$\proj_{\{0\}}\{\lambda\}=\phi_J$); the general case is obtained by
retracting $\delta$ onto $\mathscr A^\aff$ from the fundamental alcove,
see Lemma~3.6 in \cite{Tits74}.

We finally consider the map $f:d\mapsto\Delta_0$ from $F$ to $S$.
Corollary~3.4 in \cite{Tits74} implies that $f$ is injective, because
in any apartment, there is only one non-stammering gallery of the same
type as $\gamma_\lambda$ that starts from a given chamber $\Delta_0$.
On the other side, $f$ is $H$-equivariant; it is thus surjective, for
$P_J^-$ acts transitively on the codomain. We conclude that $f$ is an
isomorphism from $F$ onto $S$.
\end{proof}

This proposition implies the following equality, which we record
for later use:
\begin{equation}
\label{eq:DimVarBS}
\bigl|\Phi_+\bigr|+p=
\dim\hat\Sigma(\gamma_\lambda)=
\dim\mathscr G_\lambda+\dim(P_\lambda/B^+)=
\height(\lambda-w_0\lambda)+\dim(P_\lambda/B^+).
\end{equation}

Our next task is to obtain a Bia\l ynicki-Birula decomposition of
the Bott-Samelson variety. The torus $T$ acts on the latter by left
multiplication on the first factor. If we represent an element
$d\in\hat\Sigma(\gamma_\lambda)$ by a gallery $\delta$ as in
(\ref{eq:Gallery}), then $d$ is fixed by $T$ if and only if all the
faces $\Delta_j$ and $\Delta'_j$ are in the Coxeter complex $\mathscr
A^\aff\cong(\mathscr I^\aff)^T$. We devote a word to this situation:
a gallery $\delta$ as in~(\ref{eq:Gallery}), of the same type as
$\gamma_\lambda$, all of whose faces are in $\mathscr A^\aff$, is
called a combinatorial gallery. The weight $\nu$ such that
$\Delta'_{p+1}=\{\nu\}$ is called the weight of $\delta$; it belongs to
$\lambda+\mathbb Z\Phi^\vee$, because $\{\nu\}$ has the same type as
$\{\lambda\}$.

We denote the set of all combinatorial galleries by
$\Gamma(\gamma_\lambda)$. This set is in bijection with $W\times
W_{i_1}\times\cdots\times W_{i_p}$; indeed the map
$(\delta_0,\delta_1,\ldots,\delta_p)\mapsto\bigl[\,\overline{\delta_0},
\overline{\delta_1},\ldots,\overline{\delta_p}\,\bigr]$ from $W\times
W_{i_1}\times\cdots\times W_{i_p}$ to $\hat\Sigma(\gamma_\lambda)$ is
injective and its image is the set of $T$-fixed points in the codomain.
Concretely this correspondence maps $(\delta_0,\delta_1,\ldots,
\delta_p)\in W\times W_{i_1}\times\cdots\times W_{i_p}$ to the
combinatorial gallery whose faces are
\begin{equation}
\label{eq:SeqToFaces}
\Delta^{}_j=\delta_0\cdots\delta_j(A_\fund)\quad\text{and}
\quad\Delta'_j=\delta_0\cdots\delta_{j-1}\bigl(\phi_{\{i_j\}}\bigr)
\end{equation}
and whose weight is
\begin{equation}
\label{eq:RecTarget}
\nu=\delta_0\delta_1\cdots\delta_p\lambda_\fund.
\end{equation}

The retraction $r_\varnothing$ from $\mathscr G$ onto $\mathscr G^T
\cong\Lambda$ can be extended to a map of polysimplicial complexes
from $\mathscr I^\aff$ onto $(\mathscr I^\aff)^T\cong\mathscr A^\aff$.
Following Section~7 in~\cite{GaussentLittelmann05}, we further extend
this retraction to a map from $\hat\Sigma(\gamma_\lambda)$ onto
$\hat\Sigma(\gamma_\lambda)^T\cong\Gamma(\gamma_\lambda)$ by applying
it componentwise to galleries. The preimage by this map of a
combinatorial gallery $\delta$ will be denoted by $C(\delta)$.

Our aim now is to describe precisely the cell $C(\delta)$ associated
to a combinatorial gallery $\delta$. Representing the latter as in
(\ref{eq:Gallery}), we introduce the notation
$$\Stab_+(\delta)=\Stab_+(\Delta'_0,\Delta^{}_0)\times
\Stab_+(\Delta'_1,\Delta^{}_1)\times\cdots\times
\Stab_+(\Delta'_p,\Delta^{}_p).$$

\begin{proposition}
\label{pr:C(Delta)}
Let $\delta$ be a combinatorial gallery and let $(\delta_0,\delta_1,
\ldots,\delta_p)$ be the sequence in $W\times W_{i_1}\times\cdots\times
W_{i_p}$ associated to $\delta$ by Equations~(\ref{eq:SeqToFaces}).
Then the map
$$(v_0,v_1,\ldots,v_p)\mapsto\bigl[v_0\;\overline{\delta_0}\,,\;
\overline{\delta_0}^{-1}\,v_1\;\overline{\delta_0\delta_1}\,,\;
\overline{\delta_0\delta_1}^{-1}\,v_2\;\overline{\delta_0\delta_1
\delta_2}\,,\ldots,\;\overline{\delta_0\cdots\delta_{p-1}}^{-1}\,
v_p\;\overline{\delta_0\cdots\delta_p}\,\bigr]$$
from $\Stab_+(\delta)$ to $\hat\Sigma(\gamma_\lambda)$ is injective
and its image is $C(\delta)$.
\end{proposition}
\begin{proof}
Set
$$\widetilde{\Stab_+(\delta)}=\Stab_+(\Delta'_0)\underset{\Stab_+
(\Delta^{}_0)}\times\Stab_+(\Delta'_1)\underset{\Stab_+(\Delta^{}_1)}
\times\cdots\underset{\Stab_+(\Delta^{}_{p-1})}\times\Stab_+
(\Delta'_p)/\Stab_+(\Delta^{}_p).$$
From the inclusions
\begin{alignat*}2
\Stab_+(\Delta^{}_j)&\subseteq\overline{\delta_0\cdots\delta_j}\;\hat B\;
\overline{\delta_0\cdots\delta_j}^{-1}&&\text{(for
$0\leqslant j\leqslant p$),}\\[4pt]
\Stab_+(\Delta'_0)&\subseteq G(\mathscr O)\overline{\delta_0}^{-1},&&\\[4pt]
\Stab_+(\Delta'_j)&\subseteq\overline{\delta_0\cdots\delta_{j-1}}\;\hat
P_{i_j}\;\overline{\delta_0\cdots\delta_j}^{-1}&\qquad&\text{(for
$1\leqslant j\leqslant p$),}
\end{alignat*}
standard arguments imply that the map
$$f:[v_0,v_1,\ldots,v_p]\mapsto\bigl[v_0\;\overline{\delta_0}\,,\;
\overline{\delta_0}^{-1}\,v_1\;\overline{\delta_0\delta_1}\,,\;
\overline{\delta_0\delta_1}^{-1}\,v_2\;\overline{\delta_0\delta_1
\delta_2}\,,\ldots,\;\overline{\delta_0\cdots\delta_{p-1}}^{-1}\,
v_p\;\overline{\delta_0\cdots\delta_p}\,\bigr]$$
from $\widetilde{\Stab_+(\delta)}$ to $\hat\Sigma(\gamma_\lambda)$ is
well-defined.

The proof of Proposition~6 in~\cite{GaussentLittelmann05} says that
an element $d=[g_0,g_1,\ldots,g_p]$ in the Bott-Samelson variety
belongs to the cell $C(\delta)$ if and only if there exists
$u_0,u_1,\ldots,u_p\in U^+(\mathscr K)$ such that
$$g_0g_1\cdots g_jA_\fund\,=u_j\Delta^{}_j\quad\text{and}
\quad u_{j-1}\Delta'_j=u_j\Delta'_j$$
for each $j$. Setting $v_0=u_0$ and $v_j=u_{j-1}^{-1}u_j$ for
$1\leqslant j\leqslant p$, the conditions above can be rewritten
$$g_0g_1\cdots g_j\hat B=v_0v_1\cdots v_j\;\overline{\delta_0\delta_1
\cdots\delta_j}\;\hat B\quad\text{and}\quad v_j\in\Stab_+(\Delta'_j),$$
which shows that $f([v_0,v_1,\ldots,v_p])=d$. Therefore the image
of $f$ contains the cell $C(\delta)$. The reverse inclusion can be
established similarly.

The map $f$ is injective. Indeed suppose that two elements
$v=[v_0,v_1,\ldots,v_p]$ and $v'=[v'_0,v'_1,\ldots,v'_p]$ in
$\widetilde{\Stab_+(\delta)}$ have the same image. Then
$$v_0v_1\cdots v_j\;\overline{\delta_0\delta_1\cdots\delta_j}\;\hat B=
v'_0v'_1\cdots v'_j\;\overline{\delta_0\delta_1\cdots\delta_j}\;\hat B$$
for each $j\in\{0,\ldots,p\}$. This means geometrically that
$$v_0v_1\cdots v_j\,\overline{\delta_0\delta_1\cdots\delta_j}\,
A_\fund\,=v'_0v'_1\cdots v'_j\;\overline{\delta_0\delta_1
\cdots\delta_j}\,A_\fund;$$
in other words, $v_0v_1\cdots v_j$ and $v'_0v'_1\cdots v'_j$ are equal
in $U^+(\mathscr K)/\Stab_+(\Delta^{}_j)$. Since this holds for each $j$,
the two elements $v$ and $v'$ are equal in $\widetilde{\Stab_+(\delta)}$.
We conclude that $f$ induces a bijection from
$\widetilde{\Stab_+(\delta)}$ onto $C(\delta)$.

It then remains to observe that the map
$(v_0,v_1,\ldots,v_p)\mapsto[v_0,v_1,\ldots,v_p]$ from
$\Stab_+(\delta)$ to $\widetilde{\Stab_+(\delta)}$ is bijective. This
follows from Proposition~\ref{pr:DescStab+}~\ref{it:PrDS+b}: indeed
for each $[a_0,a_1,\ldots,a_p]\in\widetilde{\Stab_+(\delta)}$, the
element $(v_0,v_1,\ldots,v_p)\in\Stab_+(\delta)$ such that
$[v_0,v_1,\ldots,v_p]=[a_0,a_1,\ldots,a_p]$ is uniquely determined by
the condition that for all $j\in\{0,1,\ldots,p\}$,
$$v_j\in\bigl((v_0\cdots v_{j-1})^{-1}(a_0\cdots a_j)\Stab_+(\Delta_j)
\bigr)\cap\Stab_+(\Delta'_j,\Delta_j).$$
\end{proof}

The definition of the map $\pi$, Equation~(\ref{eq:RecTarget}),
Proposition~\ref{pr:DescStab+}~\ref{it:PrDS+b} and
Proposition~\ref{pr:C(Delta)} yield the following explicit description
of the image of the cell $C(\delta)$ by the map $\pi$.
\begin{corollary}
\label{co:Pi(C(Delta))}
Let $\delta$ be a combinatorial gallery of weight $\nu$, as in
(\ref{eq:Gallery}), and equip the set $\Phi_+^\aff(\Delta'_0,\Delta^{}_0)$
with a total order. Then $\pi(C(\delta))$ is the image of the map
$$(a_{j,\beta})\mapsto\prod_{j=0}^p\left(\prod_{\beta\in
\Phi_+^\aff(\Delta'_j,\Delta^{}_j)}x_\beta(a_{j,\beta})\right)[t^\nu]$$
from $\prod_{j=0}^p\mathbb C^{\Phi_+^\aff(\Delta'_j,\Delta^{}_j)}$
to $\mathscr G$.
\end{corollary}

Certainly the notation used in Corollary~\ref{co:Pi(C(Delta))} is more
complicated than really needed. Indeed except perhaps for $j=0$, each set
$\Phi_+^\aff(\Delta'_j,\Delta^{}_j)$ has at most one element. Each inner
product is therefore almost always empty or reduced to one factor.
Keeping this fact in mind may help understand the proofs of
Lemma~\ref{le:ControlCommut} and Proposition~\ref{pr:InclMVCyc} in
Section~\ref{ss:CompTh}.

We now endow $\Gamma(\gamma_\lambda)$ with the structure of a crystal.
To do that, we introduce ``root operators'' $e_\alpha$ and $f_\alpha$
for each simple root $\alpha$ of the root system $\Phi$. These operators
act on $\Gamma(\gamma_\lambda)$ and are defined by the following recipe
(see Section~6 in~\cite{GaussentLittelmann05}).

Let $\delta$ be a combinatorial gallery, as in Equation~(\ref{eq:Gallery}).
We call $m\in\mathbb Z$ the smallest integer such that the hyperplane
$H_{\alpha,m}$ contains a face $\Delta'_j$, where $0\leqslant j\leqslant
p+1$.
\begin{itemize}
\item If $m=0$, then $e_\alpha\delta$ is not defined. Otherwise
we find $k\in\{1,\ldots,p+1\}$ minimal such that $\Delta'_k\subseteq
H_{\alpha,m}$, we find $j\in\{0,\ldots,k-1\}$ maximal such that
$\Delta'_j\subseteq H_{\alpha,m+1}$, and we define the combinatorial
gallery $e_\alpha\delta$ as
\begin{multline*}
\bigl(\{0\}=\Delta'_0\subset\overline{\Delta^{}_0}\supset
\Delta'_1\subset\overline{\Delta^{}_1}\supset\cdots\supset
\Delta'_j\subset\\s_{\alpha,m+1}\bigl(\,\overline{\Delta^{}_j}\,
\bigr)\supset s_{\alpha,m+1}\bigl(\Delta'_{j+1}\bigr)\subset\cdots
\supset s_{\alpha,m+1}\bigl(\Delta'_{k-1}\bigr)\subset
s_{\alpha,m+1}\bigl(\,\overline{\Delta^{}_{k-1}}\,\bigr)\\\supset
\tau_{\alpha^\vee}\bigl(\Delta'_k\bigr)\subset\tau_{\alpha^\vee}
\bigl(\,\overline{\Delta^{}_k}\,\bigr)\supset\cdots\subset
\tau_{\alpha^\vee}\bigl(\,\overline{\Delta^{}_p}\,\bigr)\supset
\tau_{\alpha^\vee}\bigl(\Delta'_{p+1}\bigr)=\{\nu+\alpha^\vee\}\bigr).
\end{multline*}
Thus we reflect all faces between $\Delta'_j$ and $\Delta'_k$ across
the hyperplane $H_{\alpha,m+1}$ and we translate all faces after
$\Delta'_k$ by $\alpha^\vee$. (Note here that $s_{\alpha,m+1}(\Delta'_j)
=\Delta'_j$ and that $s_{\alpha,m+1}(\Delta'_k)=\tau_{\alpha^\vee}
(\Delta'_k)$.)
\item If $m=\langle\alpha,\nu\rangle$, then $f_\alpha\delta$ is not
defined. Otherwise we find $j\in\{0,\ldots,p\}$ maximal such that
$\Delta'_j\subseteq H_{\alpha,m}$, we find $k\in\{j+1,\ldots,p+1\}$
minimal such that $\Delta'_k\subseteq H_{\alpha,m+1}$, and we define
the combinatorial gallery $f_\alpha\delta$ as
\begin{multline*}
\bigl(\{0\}=\Delta'_0\subset\overline{\Delta^{}_0}\supset
\Delta'_1\subset\overline{\Delta^{}_1}\supset\cdots\supset
\Delta'_j\subset\\s_{\alpha,m}\bigl(\,\overline{\Delta^{}_j}\,
\bigr)\supset s_{\alpha,m}\bigl(\Delta'_{j+1}\bigr)\subset\cdots
\supset s_{\alpha,m}\bigl(\Delta'_{k-1}\bigr)\subset
s_{\alpha,m}\bigl(\,\overline{\Delta^{}_{k-1}}\,\bigr)\\\supset
\tau_{-\alpha^\vee}\bigl(\Delta'_k\bigr)\subset\tau_{-\alpha^\vee}
\bigl(\,\overline{\Delta^{}_k}\,\bigr)\supset\cdots\subset
\tau_{-\alpha^\vee}\bigl(\,\overline{\Delta^{}_p}\,\bigr)\supset
\tau_{-\alpha^\vee}\bigl(\Delta'_{p+1}\bigr)=\{\nu-\alpha^\vee\}\bigr).
\end{multline*}
Thus we reflect all faces between $\Delta'_j$ and $\Delta'_k$ across
the hyperplane $H_{\alpha,m}$ and we translate all faces after
$\Delta'_k$ by $-\alpha^\vee$. (Note here that $s_{\alpha,m}(\Delta'_j)
=\Delta'_j$ and that $s_{\alpha,m}(\Delta'_k)=\tau_{-\alpha^\vee}
(\Delta'_k)$.)
\end{itemize}
With the notations above, the maximal integer $n$ such that
$(e_\alpha)^n\delta$ is defined is equal to $-m$, and the maximal
integer $n$ such that $(f_\alpha)^n\delta$ is defined is equal to
$\langle\alpha,\nu\rangle-m$.

The crystal structure on $\Gamma(\gamma_\lambda)$ is then defined
as follows. Given $\delta\in\Gamma(\gamma_\lambda)$, written as
in~(\ref{eq:Gallery}), and $i\in I$, we set
$$\wt(\delta)=\nu,\quad\varepsilon_i(\delta)=-m\quad\text{and}\quad
\varphi_i(\delta)=\langle\alpha_i,\nu\rangle-m,$$
where $\nu$ is the weight of $\delta$ and $m\in\mathbb Z$ is the
smallest integer such that the hyperplane $H_{\alpha_i,m}$ contains
a face $\Delta'_j$, with $0\leqslant j\leqslant p+1$. Finally
$\tilde e_i$ and $\tilde f_i$ are given by the root operators
$e_{\alpha_i}$ and $f_{\alpha_i}$.

Let $\delta$ be a combinatorial gallery, written as in
(\ref{eq:Gallery}). We say that $\delta$ is positively folded if
$$\forall j\in\{1,\ldots,p\},\quad\Delta_{j-1}=\Delta_j\
\Longrightarrow\ \Phi_+^\aff(\Delta'_j,\Delta^{}_j)\neq\varnothing.$$
We define the dimension of $\delta$ as
$$\dim\delta=\sum_{j=0}^p\bigl|\Phi_+^\aff(\Delta'_j,\Delta^{}_j)
\bigr|.$$
(These are Definitions~16 and~17 in \cite{GaussentLittelmann05}.)
Thus for instance the gallery $\gamma_\lambda$ is positively folded
of dimension
\begin{equation}
\label{eq:DimGammaLambda}
\dim\gamma_\lambda=\bigl|\Phi_+\bigr|+p=
\height(\lambda-w_0\lambda)+\dim(P_\lambda/B^+),
\end{equation}
by Equation~(\ref{eq:DimVarBS}). We denote the set of positively
folded combinatorial gallery by $\Gamma^+(\gamma_\lambda)$. Arguing
as in the proof of Proposition~4 in \cite{GaussentLittelmann05}, one
shows that for each $\delta\in\Gamma^+(\gamma_\lambda)$ of weight~$\nu$,
$$\dim\gamma_\lambda-\dim\delta\geqslant\height(\lambda-\nu).$$

We say that a positively folded combinatorial gallery $\delta$ is
an LS gallery if this inequality is in fact an equality. The set of LS
galleries is denoted by $\Gamma^+_\LS(\gamma_\lambda)$. Then Corollary~2
in \cite{GaussentLittelmann05} says that $\Gamma^+_\LS(\gamma_\lambda)$
is a subcrystal of $\Gamma(\gamma_\lambda)$ and that for any gallery
$\delta\in\Gamma^+_\LS(\gamma_\lambda)$, there is a sequence
$(\alpha_1,\ldots,\alpha_t)$ of simple roots such that
$\delta=f_{\alpha_1}\cdots f_{\alpha_t}\gamma_\lambda$. Moreover Lemma~7
and Definition~21 in~\cite{GaussentLittelmann05} say that if $\delta$ is
an LS gallery, written as in~(5.3), if $\alpha$ is a simple root, and if
$m\in\mathbb Z$ is the smallest integer such that the hyperplane
$H_{\alpha,m}$ contains a face $\Delta'_j$, where $0\leqslant j\leqslant
p+1$, then $\delta$ does not cross $H_{\alpha,m}$; this implies that
$\Delta_{j-1}=\Delta_j$ for all $j\in\{1,\ldots,p\}$ such that
$\Delta'_j\subseteq H_{\alpha,m}$.

The following proposition makes the link between LS galleries and MV
cycles; it is equivalent to Corollary~5 in \cite{GaussentLittelmann05}
when $\lambda$ is regular.
\begin{proposition}
\label{pr:LSGalMVCyc}
The map $Z:\delta\mapsto\overline{\pi(C(\delta))}$ is a bijection
from $\Gamma^+_\LS(\gamma_\lambda)$ onto $\mathscr Z(\lambda)$;
it maps a combinatorial gallery of weight $\nu$ to a MV cycle in
$\mathscr Z(\lambda)_\nu$.
\end{proposition}
\begin{proof}
We fix $\nu\in\Lambda$. We denote the set of combinatorial galleries of
weight $\nu$ by $\Gamma(\gamma_\lambda,\nu)$ and we set
$\Gamma^+(\gamma_\lambda,\nu)=\Gamma^+(\gamma_\lambda)\cap
\Gamma(\gamma_\lambda,\nu)$. By construction,
$$\pi^{-1}(S_\nu^+)=\bigsqcup_{\delta\in\Gamma(\gamma_\lambda,\nu)}
C(\delta).$$

We set $\mathring\Sigma=\pi^{-1}\bigl(\mathscr G_\lambda\bigr)$
and $X=\pi^{-1}\bigl(S_\nu^+\cap\mathscr G_\lambda\bigr)$.
Since $S_\nu^+\cap\mathscr G_\lambda$ is of pure dimension
$\height(\nu-w_0\lambda)$, Proposition~\ref{pr:PiAboveOrbit}
and Equation~(\ref{eq:DimGammaLambda}) imply that $X$ is of pure
dimension
$$\height(\nu-w_0\lambda)+\dim(P_\lambda/B^+)=
\dim\gamma_\lambda-\height(\lambda-\nu).$$
Proposition~\ref{pr:PiAboveOrbit} implies also that the map
$Z\mapsto\pi^{-1}(Z)$ is a bijection from the set of irreducible
components of $S_\nu^+\cap\mathscr G_\lambda$ onto the set
of irreducible components of $X$.

By Lemma~11 in \cite{GaussentLittelmann05}, a cell $C(\delta)$ meets
$\mathring\Sigma$ only if $\delta$ is positively folded. Therefore
$$X=\pi^{-1}\bigl(S_\nu^+\bigr)\cap\mathring\Sigma=
\bigsqcup_{\delta\in\Gamma^+(\gamma_\lambda,\nu)}
\bigl(C(\delta)\cap\mathring\Sigma\bigr).$$
Now let $\delta\in\Gamma^+(\gamma_\lambda,\nu)$. Proposition
\ref{pr:C(Delta)} says that the cell $C(\delta)$ is isomorphic to
$\Stab_+(\delta)$, thus is an affine space of dimension $\dim\delta$.
The intersection $C(\delta)\cap\mathring\Sigma$, as a non-empty open
subset of $C(\delta)$, is then irreducible of dimension
$\dim\delta\leqslant\dim\gamma_\lambda-\height(\lambda-\nu)$.
It follows that the irreducible components of $X$ are the
closures in $X$ of the subsets $C(\delta)\cap\mathring\Sigma$,
for $\delta$ running over the set of LS galleries of weight $\nu$.

To conclude the proof, it remains to observe that
$$\overline{\pi\bigl(C(\delta)\cap\mathring\Sigma\bigr)}=
\overline{\pi\bigl(C(\delta)\bigr)}$$
for each $\delta\in\Gamma^+(\gamma_\lambda,\nu)$, since
$C(\delta)\cap\mathring\Sigma$ is dense in $C(\delta)$.
\end{proof}

\subsection{The comparison theorem}
\label{ss:CompTh}
The aim of this section is to show the following property of the map
$Z$ defined in Proposition~\ref{pr:LSGalMVCyc}.
\begin{theorem}
\label{th:CompLSBG}
The bijection $Z:\Gamma^+_\LS(\gamma_\lambda)\to\mathscr Z(\lambda)$
is an isomorphism of crystals.
\end{theorem}

The existence of an isomorphism of crystals from $\mathbf B(\lambda)$
onto $\Gamma^+_\LS(\gamma_\lambda)$ was already known; see for instance
Theorem~2 in~\cite{GaussentLittelmann05} for the case $\lambda$ regular.
The theorem above says that the map $Z^{-1}\circ\Xi(\lambda)$ is actually
such an isomorphism. For its proof, we need two propositions and a lemma.
\begin{proposition}
\label{pr:ImCellRetr}
Let $\delta$ be a combinatorial gallery of weight $\nu$, written as in
(\ref{eq:Gallery}), and let $i\in I$. Call $m$ the smallest integer
such that the hyperplane $H_{\alpha_i,m}$ contains a face $\Delta'_j$
of the gallery, where $0\leqslant j\leqslant p+1$, and set
$\rho=\nu-(\langle\alpha_i,\nu\rangle-m)\,\alpha_i^\vee$. Then
$$r_{\{i\}}(\pi(C(\delta)))=S_{\nu,\{i\}}^+\cap\overline{
S_{\rho,\{i\}}^-}\quad\text{and}\quad s_i\mu_+\Bigl(\overline{s_i}^{-1}
\;\overline{\pi(C(\delta))}\Bigr)=\rho.$$
\end{proposition}
\begin{proof}
We collect in a set $J$ the indices $j\in\{0,\ldots,p\}$ such that
$\Phi_+^\aff(\Delta'_j,\Delta^{}_j)$ contains an affine root of the
form $(\alpha_i,n)$, with $n\in\mathbb Z$. For each $j\in J$, there
is a unique integer, say $n_j$, so that $(\alpha_i,n_j)\in\Phi_+^\aff
(\Delta'_j,\Delta^{}_j)$. (Thus $n_j=f_{\Delta'_j}(\alpha_i)$ in
the notation of Section~\ref{ss:CoxCo&BTBui}.)

All these integers $n_j$ are larger or equal than $m$. We claim that
\begin{equation}
\label{eq:ImCellRetr}
\{m,m+1,m+2,\ldots\}\supseteq\{n_j\mid j\in J\}\supseteq
\{m,m+1,\ldots,\langle\alpha_i,\nu\rangle-1\}.
\end{equation}
Consider indeed an integer $n$ in the right-hand side above. Since
the gallery $\delta$ must go from the wall $H_{\alpha_i,m}$ to the
point $\nu$, it must cross the wall $H_{\alpha_i,n}$. More exactly,
there is an index $j\in\{0,\ldots,p\}$ such that $\Delta'_j\subseteq
H_{\alpha_i,n}$ and $\Delta^{}_j\not\subseteq H^-_{\alpha_i,n}$;
this implies that $(\alpha_i,n)\in\Phi_+^\aff(\Delta'_j,\Delta^{}_j)$,
and thus that $j\in J$ and $n=n_j$.

We apply now the parabolic retraction $r_{\{i\}}$ to the expression
given in Corollary~\ref{co:Pi(C(Delta))}. Equation~(\ref{eq:CommutInRetr})
allows us to remove all factors in the product that belong to the
unipotent radical of $P_{\{i\}}(\mathscr K)$. We deduce that
$r_{\{i\}}(\pi(C(\delta)))$ is the image of the map
$$(a_j)\mapsto\prod_{j\in J}x_{\alpha_i,n_j}(a_j)[t^\nu]$$
from $\mathbb C^J$ to $\mathscr M_{\{i\}}$. Using~(\ref{eq:ImCellRetr})
and the fact that $[t^\nu]$ is fixed by all subgroups
$x_{\alpha_i,n}(\mathbb C)$ with $n\geqslant\langle\alpha_i,\nu\rangle$,
we then get
$$r_{\{i\}}(\pi(C(\delta)))=\bigl\{x_{\alpha_i}\bigl(pt^{\langle\alpha,
\nu\rangle}\bigr)[t^\nu]\bigm|p\in\mathbb C[t^{-1}]_{\langle\alpha,
\nu\rangle-m}^+\bigr\}.$$
From there, the proposition follows easily using Proposition
\ref{pr:GpsOfRankOne} (with $+$ and $-$ exchanged) and
Lemma~\ref{le:RetractStrata}.
\end{proof}

For a combinatorial gallery $\delta$, written as in Equation
(\ref{eq:Gallery}), and an integer $k\in\{0,\ldots,p+1\}$, we set
\begin{align*}
\Stab_+(\delta)_{\geqslant k}&=\Stab_+(\Delta'_k,\Delta^{}_k)\times
\Stab_+(\Delta'_{k+1},\Delta^{}_{k+1})\times\cdots\times
\Stab_+(\Delta'_p,\Delta^{}_p),\\[5pt]
\pi(C(\delta))_{\geqslant k}&=\bigl\{v_kv_{k+1}\cdots v_p[t^\nu]\bigm|
(v_k,v_{k+1},\ldots,v_p)\in\Stab_+(\delta)_{\geqslant k}\bigr\}.
\end{align*}
\begin{lemma}
\label{le:ControlCommut}
Let $\delta$ be a combinatorial gallery, as in Equation
(\ref{eq:Gallery}), and let $k\in\{0,\ldots,p+1\}$.
\begin{enumerate}
\item\label{it:LeCCa}
Let $u\in\Stab_+(\Delta'_k)$. Then the left action of $u$ on
$\mathscr G$ leaves $\pi(C(\delta))_{\geqslant k}$ stable. More
precisely, for each $(v_k,\ldots,v_p)\in\Stab_+(\delta)_{\geqslant k}$,
there exists $(v'_k,\ldots,v'_p)\in\Stab_+(\delta)_{\geqslant k}$
such that $v'_k\cdots v'_p[t^\nu]=uv_k\cdots v_p[t^\nu]$ and
$$\forall j\in\{k+1,\ldots,p\},\quad\Delta^{}_{j-1}=\Delta^{}_j
\ \Longrightarrow\ v^{}_j=v'_j;$$
moreover if $k>0$ and $u\in\Stab_+(\Delta^{}_k)$, then one can manage
so that $v^{}_k=v'_k$.
\item\label{it:LeCCb}
Assume that $k>0$, let $p\in\mathscr O^\times$ and let $\mu\in\Lambda$.
Then the left action of $p^\mu$ on $\mathscr G$ leaves
$\pi(C(\delta))_{\geqslant k}$ stable. Suppose moreover that
$p\in1+t\mathscr O$ and let $(v_k,\ldots,v_p)\in
\Stab_+(\delta)_{\geqslant k}$. Then there exists
$(v'_k,\ldots,v'_p)\in\Stab_+(\delta)_{\geqslant k}$ such that
$v'_k\cdots v'_p[t^\nu]=p^\mu v^{}_k\cdots v^{}_p[t^\nu]$ and
$$\forall j\in\{k,\ldots,p\},\quad\Delta^{}_{j-1}=\Delta^{}_j
\ \Longrightarrow\ v^{}_j=v'_j.$$
\item\label{it:LeCCc}
Assume that $k>0$ and that $\delta$ is an LS gallery. Let
$(v_k,\ldots,v_p)\in\Stab_+(\delta)_{\geqslant k}$, let $\alpha$ be
a simple root of the root system $\Phi$, and let $c\in\mathbb
C^\times$. Call $m$ the smallest integer such that the
hyperplane $H_{\alpha,m}$ contains a face $\Delta'_j$, where
$0\leqslant j\leqslant p+1$, form the list $(k_1,k_2,\ldots,k_r)$ in
increasing order of all indices $l\in\{k,\ldots,p\}$ such that
$\Phi_+^\aff(\Delta'_l,\Delta^{}_l)=\{(\alpha,m)\}$, and find the
complex numbers $c_1$, $c_2$, \dots, $c_r$ such that $v_{k_s}=
x_{\alpha,m}(c_s)$. Assume that $c+c_1+c_2+\cdots+c_s\neq0$ for each
$s\in\{1,\ldots,r\}$. Then $x_{-\alpha,-m}(1/c)v^{}_k\cdots
v^{}_p[t^\nu]$ belongs to $\pi(C(\delta))_{\geqslant k}$.
\end{enumerate}
\end{lemma}
\begin{proof}
The proof of these three assertions proceeds by decreasing induction
on $k$. For $k=p+1$, all of them hold: indeed the element $u$ in
Assertion~\ref{it:LeCCa}, the element $p^\mu$ in
Assertion~\ref{it:LeCCb} and the element $x_{-\alpha,-m}(1/c)$ in
Assertion~\ref{it:LeCCc} fix the point $[t^\nu]$.

Now assume that $k\leqslant p$ and that the result holds for $k+1$.
If $\Phi_+^\aff(\Delta'_k,\Delta^{}_k)$ is empty, then
$\Stab_+(\Delta'_k,\Delta^{}_k)=\{1\}$. Assertions~\ref{it:LeCCa},
\ref{it:LeCCb} and \ref{it:LeCCc} follow then immediately from the
inductive assumption, after one has observed that the element $u$
in Assertion~\ref{it:LeCCa} belongs by assumption to
$\Stab_+(\Delta'_k)$ and that $\Stab_+(\Delta'_k)=\Stab_+(\Delta_k)
\subseteq\Stab_+(\Delta'_{k+1})$. In the rest of the proof, we assume
that $\Phi_+^\aff(\Delta'_k,\Delta^{}_k)$ is not empty. Let
$(v_k,\ldots,v_p)\in\Stab_+(\delta)_{\geqslant k}$. Except in the
case $k=0$ (dealt with only in Assertion~\ref{it:LeCCa}),
$\Phi_+^\aff(\Delta'_k,\Delta^{}_k)$ has a unique element,
say $(\zeta,n)$ with $\zeta\in\Phi_+$, and there exists $b\in\mathbb C$
such that $v_k=x_{\zeta,n}(b)$.

Consider first Assertion~\ref{it:LeCCa}. The element $uv_k$ belongs to
$\Stab_+(\Delta'_k)$. By Proposition~\ref{pr:DescStab+}~\ref{it:PrDS+b},
there exists $v'_k\in\Stab_+(\Delta'_k,\Delta^{}_k)$ and
$u'\in\Stab_+(\Delta^{}_k)$ such that $uv^{}_k=v'_ku'$. The
inductive assumption applied to $u'$ and $(v_{k+1},\ldots,v_p)
\in\Stab_+(\delta)_{\geqslant k+1}$ asserts the existence of
$(v'_{k+1},\ldots,v'_p)\in\Stab_+(\delta)_{\geqslant k+1}$ such that
$u'v_{k+1}\cdots v_p[t^\nu]=v'_{k+1}\cdots v'_p[t^\nu]$, with the
further property that $v^{}_j=v'_j$ for all $j>k$ verifying
$\Delta^{}_{j-1}=\Delta^{}_j$. Certainly then $uv^{}_kv^{}_{k+1}\cdots
v^{}_p[t^\nu]=v'_kv'_{k+1}\cdots v'_p[t^\nu]$. Now assume that $k>0$ and
that $u\in\Stab_+(\Delta^{}_k)$. By Proposition~\ref{pr:DescStab+}
\ref{it:PrDS+a}, we may write $u$ as a product of elements of
the form $x_{\beta,n}(q)$ with $q\in\mathscr O$ and $(\beta,n)\in
\Phi_+\times\mathbb Z$ such that $\Delta^{}_k\subseteq H^-_{\beta,n}$.
Lemma~\ref{le:CommStab+} now implies that $uv^{}_k\in v^{}_k
\Stab_+(\Delta^{}_k)$, which establishes $v'_k=v^{}_k$. This shows
that Assertion~\ref{it:LeCCa} holds at $k$.

Consider now Assertion~\ref{it:LeCCb}. Let $a\in\mathbb C^\times$ be
the constant term coefficient of $p$ and set $q=\bigl(p^{\langle\zeta,
\mu\rangle}-a^{\langle\zeta,\mu\rangle}\bigr)/t$. Then
$$p^\mu v^{}_k=x_{\zeta,n}\bigl(bp^{\langle\zeta,\mu\rangle}\bigr)
p^\mu=x_{\zeta,n}(b')u'p^\mu=v'_ku'p^\mu,$$
where $b'=ba^{\langle\zeta,\mu\rangle}$, $u'=x_{\zeta,n+1}(bq)$ and
$v'_k=x_{\zeta,n}(b')$. Observing that $u'\in\Stab_+(\Delta^{}_k)$
and using the inductive assumption and Assertion~\ref{it:LeCCa}, we find
$(v'_{k+1},\ldots,v'_p)\in\Stab_+(\delta)_{\geqslant k+1}$ such that
$u'p^\mu v^{}_{k+1}\cdots v^{}_p[t^\nu]=v'_{k+1}\cdots v'_p[t^\nu]$; in
the case $a=1$, we may even demand that $v^{}_j=v'_j$ for all $j>k$
verifying $\Delta^{}_{j-1}=\Delta^{}_j$. Then $p^\mu v^{}_kv^{}_{k+1}
\cdots v^{}_p[t^\nu]=v'_kv'_{k+1}\cdots v'_p[t^\nu]$, which shows that
Assertion~\ref{it:LeCCb} holds at $k$.

It remains to prove Assertion~\ref{it:LeCCc}. We distinguish
several cases.

Suppose first that $\zeta\neq\alpha$. By Lemma~\ref{le:CommStab+},
the element
$$u=x_{-\alpha,-m}(-1/c)\;(v^{}_k)^{-1}\;x_{-\alpha,-m}(1/c)\;v^{}_k$$
belongs to $\Stab_+(\Delta^{}_k)$. Using Assertion~\ref{it:LeCCa},
we find $(v'_{k+1},\ldots,v'_p)\in\Stab_+(\delta)_{\geqslant k+1}$
such that $uv^{}_{k+1}\cdots v^{}_p[t^\nu]=v'_{k+1}\cdots v'_p[t^\nu]$.
Moreover, since $\delta$ is an LS gallery, we know that
$\Delta_{k_s-1}=\Delta_{k_s}$ for each $s\in\{1,\ldots,r\}$,
and we may thus demand that $v'_{k_s}=v^{}_{k_s}=x_{\alpha,m}(c_s)$.
Applying the inductive assumption, we find a tuple
$(v''_{k+1},\ldots,v''_p)\in\Stab_+(\delta)_{\geqslant k+1}$ such that
$x_{-\alpha,-m}(1/c)\,v'_{k+1}\cdots v'_p[t^\nu]=v''_{k+1}\cdots
v''_p[t^\nu]$. Then
$$x_{-\alpha,-m}(1/c)\,v^{}_kv^{}_{k+1}\cdots v^{}_p[t^\nu]=v_k
v''_{k+1}\cdots v''_p[t^\nu],$$
which establishes that Assertion~\ref{it:LeCCc} holds at $k$ in this
first case.

The second case is when $\zeta=\alpha$ but $n\neq m$. Then
$n>m$, by the minimality of $m$. Let $p$ be the square root in
$1+t\mathscr O$ of $1+t^{n-m}b/c$. Equation~(\ref{eq:EqCommut3})
implies that
\begin{align*}
x_{-\alpha,-m}(1/c)v^{}_k&=x_{-\alpha}(1/ct^m)x_\alpha(bt^n)\\[4pt]
&=p^{-\alpha^\vee}x_\alpha(bt^n)x_{-\alpha}(1/ct^m)p^{-\alpha^\vee}\\[4pt]
&=p^{-\alpha^\vee}v^{}_kx_{-\alpha,-m}(1/c)p^{-\alpha^\vee}.
\end{align*}
Assertion~\ref{it:LeCCb} allows us to find $(v'_{k+1},\ldots,v'_p)\in
\Stab_+(\delta)_{\geqslant k+1}$ such that $p^{-\alpha^\vee}
v_{k+1}\cdots v_p[t^\nu]=v'_{k+1}\cdots v'_p[t^\nu]$, with the
further property that $v'_{k_s}=v^{}_{k_s}=x_{\alpha,m}(c_s)$
for each $s\in\{1,\ldots,r\}$. We apply then the inductive assumption
and find $(v''_{k+1},\ldots,v''_p)\in\Stab_+(\delta)_{\geqslant k+1}$
such that $x_{-\alpha,-m}(1/c)\,v'_{k+1}\cdots v'_p[t^\nu]=v''_{k+1}\cdots
v''_p[t^\nu]$. Then
$$x_{-\alpha,-m}(1/c)\,v^{}_kv^{}_{k+1}\cdots v^{}_p[t^\nu]
=p^{-\alpha^\vee}v^{}_kv''_{k+1}\cdots v''_p[t^\nu],$$
and a final application of Assertion~\ref{it:LeCCb} concludes the
proof of Assertion~\ref{it:LeCCc} at $k$ in this second case.

The last case is $(\zeta,n)=(\alpha,m)$. In this case, $k_1=k$
and $b=c_{k_1}$. The assumptions of the lemma imply that $b+c\neq0$.
Equation~(\ref{eq:EqCommut3}) says then that
$$x_{-\alpha,-m}(1/c)v^{}_k=x_{\alpha,m}(bc/(b+c))
(1+b/c)^{-\alpha^\vee}x_{-\alpha,-m}(1/(b+c)).$$
Applying the inductive assumption, we find $(v'_{k+1},\ldots,v'_p)
\in\Stab_+(\delta)_{\geqslant k+1}$ such that
$$x_{-\alpha,-m}(1/(b+c))\,v^{}_{k+1}\cdots v^{}_p[t^\nu]=v'_{k+1}\cdots
v'_p[t^\nu].$$
Using now Assertion~\ref{it:LeCCb}, we see that
$$x_{-\alpha,-m}(1/c)\,v^{}_kv^{}_{k+1}\cdots v^{}_p[t^\nu]=x_{\alpha,m}
(bc/(b+c))\,(1+b/c)^{-\alpha^\vee}v'_{k+1}\cdots v'_p[t^\nu]$$
belongs to $\pi(C(\delta))_{\geqslant k}$. This concludes the proof
of Assertion~\ref{it:LeCCc} at $k$.
\end{proof}

At the end of their paper~\cite{GaussentLittelmann05}, Gaussent and
Littelmann describe several cases where the crystal structure on
$\Gamma^+_\LS(\gamma_\lambda)$ controls inclusions between MV cycles.
The next proposition presents a general result.
\begin{proposition}
\label{pr:InclMVCyc}
Let $\delta$ be an LS gallery and let $\alpha$ be a simple root of
the system $\Phi$. If the gallery $e_\alpha\delta$ is defined, then
$Z(\delta)\subseteq Z(e_\alpha\delta)$.
\end{proposition}
\begin{proof}
We represent $\delta$ as in~(\ref{eq:Gallery}). We assume that
$e_\alpha\delta$ is defined and we let $m\in\mathbb Z$ and $0\leqslant
j<k\leqslant p+1$ be as in the definition of $e_\alpha\delta$. We
call $(k=k_0,k_1,\ldots,k_r)$ the list in increasing order of all indices
$l\in\{1,\ldots,p\}$ such that $\Phi_+^\aff(\Delta'_l,\Delta^{}_l)=
\bigl\{(\alpha,m)\bigr\}$. Finally we equip $\Phi_+^\aff(\Delta'_0,
\Delta^{}_0)$ with a total order.

Let $(a_{l,\beta})\in\prod_{l=0}^p\mathbb C^{\Phi_+^\aff
(\Delta'_l,\Delta^{}_l)}$ be a family of complex numbers such that
$$a_{k_0,(\alpha,m)}+a_{k_1,(\alpha,m)}+\cdots+a_{k_s,(\alpha,m)}
\neq0$$
for each $s\in\{0,1,\ldots,r\}$ and set
\begin{align*}
v_l&=\prod_{\beta\in\Phi_+^\aff(\Delta'_l,\Delta^{}_l)}
x_\beta(a_{l,\beta})\quad\text{for each $l\in\{0,1,\ldots,p\}$},\\[2pt]
A&=\prod_{l=0}^{j-1}v_l\qquad\text{and}\qquad B=\prod_{l=j}^pv_l.
\end{align*}
By Corollary~\ref{co:Pi(C(Delta))}, the element $AB[t^\nu]$ describes
a dense subset of $Z(\delta)$ when the parameters $a_{l,\beta}$ vary.
To establish the proposition, it therefore suffices to show that
$AB[t^\nu]$ belongs to $Z(e_\alpha\delta)$. What we will now show is
more precise:\\[5pt]
\textit{For any non-zero complex number $h$, the element
$Ax_{-\alpha,-m-1}(h)B[t^\nu]$ belongs to $\pi(C(e_\alpha\delta))$.}
\vspace{5pt}

We first observe that
$x_{\alpha,m+1}(1/h)\in\Stab_+(\Delta'_j)$, for $\Delta'_j\subseteq
H_{\alpha,m+1}$. Using Lemma~\ref{le:ControlCommut}~\ref{it:LeCCa},
we find $(v'_j,v'_{j+1},\ldots,v'_p)\in\Stab_+(\delta)_{\geqslant j}$
such that
$$x_{\alpha,m+1}(1/h)B[t^\nu]=v'_jv'_{j+1}\cdots v'_p[t^\nu].$$
We may moreover demand that $v'_{k_s}=v^{}_{k_s}=x_{\alpha,m}
(a_{k_s,(\alpha,m)})$ for all $s\in\{0,1,\ldots,r\}$, for
$\Delta_{k_s-1}=\Delta_{k_s}$. We set
$$C=\prod_{l=j}^{k-1}v'_l\qquad\text{and}\qquad D=\prod_{l=k+1}^pv'_l,$$
and then $B[t^\nu]=x_{\alpha,m+1}(-1/h)Cv'_kD[t^\nu]$. Using
Lemma~\ref{le:ControlCommut}~\ref{it:LeCCc}, we now find
$(v''_{k+1},v''_{k+2},\linebreak\ldots,v''_p)\in\Stab_+(\delta
)_{\geqslant k+1}$ such that
$$x_{-\alpha,-m}\bigl(1/a_{k,(\alpha,m)}\bigr)D[t^\nu]=
v''_{k+1}v''_{k+2}\cdots v''_p[t^\nu].$$
We finally set
\begin{align*}
E&=x_{\alpha,m}\bigl(a_{k,(\alpha,m)}\bigr)
x_{-\alpha,-m}\bigl(-1/a_{k,(\alpha,m)}\bigr)
x_{\alpha,m}\bigl(a_{k,(\alpha,m)}\bigr),\\[5pt]
F&=x_{\alpha,m}\bigl(-a_{k,(\alpha,m)}\bigr)\prod_{l=k+1}^pv''_l,\\[5pt]
K&=x_{-\alpha,-m-1}(h)x_{\alpha,m+1}(-1/h).
\end{align*}
Then $Ax_{-\alpha,-m-1}(h)B[t^\nu]=AKCEF[t^\nu]$.

We now observe that
$$\Phi_+^\aff\bigl(s_{\alpha,m+1}(\Delta'_l),s_{\alpha,m+1}
(\Delta^{}_l)\bigr)=\begin{cases}
\bigl\{(\alpha,m+1)\bigr\}\sqcup s_{\alpha,m+1}
\bigl(\Phi_+^\aff(\Delta'_j,\Delta^{}_j)\bigr)
&\text{if $l=j$,}\\[5pt]
s_{\alpha,m+1}\bigl(\Phi_+^\aff(\Delta'_l,\Delta^{}_l)\bigr)
&\text{if $j<l<k$,}
\end{cases}$$
and that
$$\Phi_+^\aff\bigl(\tau_{\alpha^\vee}(\Delta'_l),\tau_{\alpha^\vee}
(\Delta^{}_l)\bigr)=\tau_{\alpha^\vee}\bigl(\Phi_+^\aff
(\Delta'_l,\Delta^{}_l)\bigr)\qquad\text{if $l\geqslant k$.}$$
These equalities, the definition of $e_\alpha\delta$, Equation
(\ref{eq:EqCommut6}) and Proposition~\ref{pr:DescStab+}~\ref{it:PrDS+b}
imply that the sequence
\begin{multline*}
\Bigl(v^{}_0,\;\ldots,\;v^{}_{j-1},\;
x_{\alpha,m+1}(h)\bigl(t^{(m+1)\alpha^\vee}\overline{s_\alpha}\,\bigr)
v'_j\bigl(t^{(m+1)\alpha^\vee}\overline{s_\alpha}\,\bigr)^{-1},\\[4pt]
\bigl(t^{(m+1)\alpha^\vee}\overline{s_\alpha}\,\bigr)v'_{j+1}
\bigl(t^{(m+1)\alpha^\vee}\overline{s_\alpha}\,\bigr)^{-1},\;
\ldots,\,\bigl(t^{(m+1)\alpha^\vee}\overline{s_\alpha}\,\bigr)
v'_{k-1}\bigl(t^{(m+1)\alpha^\vee}\overline{s_\alpha}\,\bigr)^{-1},\\[4pt]
t^{\alpha^\vee}x_{\alpha,m}\bigl(-a_{k,(\alpha,m)}\bigr)t^{-\alpha^\vee},\;
t^{\alpha^\vee}v''_{k+1}t^{-\alpha^\vee},\;\ldots,\;
t^{\alpha^\vee}v''_pt^{-\alpha^\vee}\Bigr)
\end{multline*}
belongs to $\Stab_+(e_\alpha\delta)$. Proposition~\ref{pr:C(Delta)},
Equation~(\ref{eq:RecTarget}) and the definition of the map $\pi$ then
say that
$$A\;x_{\alpha,m+1}(h)\;\bigl(t^{(m+1)\alpha^\vee}\overline{s_\alpha}\,
\bigr)\;C\;\bigl(t^{(m+1)\alpha^\vee}\overline{s_\alpha}\,\bigr)^{-1}\;
t^{\alpha^\vee}\,F\,[t^\nu]$$
belongs to $\pi(C(e_\alpha\delta))$. An appropriate application of
Lemma~\ref{le:ControlCommut}~\ref{it:LeCCb} shows that the element
obtained by inserting extra factors $(-h)^{-\alpha^\vee}$ and
$\bigl(-a_{k,(\alpha,m)}\bigr)^{-\alpha^\vee}$ in this expression,
respectively after $A$ and before $t^{\alpha^\vee}$, also belongs to
$\pi(C(e_\alpha\delta))$. Now Equation~(\ref{eq:EqCommut4}) allows
to rewrite
$$K=(-h)^{-\alpha^\vee}\,x_{\alpha,m+1}(h)\;\bigl(t^{(m+1)\alpha^\vee}
\overline{s_\alpha}\,\bigr)$$
and
$$E=\bigl(t^{(m+1)\alpha^\vee}\overline{s_\alpha}\,\bigr)^{-1}\;
\bigl(-a_{k,(\alpha,m)}\bigr)^{-\alpha^\vee}\;t^{\alpha^\vee},$$
and we conclude that $AKCEF[t^\nu]=Ax_{-\alpha,-m-1}(h)B[t^\nu]$
belongs to $\pi(C(e_\alpha\delta))$, as announced.
\end{proof}

\trivlist
\item[\hskip\labelsep{\itshape Proof of Theorem~\ref{th:CompLSBG}.}]
\upshape
Obviously $Z$ preserves the weight. Comparing Proposition
\ref{pr:ImCellRetr} with Equation~(\ref{eq:DefPhiBFG}), we see that $Z$
is compatible with the structure maps $\varphi_i$. The axioms of a
crystal imply then that $Z$ is compatible with the structure maps
$\varepsilon_i$. Now let $\delta$ be an LS gallery of weight $\nu$,
let $i\in I$, and assume that the LS gallery $e_{\alpha_i}\delta$ is
defined. Then the two MV cycles $Z(\delta)$ and $Z(e_{\alpha_i}\delta)$
satisfy the four conditions of Proposition~\ref{pr:CritOpCrys}. Indeed
the first and the third conditions follow immediately from the fact
that $Z(\delta)\in\mathscr Z(\lambda)_\nu$ and $Z(e_{\alpha_i}\delta)
\in\mathscr Z(\lambda)_{\nu+\alpha_i^\vee}$; the second condition comes
from Proposition~\ref{pr:ImCellRetr} and from the second assertion of
Lemma~6~(iii) in~\cite{GaussentLittelmann05}; the fourth condition comes
from Proposition~\ref{pr:InclMVCyc}. Therefore $Z(e_{\alpha_i}\delta)=\tilde
e_iZ(\delta)$; in other words, $Z$ intertwines the action of the root
operators on $\Gamma^+_\LS(\gamma_\lambda)$ with the action of Braverman
and Gaitsgory's crystal operators on $\mathscr Z(\lambda)$. This
concludes the proof that $Z$ is a morphism of crystals. Since $Z$ is
bijective and both crystals $\Gamma^+_\LS(\gamma_\lambda)$ and
$\mathscr Z(\lambda)$ are normal, $Z$ is an isomorphism.
\nobreak\noindent$\square$
\endtrivlist

\bigskip
\noindent\parbox[t]{214pt}{\noindent
Pierre Baumann\\
Institut de Recherche Math\'ematique Avanc\'ee\\
Universit\'e Louis Pasteur et CNRS\\
7, rue Ren\'e Descartes\\
67084 Strasbourg Cedex\\
France\\[2pt]
E-mail: \texttt{baumann@math.u-strasbg.fr}}
\bigskip

\noindent\parbox[t]{226pt}{\noindent
St\'ephane Gaussent\\
Institut \'Elie Cartan\\
Unit\'e Mixte de Recherche 7502\\
Nancy-Universit\'e, CNRS, INRIA\\
Boulevard des Aiguillettes\\
B.P.~239\\
54506 Vand\oe uvre-l\`es-Nancy Cedex\\
France\\[2pt]
E-mail: \texttt{Stephane.Gaussent@iecn.u-nancy.fr}}
\end{document}